\documentclass[11pt]{article}
\usepackage{amssymb}
\usepackage{amsmath}
\usepackage{amsthm}
\newcommand{\cB}{{\cal B}}
\newcommand{\cC}{{\cal C}}

\newcommand{\cH}{{\cal H}}
\newcommand{\cE}{{\cal E}}
\newcommand{\cG}{{\cal G}}

\newcommand{\cJ}{{\cal J}}
\newcommand{\cO}{{\cal O}}
\newcommand{\cL}{{\cal L}}

\newcommand{\cF}{{\cal F}}

\newcommand{\cT}{{\cal T}}
\newcommand{\cU}{{\cal U}}
\newcommand{\cV}{{\cal V}}
\newcommand{\cW}{{\cal W}}
\newcommand{\cX}{{\cal X}}
\newcommand{\cY}{{\cal Y}}
\newcommand{\cZ}{{\cal Z}}

\renewcommand{\AA}{{\mathbb A}}

\newcommand{\ZZ}{{\mathbb Z}}
\newcommand{\QQ}{{\mathbb Q}}

\newcommand{\PP}{{\mathbb P}}
\newcommand{\OO}{{\mathbb O}}
\newcommand{\SSS}{{\mathbb S}}
\newcommand{\EE}{{\mathbb E}}
\newcommand{\FF}{{\mathbb F}}
\newcommand{\WW}{{\mathbb W}}

\newcommand{\gt}{\mathfrak{t}}

\newcommand{\gs}{\mathfrak{s}}

\newcommand{\on}{\operatorname}

\newcommand{\mult}{\on{mult}}

\newcommand{\Qlb}{\mathbb{\bar Q}_\ell}
\newcommand{\Gm}{\mathbb{G}_m}
\newcommand{\Ga}{\mathbb{G}_a}
\newcommand{\A}{\mathbb{A}}

\newcommand{\toup}[1]{\stackrel{#1}{\to}}
\newcommand{\hook}[1]{\stackrel{#1}{\hookrightarrow}}

\newcommand{\getsup}[1]{\stackrel{#1}{\gets}}
\newcommand{\Sp}{\on{\mathbb{S}p}}

\newcommand{\Hom}{\on{Hom}}

\newcommand{\Sym}{\on{Sym}}
\newcommand{\SO}{\on{S\mathbb{O}}}

\newcommand{\Aut}{\on{Aut}}

\newcommand{\RG}{\on{R\Gamma}}

\newcommand{\Spec}{\on{Spec}}

\newcommand{\GL}{\on{GL}}

\newcommand{\Fr}{{\on{Fr}}}

\newcommand{\pr}{\on{pr}}
\newcommand{\id}{\on{id}}
\newcommand{\tr}{\on{tr}}
\newcommand{\QED}{$\square$} 
  
\newcommand{\iso}{{\widetilde\to}}

\newcommand{\comp}{\circ}
\newcommand{\Four}{\on{Four}}
\renewcommand{\H}{{\on{H}}}   

\newcommand{\DD}{\mathbb{D}}  
\newcommand{\D}{\on{D}}       

\newcommand{\ov}[1]{\overline{#1}}
\newcommand{\select}[1]{{\it{#1}}}

\renewcommand{\P}{{\on{P}}}

\newcommand{\<}{\langle}
\renewcommand{\>}{\rangle}

\newcommand{\ev}{\mathit{ev}}

\newcommand{\ttimes}{\tilde\times}

\newcommand{\act}{\on{act}}
\newcommand{\dimrel}{\on{dim.rel}}

\renewcommand{\Im}{\on{Im}}

\newcommand{\Cov}{\on{Cov}}
\newcommand{\sign}{\on{sign}}
\newcommand{\SL}{\on{SL}}

\newcommand{\diag}{\on{diag}}

\newcommand{\Disc}{\on{Disc}}
\newcommand{\Mp}{\on{Mp}}
\newcommand{\ASp}{\on{A\mathbb{S}p}}
\newcommand{\AMp}{\on{AMp}}
\newcommand{\Mat}{\on{Mat}}

\newcommand{\Ad}{\on{Ad}}
\newcommand{\disc}{\on{disc}}

\newtheorem{Lm}{Lemma}
\newtheorem{Th}{Theorem}
\newtheorem{Pp}{Proposition}
\newtheorem{Cor}{Corollary}

\theoremstyle{remark}
\newtheorem{Rem}{Remark}

\theoremstyle{definition}
\newtheorem{Def}{Definition}

\newenvironment{Prf}{\par\noindent {\it Proof }}{\QED}
\newcommand{\Step}[1]{\par\noindent{\bf Step {#1}}.}

\begin{document}
\author{Alain Genestier and Sergey Lysenko}
\title{Geometric Weil representation in characteristic two}
\date{}
\maketitle
\begin{abstract}
\noindent{\scshape Abstract}\hskip 0.8 em Let $k$ be an algebraically closed field of characteristic two. Let $R$ be the ring of Witt vectors of length two over $k$. We construct a group stack $\hat G$ over $k$, the metaplectic extension of the Greenberg realization of $\Sp_{2n}(R)$. We also construct a geometric analog of the Weil representation of $\hat G$, this is a triangulated category on which $\hat G$ acts by functors. This triangulated category and the action are geometric in a suitable sense.
\end{abstract} 

\medskip

{\centerline{\scshape 1. Introduction}}

\bigskip\noindent
1.1 Apparently, a version of the Weil representation in characteristic two first appeared in 1958 paper by D. A. Suprunenko (\cite{Su}, Theorem~11) (before the celebrated 1964 paper by A. Weil \cite{W}). This representation and its character were also studied in \cite{Is, HZ}. Being inspired mostly by \cite{GH} and \cite{L2}, in this paper we propose a geometric analog of this representation. 

  Let $k$ be a finite field of characteristic two. Let $R$ be the ring of Witt vectors of length 2 over $k$. Given a free $R$-module $\tilde V$ of rank $2n$ with symplectic form $\tilde\omega: \tilde V\times\tilde V\to R$, set $V=\tilde V\otimes_R k$. Write $\Sp(\tilde V)$ for the group of isometries of the form $\tilde\omega$. Pick a bilinear form $\tilde\beta: \tilde V\times\tilde V\to R$ such that $\tilde\beta(\tilde x, \tilde y)-\tilde\beta(\tilde y, \tilde x)=\tilde\omega(\tilde x, \tilde y)$ for all $\tilde x, \tilde y\in \tilde V$. Let $\beta: V\times V\to 2R\subset R$ be the map $(x,y)\mapsto 2\tilde\beta(\tilde x,\tilde y)$ for any $\tilde x,\tilde y\in\tilde V$ over $x,y\in V$. It gives rise to the Heisenberg group $H(V)=V\times R$ with operation
$$ 
(v_1,z_1)(v_2,z_2)=(v_1+v_2, z_1+z_2+\beta(v_1,v_2))
$$
The reason for using $R$ instead of $2R$ in the definition of $H(V)$ is that in this way it acquires a larger group of automorphisms acting trivially on the center. The group $\Sp(\tilde V)$ maps naturally to this group. 

  Fix a prime $\ell\ne 2$ and a faithful character $\psi: \ZZ/4\ZZ\to \Qlb^*$. A version of the Stone-von Neumann theorem holds in this setting giving rise to the metaplectic extension of $\Sp(\tilde V)$ and its Weil representation $\cH_{\psi}$ (cf. Section~2 for details). According to \cite{GH}, it can be seen as a group $\Mp(\tilde V)$ that fits into an exact sequence
\begin{equation}
\label{ext_introd_one}
1\to \ZZ/2\ZZ\to \Mp(\tilde V)\to \Sp(\tilde V)\to 1
\end{equation}
 
  In the geometric setting, assume $k$ to be an algebraically closed field of characteristic two. We propose geometric analogs of $\Mp(\tilde V)$ and $\cH_{\psi}$. 
Let $\tilde V$ be a free $R$-module of rank $2n$ with a symplectic form. Write $G$ for the Greenberg realization of the $R$-scheme $\Sp(\tilde V)$. View $H(V)$ as a group scheme over $k$, an extension of $V$ by the 
Greenberg realization of $R$.

 Though we mostly follow the strategy of \cite{L2}, there are new difficulties and phenomena in characteristic two. To the difference with the case of other characteristics, the metaplectic extension (\ref{ext_introd_one}) is nontrivial. The geometric analog of (\ref{ext_introd_one}) is an algebraic group stack $\hat G$ over $k$ that fits into an exact sequence 
\begin{equation}
\label{seq_first_for_hat_G}
1\to B(\ZZ/4\ZZ)\to \hat G\to G\to 1
\end{equation}
of group stacks over $k$. Here for an algebraic group $H$ over $k$ we write $B(H)$ for the classifying stack of $H$ over $k$. Actually, from our Remark~\ref{Rem_ELag_so_on} ii) it follows that there is a group stack $\hat G_b$ over $k$ included into an exact sequence $1\to B(\ZZ/2\ZZ)\to \hat G_b\to G\to 1$ such that (\ref{seq_first_for_hat_G}) is its push-forward via the natural map $B(\ZZ/2\ZZ)\to B(\ZZ/4\ZZ)$. More properly, $\hat G_b$ is the geometric analog of $\Mp(\tilde V)$, but $\hat G_b$ will not appear in this paper.

We don't know if $\hat G$ admits a presentation as the stack quotient $G_1/G_0$ for a morphism $G_0\to G_1$ of algebraic groups over $k$, where $G_0$ is abelian and maps to the center of $G_1$  (we would rather expect that $\hat G$ corresponds to a nontrivial crossed module). 

 We have not found a relation with the K-theory (or the universal central extension of $\Sp_{2n}$ by $K_2$ constructed by Brylinski-Deligne \cite{BD}). Instead, our construction of $\hat G$ goes as follows. 
 
 Let $\cL(\tilde V)$ be the Greenberg realization of the $R$-scheme of free lagrangian submodules in $\tilde V$. First, we define a certain $\ZZ/4\ZZ$-gerb 
$$
\hat\cL(\tilde V)\to \cL(\tilde V)
$$ 
via the geometric Maslov index (cf. Section~6). 
It turns out that the corresponding class in $\H^2(\cL(\tilde V), \ZZ/4\ZZ)$ is invariant under $G$, but the gerb itself is not $G$-equivariant. Then $\hat G$ is defined as the stack  
of pairs $(g,\sigma)$, where $g\in G$ and $\sigma: g^*\hat\cL(\tilde V)\,\iso\, \hat\cL(\tilde V)$ is an isomorphism of $\ZZ/4\ZZ$-gerbs over $\cL(\tilde V)$. 

 We generalize the theory of canonical interwining operators from \cite{L2} to the case of characteristic two (cf. Section~7). This allows us to come up with a construction of the Weil category $W(\tilde V)$, which is a geometric analog of $\cH_{\psi}$. Here $W(\tilde V)$ is a category of certain perverse sheaves on $\hat\cL(\tilde V)\times H(V)$. The group stack $\hat G$ acts on $W(\tilde V)$ by functors. This action is geometric in the sense that it comes from the natural action of $\hat G$ on $\hat\cL(\tilde V)\times H(V)$.
 
 Similarly to the case of other characteristics, we also construct \select{the finite-dimensional theta-sheaf} $S_{\tilde V,\psi}$, which is a geometric analog of some matrix coefficient of the representation $\cH_{\psi}$. It serves as the key ingredient for the construction of $W(\tilde V)$. 
 
\medskip\noindent 
{\scshape Acknowledgements.} We are very grateful to P. Deligne for his comments and suggestions about the first version of this paper.  
   

\medskip\noindent
1.2 {\scshape Notation} Let $A$ be a commutative ring and $V$ a free $A$-module of rank $d$. As in (\cite{G}, Section~5.1), denote by $\Sym^{!2}(V)\subset V\otimes V$ the submodule of $S_2$-invariant vectors. It is the submodule spanned by the vectors of the form $v\otimes v$, $v\in V$. Let $\wedge^2(V)\subset V\otimes V$ be the submodule spanned by the vectors of the form $v\otimes u-u\otimes v$. Let $\Sym^{*2}(V)$ be the quotient of $V\otimes V$ by $\wedge^2(V)$. For any $n$ let
$$
\wedge^n(V)=\mathop{\cap}\limits_{1\le i\le n-1} V^{\otimes i-1}\otimes (\wedge^2 V)\otimes V^{\otimes d-1-i}
$$

 For any $n$ write $\Sym^{!n}(V)\subset V^{\otimes n}$ for the submodule of $S_n$-invariant tensors and $\Sym^{*n}(V)$ for the $A$-module of $S_n$-coinvariants of $V^{\otimes n}$.
We have canonically 
$$
\Sym^{!n}(V^*)\,\iso\, (\Sym^{*n}(V))^* \;\;\;\;\mbox{and}\;\;\;\; \wedge^n(V^*)\,\iso\, (\wedge^n V)^*
$$  
 
 The space $\Sym^{!2}(V^*)$ can be seen as the space of symmetric bilinear forms on $V$. This is the space of $A$-linear maps $\phi: V\to V^*$ such that $\phi^*=\phi$. 
 
 Let $V$ be a free $A$-module of rank $2n$. Say that $V$ is a symplectic if it is equipped with a non-degenerate bilinear form $\omega: V\times V\to A$ such that in a suitable base $(e_i, e_{-i}, 1\le i\le n)$ we have $\omega(e_i, e_{-j})=\delta_{ij}$, $\omega(e_i, e_j)=\omega(e_{-i}, e_{-j})=0$ for $i,j\in \{1,\ldots,n\}$ and $\omega(x,y)=-\omega(y,x)$. We call such base \select{a symplectic base}. If $V$ is a symplectic $A$-module then we trivialize $\omega_{\wedge}: \det V\,\iso\, A$ by $e\mapsto 1$, where $e=e_1\wedge\ldots\wedge e_n\wedge e_{-1}\wedge\ldots\wedge e_{-n}$ does not depend on a symplectic base $(e_i, e_{-i})$. For $\omega^n\in \wedge^{2n} M^*$ we have $n! e=\pm\omega^n$.
 
  The pairing between a free $A$-module and its dual is usually denoted by $\<\cdot,\cdot\>$. If $A$ is of characteristic $p$, write $V^{(p)}=V\otimes_A A$, where $A\to A$, $a\mapsto a^p$ is the Frobenius map.

\medskip\noindent
1.3 {\scshape Generalities on quadratic forms} Let $k$ be a field of  characteristic two, which is either finite or algebraically closed. Let $R$ be the ring of Witt vectors of length two over $k$. 

 Let $L$ be a finite-dimensional $k$-vector space.
Let $B_a(L^*)\subset (L\otimes L)^*$ be the subspace of bilinear forms $\phi$ on $L$ satisfying $\phi(x,x)=0$ for all $x\in L$. We call them \select{alternating bilinear forms}. Write $Q(L^*)$ for the $k$-space of quadratic forms on $L$. By definition, it is included into an exact sequence $0\to B_a(L^*)\to (L\otimes L)^*\to Q(L^*)\to 0$, where the second map sends a bilinear form $\phi$ to the quadratic form $x\mapsto \phi(x,x)$.

Write $Q_a(L^*)$ for the $k$-space of additive quadratic forms on $L$. The map $L^*\to Q_a(L^*)$ sending $y^*$ to the quadratic form $y\mapsto \<y,y^*\>^2$, $y\in L$ yields an isomorphism $(L^*)^{(2)}\,\iso\, Q_a(L^*)$. One has an exact sequence 
$
0\to B_a(L^*)\to\Sym^{!2}(L^*)\to Q_a(L^*)\to 0
$
of $k$-vector spaces. One also has an exact sequence 
$$
0\to Q_a(L^*)\to Q(L^*)\to B_a(L^*)\to 0,
$$
where the second map sends $q$ to the the bilinear form $(x,y)\mapsto q(x+y)-q(x)-q(y)$.

 Assume given a free $R$-module $\tilde L$ with an isomorphism $\tilde L\otimes_R k\,\iso\, L$ of $k$-vector spaces. Note that $B_a(L^*)\subset \Sym^{!2}(\tilde L^*)$ can be seen as the $R$-submodule consisting of $\varphi$ satisfying $\varphi(\tilde x,\tilde x)=0$ for all $\tilde x\in\tilde L$. Let $Q^!(L^*)$ be the quotient of $\Sym^{!2}(\tilde L^*)$ by $B_a(L^*)$, this is an $R$-module included into an exact sequence $0\to Q_a(L^*)\to Q^!(L^*)\to \Sym^{!2}(L^*)\to 0$ of $R$-modules. The action of $\GL(\tilde L)$ on $Q^!(L^*)$ factors through an action of $\GL(L)$.
 
  Equivalently, one may define $Q^!(L^*)$ as the $R$-module of maps $q: L\to R$ such that 
\begin{itemize} 
\item the map $b_q: L\times L\to R$ given by $b_q(x_1,x_2)=q(x_1+x_2)-q(x_1)-q(x_2)$ is bi-additive, and $b_q(ax_1, x_2)=\tilde ab_q(x_1,x_2)$ for any $\tilde a\in R$ over $a\in k$;
\item $q(ax)=\tilde a^2 q(x)$ for $x\in L$ and $\tilde a\in R$ over $a\in k$.
\end{itemize} 
These $R$-valued `quadratic forms' on $L$ have been considered, for example, in \cite{GH, Woo, LoRu}.

\medskip\noindent
1.4 We refer the reader to Section~2 for motivations, in Section~3 we give a detailed description of the main results and explain the structure of the paper.

\bigskip\medskip
\centerline{\scshape 2. Classical Weil representation and motivations}

\bigskip\noindent
2.1. In this section we remind the construction of the Weil representation in characteristic two following essentially \cite{GH}. This is our subject to geometrize.

\medskip\noindent
2.2. Let $k$ be a finite field of characteristic two with $q$ elements, $R$ be the ring of Witt vectors of length two over $k$. Let $\tilde V$ be a free $R$-module of rank $2n$ with symplectic form $\tilde\omega: \tilde V\otimes\tilde V\to R$. Set $V=\tilde V\otimes_R k$. Let $\omega: V\times V\to 2R$ be given by $\omega(x,y)=2\tilde\omega(\tilde x, \tilde y)$ for any $\tilde x,\tilde y\in \tilde V$ over $x,y\in V$.

Pick a bilinear form $\tilde\beta: \tilde V\times\tilde V\to R$ such that 
\begin{equation}
\label{def_tilde_beta}
\tilde\beta(\tilde x,\tilde y)-\tilde\beta(\tilde y,\tilde x)=\tilde\omega(\tilde x,\tilde y)
\end{equation}
for all $\tilde x,\tilde y\in \tilde V$. Let $\beta: V\times V\to R$ be the map $(x,y)\mapsto 2\tilde\beta(\tilde x,\tilde y)$ for any $\tilde x,\tilde y\in\tilde V$ over $x,y$. It gives rise to the Heisenberg group $H(V)=V\times R$ with operation 
\begin{equation}
\label{def_Heis_group}
(v_1,z_1)(v_2,z_2)=(v_1+v_2, z_1+z_2+\beta(v_1,v_2)),\;\;\;\; v_i\in V, z_i\in R
\end{equation}
Its center is $Z(H(V))=\{(0,z)\in H(V)\mid z\in R\}$. 

  Gurevich and Hadani consider the group of all automorphisms of $H(V)$ acting trivially on the center $Z(H(V))$. For the purposes of geometrisation, we modify their definition slightly as follows. Let $\ASp(V)$ be the set of pairs $(g,\alpha)$, where $g\in \Sp(V)$ and $\alpha: V\to R$ satisfies
\begin{itemize}
\item $\alpha(v_1+v_2)-\alpha(v_1)-\alpha(v_2)=\beta(g(v_1), g(v_2))-\beta(v_1,v_2)$ for all $v_i\in V$;
\item $\alpha(av)=\tilde a^2\alpha(v)$ for any $v\in V$ and $\tilde a\in R$ over $a\in k$.
\end{itemize}
An element $(g,\alpha)\in\ASp(V)$ yields an automorphism of $H(V)$ given by $(v,z)\mapsto (gv, z+\alpha(v))$. In this way $\ASp(V)$ maps injectively into the group of automorphisms of $H(V)$ acting trivially on $Z(H(V))$. The composition in $\ASp(V)$ is given by 
$$
(g,\alpha_g)(h,\alpha_h)=(gh, h^{-1}(\alpha_g)+\alpha_h)
$$
with $h^{-1}(\alpha_g)(v)=\alpha_g(hv)$ for all $v\in V$.
We will refer to $\ASp(V)$ as \select{the affine symplectic group}. For a $k$-vector space $L$ write
$$
FQ_a(L^*)=\{\alpha: L\to R\mid \alpha(x_1+x_2)=\alpha(x_1)+\alpha(x_2)\;\mbox{and}\; \alpha(ax)=\tilde a^2\alpha(x), \tilde a\in R\; \mbox{over}\; a\in k, x\in L\}
$$
An element of $FQ_a(L^*)$ writes in Witt coordinates as $(0, \alpha_1)$, where $\alpha_1: L\to k$ is additive and $\alpha_1(ax)=a^4\alpha_1(x)$ for all $a\in k, x\in L$. So, $FQ_a(L^*)\,\iso\, (L^*)^{(4)}$. The group $\ASp(V)$ fits into an exact sequence
$$
1\to FQ_a(V^*)\to \ASp(V)\to \Sp(V)\to 1
$$
Though it is not reflected in the notation, $\ASp(V)$ depends not only on $\tilde\omega$ but also on $\tilde\beta$. 
  
  Let $G=\Sp(\tilde V)$. We have a surjective homomorphism $\xi: G\to \ASp(V)$ sending $\tilde g$ to $(g, \alpha_{\tilde g})$, where $g\in \Sp(V)$ is the image of $\tilde g$, and $\alpha_{\tilde g}: V\to R$ is given by
$$
\alpha_{\tilde g}(v)=\tilde\beta(\tilde g\tilde v, \tilde g\tilde v)-\tilde\beta(\tilde v, \tilde v)
$$
for any $\tilde v\in \tilde V$ over $v\in V$.

\medskip\noindent
2.3 Fix a prime $\ell\ne 2$. Let $\psi: \ZZ/4\ZZ\to \Qlb^*$ be a faithful character. We denote by the same symbol $\psi$ the composition $R\toup{\tr}\ZZ/4\ZZ\to\Qlb^*$. A version of the Stone-von Neumann theorem holds in this setting, namely there exists a unique (up to isomorphism) irreducible $\Qlb$-representation of $H(V)$ with central character $\psi$ (\cite{GH}).

 Given an irreducible $\Qlb$-representation $(\rho, \cH_{\psi})$ of $H(V)$ with central character $\psi$, one gets in the usual way a version of the metaplectic group
\begin{multline*}
\Mp(\cH_{\psi})=\{(g, M[g])\mid g\in \ASp(V), M[g]\in\Aut(\cH_{\psi})\\
\!\!\!\!\rho(gh)\comp M[g]=M[g]\comp \rho(h)\;\mbox{for all}\; h\in H(V)\}    
\end{multline*}    
included into an exact sequence 
\begin{equation}
\label{seq_Mp(cH_psi)}
1\to \Qlb^*\to  \Mp(\cH_{\psi})\to \ASp(V)\to 1
\end{equation}
It comes together with the Weil representation of $\Mp(\cH_{\psi})$ on $\cH_{\psi}$ given by $(g, M[g])\mapsto M[g]\in\Aut(\cH_{\psi})$. 

 Among various results of \cite{Su, Is, HZ, GH} let us cite the following two, which are our main motivation for geometrisation. 
 
\begin{Pp}[\cite{GH}] 
\label{Pp_GH_meta_extensions}
1) There exists a group $\AMp(V)$ included into an exact sequence 
$$
1\to \mu_4(\Qlb)\to \AMp(V)\to \ASp(V)\to 1
$$ 
such that  (\ref{seq_Mp(cH_psi)}) is its push-forward via $\mu_4(\Qlb)\hook{} \Qlb^*$.\\
1) There is a group $\Mp(\tilde V)$, which is an extension of $G$ by $\mu_2(\Qlb)$ and a morphism of exact sequences extending $\xi$
$$
\begin{array}{ccccccc}
1\to  & \mu_4(\Qlb) & \to & \AMp(V) & \to & \ASp(V) & \to 1\\
 & \uparrow && \uparrow && \uparrow\lefteqn{\scriptstyle\xi}\\
 1\to & \mu_2(\Qlb) & \to & \Mp(\tilde V) & \to & G &\to 1
\end{array}
$$
\end{Pp}

\smallskip\noindent
2.4 {\scshape Models of the Weil representation}  

\medskip\noindent
2.4.1 Write $\cL(V)$ for the set of lagrangians in $V$. We modify slightly the notion of an enhanced lagrangian from \cite{GH} as follows. \select{An enhanced lagrangian} in $V$ is a pair $(L,\alpha)$, where $L\in \cL(V)$ and $\alpha: L\to R$ satisfies 
$$
\alpha(l_1+l_2)-\alpha(l_1)-\alpha(l_2)=\beta(l_1,l_2)
$$ 
for $l_i\in L$ and $\alpha(al)=\tilde a^2\alpha(l)$ for $l\in L, a\in k$ and $\tilde a\in R$ over $a$. Sometimes we will refer to it as \select{the enhanced structure} on $L$. Let $ELag(V)$ be the set of enhanced lagrangians in $V$. Then $ELag(V)$ is a torsor under the vector bundle on $\cL(V)$ whose fibre at $L$ is $FQ_a(L^*)$. 

 Given $L\in \cL(V)$, the subgroup $L\times R\subset H(V)$ is abelian. The novelty in characteristic two (compared to the case of other characteristics) is that this is not a direct product of subgroups. This is an extension of $L$ by $R$ in the sense of commutative unipotent group schemes over $k$, and an enhanced structure on $L$ yields a splitting of this extension.
 
 Let $(L, \alpha)$ be an enhanced lagrangian in $V$. Then $\tau: L\to H(V)$ given by $\tau(x)=(x,\alpha(x))$ is a group homomorphism. One associates to the enhanced lagrangian $(L,\alpha)$ a model of the Weil representation
$$
\cH_L=\{f: H(V)\to\Qlb\mid f(z\tau(x)h)=\psi(z)f(h), x\in L, z\in R, h\in H(V)\},
$$ 
on which $H(V)$ acts by right translations. Write $\rho: H(V)\to \Aut(\cH_L)$ for this action.
 
 The group $\ASp(V)$ acts on the left on $ELag(V)$, namely $g\in \ASp(V)$ sends $\tau: L\to H(V)$ to $\Ad_g\tau: gL\to H(V)$ given by 
$$
(\Ad_g\tau)(x)=g\tau(g^{-1}x),
$$
$x\in gL$. This action map is denoted by $\act:\ASp(V)\times ELag(V)\to ELag(V)$. The projection $ELag(V)\to\cL(V)$ is $\ASp(V)$-equivariant. Set
\begin{multline*}
\Mp(\cH_L)=\{(g, M[g])\mid g\in \ASp(V), M[g]\in\Aut(\cH_L)\\
\!\!\!\!\rho(gh)\comp M[g]=M[g]\comp \rho(h)\;\mbox{for}\; h\in H(V)\}    
\end{multline*}    
It fits into an exact sequence
\begin{equation}
\label{seq_for_Mp(H_L)}
1\to \Qlb^*\to \Mp(\cH_L) \to \ASp(V)\to 1
\end{equation}
 
\begin{Rem}
\label{Rem_splitting_one}
Each $g\in \ASp(V)$ yields an isomorphism $\cH_L\iso \cH_{gL}$ of $\Qlb$-vector spaces sending $f$ to $gf$. Here $(gf)(h)=f(g^{-1}h)$ for $h\in H(V)$. Write $St_L$ for the stabilizor of the enhanced lagrangian $L$ in $\ASp(V)$. Then $St_L$ acts naturally in $\cH_L$, namely $g\in St_L$ sends $f$ to $gf$. This is a splitting of (\ref{seq_for_Mp(H_L)}) over $St_L$.
\end{Rem}
  
\medskip\noindent
2.4.2  Write $\cL(\tilde V)$ for the set of free lagrangian $R$-submodules in $\tilde V$. Write $\epsilon: \cL(\tilde V)\to ELag(V)$ for the map sending $\tilde L$ to $(L,\alpha_{\tilde L})$, where $\alpha_{\tilde L}: L\to R$ is given by $\alpha_{\tilde L}(x)=\tilde\beta(\tilde x, \tilde x)$ for any $\tilde x\in \tilde L$ over $x\in L$. 

  The surjection $\xi: G\to \ASp(V)$ is compatible with $\epsilon$, namely the diagram 
\begin{equation}
\label{diag_action_of_ASp(V)_and_G}
\begin{array}{ccc}
G\times \cL(\tilde V) & \toup{\act} & \cL(\tilde V)\\
\downarrow\lefteqn{\xi\times \epsilon} && \downarrow\lefteqn{\scriptstyle \epsilon}\\
\ASp(V)\times ELag(V) & \toup{\act} & ELag(V)
\end{array}
\end{equation} 
commutes, where the top horizontal map sends $(\tilde g, \tilde L)$ to $\tilde g\tilde L$. For $\tilde L\in \cL(\tilde V)$ set also $\cH_{\tilde L}=\cH_L$, where $L=\epsilon(\tilde L)$. Let $P(\tilde L)$ be the stabilizor of $\tilde L$ in $G$. Then $\xi$ resticts to a homomorphism $\xi: P(\tilde L)\to St_L$, and Remark~\ref{Rem_splitting_one} yields a canonical splitting of the pull-back of (\ref{seq_for_Mp(H_L)}) under the composition $P(\tilde L)\hook{} G\toup{\xi}\ASp(V)$. The group $P(\tilde L)$ fits into an exact sequence $1\to \Sym^{!2}(\tilde L)\to P(\tilde L)\to \GL(\tilde L)\to 1$.     
     
\medskip\noindent
2.5 {\scshape Interwining operators}  

\medskip\noindent
For a pair $L_1, L_2\in ELag(V)$ consider the \select{non-normalized} interwining operator $F_{L_1L_2}: \cH_{L_2}\to \cH_{L_1}$ given by
\begin{equation}
\label{operator_F_transverse}
(F_{L_1L_2}[f])(h)=\sum_{m\in L_1} f(\tau_1(m)h) 
\end{equation}
for $f\in \cH_{L_2}, h\in H(V)$, where $\tau_i: L_i\to H(V)$ is the enhanced structure for $L_i$. It commutes with the actions of $H(V)$. It is easy to check that $F_{L_1,L_2}$ does not vanish iff $\tau_1$ and $\tau_2$ coincide on $L_1\cap L_2$, and in this case it is an isomorphism.
     
  Let now $L_1,L_2,L_3\in ELag(V)$ with $L_1\cap L_2=L_1\cap L_3=0$. We write $\tau_i: L_i\to H(V)$ for the  enhanced structure for $L_i$. Let $r: L_2\to L_1$ be the $k$-linear map such that $L_3=\{r(x)-x\in L_1\oplus L_2\mid x\in L_2\}$. Let us calculate 
$$
F_{L_1,L_2}F_{L_2,L_3}: \cH_{L_3}\to \cH_{L_1}
$$ 
Define the quadratic form 
$Q_{L_1,L_2,L_3}: L_2\to R$ by the following equality in $H(V)$
\begin{equation}
\label{eq_def_for_Q_123}
\tau_3(r(m)-m)\tau_2(m)\tau_1(-r(m))=(0, Q_{L_1,L_2,L_3}(m))
\end{equation}
for any $m\in L_2$. In other words, 
\begin{equation}
\label{formula_for_Q_L123}
Q_{L_1,L_2,L_3}(m)=\alpha_2(m)+\alpha_1(r(m))-\alpha_3(m-r(m))+\beta(m,r(m))
\end{equation}
Then $Q_{L_1,L_2,L_3}\in Q^!(L_2^*)$, and
(\ref{formula_for_Q_L123}) implies for $m_1,m_2\in L_2$
$$
Q_{L_1,L_2,L_3}(m_1+m_2)-Q_{L_1,L_2,L_3}(m_1)-Q_{L_1,L_2,L_3}(m_2)=\omega(r(m_1), m_2)
$$ 
This is a symmetric bilinear form $L_2\times L_2\to 2R$. 
 
 Define the Gauss sum of $Q_{L_1,L_2,L_3}$ by
\begin{equation}
\label{Gauss_sum_first_classical}
C(L_1,L_2,L_3)=\sum_{m\in L_2} \psi(Q_{L_1,L_2,L_3}(m))
\end{equation}
\begin{Lm} 
Assuming only $L_1\cap L_2=L_1\cap L_3=0$ one has
$$
F_{L_1,L_2}F_{L_2,L_3}=C(L_1,L_2,L_3)F_{L_1,L_3}
$$
So, $C(L_1,L_2,L_3)F_{L_2,L_1}F_{L_1,L_3}=q^n F_{L_2,L_3}$.
\end{Lm}
\begin{Prf}
Let $f\in \cH_{L_3}$. For $h\in H(V)$ we get 
$$
(F_{L_1,L_2}F_{L_2,L_3}f)(h)=\sum_{u\in L_1, m\in L_2} f(\tau_2(m)\tau_1(u)h)=\sum_{u\in L_1, m\in L_2} f(\tau_3(r(m)-m)\tau_2(m)\tau_1(u)h)
$$
We make a change of variables replacing $(u,m)$ by $(v,m)$, where $v\in L_1$ is given by $v=r(m)+u$. The above sum equals
$$
\sum_{v\in L_1, m\in L_2} f(\tau_3(r(m)-m)\tau_2(m)\tau_1(-r(m))\tau_1(v)h)
$$
Our first assertion follows now from (\ref{eq_def_for_Q_123}).  The second assertion follows from the fact that $F_{L_2, L_1}F_{L_1,L_2}=q^n$.
\end{Prf}

\medskip

\begin{Rem} Given $L_1\in ELag(V)$, the variety $\{L\in ELag(V)\mid L\cap L_1=0\}$ is naturally a homogeneous space under $Q^!(L_1)$. Namely, given $L_2$ and $L_3$ in this variety one gets $Q_{L_1,L_2,L_3}$ as above. Identify further $L_1$ with $L_2^*$ via the form $\omega$, so $Q_{L_1,L_2,L_3}\in Q^!(L_1)$.
 
  Conversely, given $L_2$ in this variety and $Q_{L_1,L_2,L_3}\in Q^!(L_2^*)$, let $\phi$ be its image in $\Sym^{!2}(L_2^*)$. Define $r: L_2\to L_1$ by $\omega(r(m_1), m_2)=\phi(m_1,m_2)$ for all $m_i\in L_2$. Let $L_3=\{r(x)-x\in L_1\oplus L_2\mid x\in L_2\}$. Now there is a unique $\tau_3: L_3\to H(V)$, $\tau_3(y)=(y, \alpha_3(y))$ such that (\ref{eq_def_for_Q_123}) holds. The fact that this $\tau_3$ is a homomorphism is checked using (\ref{formula_for_Q_L123}).
\end{Rem}

\begin{Rem} For any pair of enhanced lagrangians $L_1,L_2\in ELag(V)$ define the interwining operator $F^{\flat}_{L_1L_2}: \cH_{L_2}\to \cH_{L_1}$ as follows. There is $w\in V$ (its image in $V/(L_1+L_2)$ is uniquely defined) such that for all $x\in L_1\cap L_2$ we have $\alpha_2(x)-\alpha_1(x)=\Fr(\omega(x,w))$, here $\Fr: R\to R$ is the Frobenius. 

 Write $\bar L_1$ for the subgroup $L_1\times R\subset H(V)$, similarly for $\ov{L_1\cap L_2}\subset H(V)$.
Let $\psi_1: \bar L_1\to\Qlb^*$ be the unique character such that $\psi_1(\tau_1(x)(0,z))=\psi(z)$ for all $x\in L_1, z\in R$. Given $f\in \cH_{L_2}$, the function sending $u\in \bar L_1$ to 
$$
f((w,0)uh)\psi_1^{-1}(u)
$$
actually depends only on the image of $u$ in $\ov{L_1\cap L_2}\backslash \bar L_1$. Then we set
\begin{equation}
\label{operator_F_general}
(F^{\flat}_{L_1L_2}f)(h)=\sum_{u\in \ov{L_1\cap L_2}\backslash \bar L_1} f((w,0)uh)\psi_1^{-1}(u)
\end{equation}
To see that $F^{\flat}_{L_1,L_2}$ is non zero, take $f\ne 0$ supported at $\bar L_2$ then $(F^{\flat}_{L_1L_2}f)(w,0)\ne 0$. Finally, if $L_1\cap L_2=0$ take $w=0$ then (\ref{operator_F_general})  coincides with (\ref{operator_F_transverse}). The operator (\ref{operator_F_general}) does depend on a choice of $w$.
\end{Rem}
   

\bigskip

\centerline{\scshape 3. Main results}

\bigskip\noindent
3.1 {\scshape Notation} From now on $k$ denotes an algebraically closed field of characteristic two. Write $W_2$ for the $k$-scheme of Witt vectors of length 2 over $k$. Let $R=W_2(k)$ be the corresponding ring of Witt vectors. Fix a prime $\ell\ne 2$. Let $\psi: \ZZ/4\ZZ\to \Qlb^*$ be a faithful character. For a $k$-stack of finite type $S$ write $\D(S)$ for the bounded derived category of $\Qlb$-sheaves on $S$. Write $\P(S)\subset \D(S)$ for the full subcategory of perverse sheaves. Write $\DD: \D(S)\to\D(S)$ for the Verdier duality. Since we are working over an algebraically closed field, we systematically ignore the Tate twists. (Certain isomorphisms that we call canonical actually contain the Tate twists, for example those of Lemma~\ref{Lm_first_true_convolution}). Write $\mu_2$ for the group scheme $\Spec k[x]/(x^2-1)$ over $k$.

 Remind the definition of the Artin-Schreier-Witt local system $\cL_{\psi}$ on $W_{2,k}$. One has the Lang isogeny $La: W_2\to W_2$ sending $x$ to $\Fr(x)-x$, where $\Fr: W_2\to W_2$ is the Frobenius morphism (\cite{S}, section~0.1.2). It can be seen as a $\ZZ/4\ZZ$-torsor over $W_2$. Denote by $\cL_{\psi}$ the smooth $\Qlb$-sheaf of rank one on $W_2$ obtained from this torsor via extension of scalars $\psi: \ZZ/4\ZZ\to\Qlb^*$. This is a character sheaf, for the sum $s: W_2\times W_2\to W_2$ we have $s^*\cL_{\psi}\,\iso\, \cL_{\psi}\boxtimes\cL_{\psi}$ canonically, and $\cL_{\psi}$ is canonically trivialized at  $0\in W_2$.
  
  For a scheme $Z$ over $R$ we denote by the same symbol $Z$ its Greenberg realization over $k$ \cite{Gr} (the precise meaning being understood from the context). 
  
  For an $R$-module of finite type $M$ the same symbol $M$ stands for the $k$-scheme whose set of $k$-points is $M$. It can be defined as follows. Pick a resolution $M_{-1}\toup{f} M_0$ of $M$ by free $R$-modules of finite type. The $k$-scheme associated to $M$ is defined as the cokernel of the corresponding morphism between the Greenberg realizations of $M_i$. One checks that $k$-scheme so obtained is defined up to a unique isomorphism.
   
  For a morphism $f: Y_1\to Y_2$ of irreducible $k$-schemes of finite type write $\dimrel(f)=\dim Y_1-\dim Y_2$.  
 
\medskip   
\noindent
3.2  In Section~4 we study the geometric analogs of Gauss sums attached to $R$-valued quadratic forms on a finite-dimensional $k$-vector space $L$. Namely, we introduce an irreducible perverse sheaf $\bar S_{\psi}$ on $Q^!(L^*)$, whose fibre at $q\in Q^!(L^*)$ is the Gauss sum attached to $q$ and the character $\psi$. The scheme $Q^!(L^*)$ is naturally stratified, we describe the restriction of $S_{\psi}$ to each stratum (Proposition~\ref{Pp_restriction_strata}). Generically, $S_{\psi}$ is a (shifted) local system of rank one and order four, and we describe generically the $\ZZ/4\ZZ$-covering, on which the corresponding local system trivializes (cf. Section~4.5). The sheaf $\bar S_{\psi}$ is $\GL(L)$-equivariant, write $\bar\SSS_{\psi}$ for the perverse sheaf on the stack quotient $Q^!(L^*)/\GL(L)$ equipped with an isomorphism $\pr^*\bar\SSS_{\psi}[\dimrel(\pr)]\,\iso\, \bar S_{\psi}$ for the projection $\pr: Q^!(L^*)\to Q^!(L^*)/\GL(L)$.

 In Section~5 we generalize the results of Thomas \cite{TT} on the Maslov index (of a finite collection of lagrangian subspaces in a symplectic space) to the case of characteristic two. More precisely, we consider a free $R$-module $\tilde V$ of rank $2n$ with a symplectic form. To a finite collection of free lagrangian $R$-submodules $\tilde L_1,\ldots,\tilde L_m\subset \tilde V$ we attach $R$-modules with symmetric bilinear forms $K_{1,\ldots,m}$ and $T_{1,\ldots,m}$, here $T_{1,\ldots,m}$ is the quotient of $K_{1,\dots,m}$ by the kernel of the form. This is the Maslov index of the collection $\{\tilde L_i\}$. In addition to the standard properties of the Maslov index (Propositions~\ref{Pp_chain_condition} and \ref{Con_Maslov_is_cocycle}), we also obtain some new isometries in the case $m=4$, 
which actually hold also in other characteristics. (These new isometries are used in Section~6.3.2 for the construction of the gerb $\hat Y$).
  
   Let $\cL(\tilde V)$ be the Greenberg realization of the $R$-scheme of lagrangian free $R$-submodules in $\tilde V$. In Section~6 we use the Maslov index to construct a $\ZZ/4\ZZ$-gerb $\hat\cL(\tilde V)\to \cL(\tilde V)$. Let $G$ be the Greenberg realization of $\Sp(\tilde V)$. The gerb $\hat\cL(\tilde V)$ is not $G$-equivariant, but rather gives rise to a central extension
$$
1\to B(\ZZ/4\ZZ)\to \hat G\to G\to 1,
$$      
of group stacks over $k$. We call $\hat G$ \select{the metaplectic group} (in a sense, this is a geometric analog of $\Mp(\tilde V)$ from Proposition~\ref{Pp_GH_meta_extensions}).

 Set $Y=\cL(\tilde V)\times \cL(\tilde V)$. Let $\hat Y\to Y$ be the $\ZZ/4\ZZ$-gerb obtained from $\hat\cL(\tilde V)\times \hat\cL(\tilde V)$ by extending the structure group via $\ZZ/4\ZZ\times \ZZ/4\ZZ\to \ZZ/4\ZZ$, $(a,b)\mapsto b-a$. We show that $\hat Y$ is naturally $G$-equivariant, that is, can be seen as a gerb over the stack quotient $Y/G$. We construct an irreducible $G$-equivariant perverse sheaf $S_{\tilde V,\psi}$ on $\hat Y$, which we call \select{the finite-dimensional theta sheaf} (cf. Definition~\ref{Def_sheaf_S_V_psi}).
 
  Let $U_{23}$ be the scheme classifying $(\tilde L_1,\tilde L_2,\tilde L_3)\in \cL(\tilde V)^3$ with $\tilde L_1\cap\tilde L_2=\tilde L_1\cap\tilde L_3=0$. We define a morphism $\hat\nu_{23}: U_{23}\to \hat Y$ extending the projection $\nu_{23}: U_{23}\to Y$, $(\tilde L_1,\tilde L_2,\tilde L_3)\mapsto (\tilde L_2, \tilde L_3)$ (cf. Section~6.3.1). 
  
  Let $\tilde L$ be a free $R$-module of rank $n$, set $L=\tilde L\otimes_R k$. The Maslov index of a triple of lagrangians yields a morphism to the stack quotient $\pi_U: U_{23}\to \Sym^{!2}(\tilde L^*)/\GL(\tilde L)$ (cf. Section~5.4).
    
  The key property of $S_{\tilde V,\psi}$ is Corollary~\ref{Cor_theta_descent}, which is the main result of Section~6. It establishes a canonical isomorphism of perverse sheaves on $U_{23}$
$$
\pi_U^*p_L^*\bar\SSS_{\psi}[\dimrel(p_L\comp\pi_U)]
\,\iso\,\hat\nu_{23}^*S_{\tilde V,\psi}[\dimrel(\nu_{23})],
$$
where 
$
p_L: \Sym^{!2}(\tilde L)/\GL(\tilde L)\to Q^!(L^*)/\GL(L)
$ 
is the natural map. So, similarly to the case of other characteristics, the Gauss sum of the Maslov index of $(\tilde L_1,\tilde L_2,\tilde L_3)\in U_{23}$ is `almost independent' of $\tilde L_1$. 

 Let $V=\tilde V\otimes_R k$. The Heisenberg group $H=H(V)=V\times R$ with operation (\ref{def_Heis_group}) is viewed as an algebraic group over $k$ (an extension of $V$ by the Greenberg realization of $R$). 
 
 In Section~7 we generalize the theory of canonical interwining operators (\cite{L2}) to the case of characteristic two. Main result here is Theorem~\ref{Th_CIO}, which establishes the existence of an irreducible perverse sheaf $F\in \P(\hat\cL(\tilde V)\times\hat\cL(\tilde V)\times H)$ of canonical interwining operators between the Heisenberg models of the Weil representation. The proof follows the strategy from \cite{L2}. However, a new technical point is that we need to consider an action of an algebraic group stack on an algebraic stack (our perverse sheaves are equivariant under the action of a group stack\footnote{The actions we consider are `free' in a suitable sense, so that the quotient could be defined as a 1-stack (\cite{D}, Section~2.4.4).}). Generalities on such actions are collected in Appendix~A. The theta-sheaf $S_{\tilde V,\psi}$ is one of the main ingredients in the construction of $F$. 
 
   Finally, we construct a category $W(\tilde V)$ of certain perverse sheaves on $\hat\cL(\tilde V)\times H$, which provides a geometric analog of the Weil representation $\cH_{\psi}$ from Section~2.3. The group stack $\hat G$ acts on $W(\tilde V)$ by functors. This action is \select{geometric} in the sense that it comes from a natural action of $\hat G$ on $\hat\cL(\tilde V)$. As in \cite{L2}, we also construct the \select{non-ramified Weil category} $W(\hat\cL(\tilde V))$, which is a category of certain perverse sheaves on $\hat\cL(\tilde V)$. This is another geometric realization of the Weil representation (one has an obvious functor $W(\tilde V)\to W(\hat\cL(\tilde V))$, we don't know if this is an equivalence).
   
  In Appendix B we turn back to the classical setting and show that the Weil representation from Section~2 is obtained by some reduction from the Weil representation over the non archimedian local field of characteristic zero and residual characteristic two.

\medskip\noindent
3.3 {\scshape Comments on open problems.} From our point of view, the following problems would be interesting to solve. Find an interpretation of the construction of $\hat G$ in terms of $K$-theory (as for characteristics different from two). Find a description of $\hat\cL(\tilde V)$ and of $\hat Y$ as the moduli stack classifying some geometric objets related to $\tilde V$. Geometrization of the Howe correspondence in the case of characteristic two. If there could be one, find a description of $\hat G$ by a crossed module.
        
\bigskip

\centerline{\scshape 4. The sheaf $\bar S_{\psi}$}

\bigskip\noindent
4.1 Let $\tilde L$ be a free $R$-module of rank $n$, set $L=\tilde L\otimes_R k$. Let $\pi: \tilde L\to \Sym^{*2}(\tilde L)$ be the map sending $x$ to $x\otimes x$, here $\Sym$ is taken over $R$. The map $\pi$ is constant along the fibres of the projection $\tilde L\to L$, so yields a map $\bar\pi: L\to \Sym^{*2}(\tilde L)$. The latter map induces a closed immersion
of $k$-schemes $\bar\pi: L^{(2)}\to \Sym^{*2}(\tilde L)$.  
   
   Since $\Sym^{*2}(\tilde L)$ and $\Sym^{!2}(\tilde L^*)$ is a dual pair of commutative unipotent group schemes over $k$, one has the Fourier transform functor 
$$
\Four_{\psi}: \D(\Sym^{*2}(\tilde L))\to \D(\Sym^{!2}(\tilde L^*))
$$
introduced in \cite{S}. In this particular case, one also has the evaluation map $\ev: \Sym^{*2}(\tilde L)\times \Sym^{!2}(\tilde L^*)\to W_2$, and $\Four_{\psi}$ is given by
$$
\Four_{\psi}(K)=\pr_{2!}(\pr_1^*K\otimes \ev^*\cL_{\psi})[n(n+1)]
$$
for $K\in \D(\Sym^{*2}(\tilde L))$. By \cite{S}, $\Four_{\psi}$ is an equivalence, commutes with Verdier duality and preserves perversity. 

\begin{Rem}
\label{Rem_G-equivariance_equivalence}
Let $\cZ$ be an algebraic stack locally of finite type, $G\to\cZ$ be a group scheme of finite type and smooth of relative dimension $m$ over $\cZ$. Let $f: \cY\to\cZ$ be a $G$-torsor over $\cZ$. Then the functor $K\mapsto 
f^*K[m]$ is an equivalence of the category of perverse sheaves on $\cZ$ with the category of $G$-equivariant perverse sheaves on $\cY$ (\cite{L1}, A.2).
\end{Rem}

\begin{Def} As in (\cite{L1}, Section~4.1)
set
$$
S_{\psi}=\Four_{\psi}(\bar\pi_!\Qlb[n])
$$
Since $\bar\pi_!\Qlb[n]$ is an irreducible perverse sheaf, so is $S_{\psi}$. The sheaf $S_{\psi}$ is naturally $\GL(\tilde L)$-equivariant (by $\GL(\tilde L)$ we mean here its Greenberg realisation over $k$). Let $q_{\tilde L}: \Sym^{!2}(\tilde L^*)\to \Sym^{!2}(\tilde L^*)/\GL(\tilde L)$ be the stack quotient. Write $\SSS_{\psi}$ for the perverse sheaf on $\Sym^{!2}(\tilde L^*)/\GL(\tilde L)$ equipped with an isomorphism 
$$
q_{\tilde L}^*\SSS_{\psi}[\dimrel(q_{\tilde L})]\,\iso\, S_{\psi}
$$ 
By Remark~\ref{Rem_G-equivariance_equivalence}, the perverse sheaf $\SSS_{\psi}$ is defined up to a unique isomorphism.
\end{Def}

 Set $Q^*(L)=\Hom_R(Q^!(L^*), R)$. 
So, $Q^*(L)\subset \Sym^{*2}(\tilde L)$ is the $R$-submodule of those $A\in \Sym^{*2}(\tilde L)$ which vanish on $B_a(L^*)$. The map $\bar\pi: L\to \Sym^{*2}(\tilde L)$ factors as
$$
L\toup{\pi_Q} Q^*(L)\hook{} \Sym^{*2}(\tilde L) 
$$
Again, $Q^*(L)$ and $Q^!(L^*)$ is a dual pair of commutative unipotent group schemes over $k$, one has the evaluation map $\ev: Q^*(L)\times Q^!(L^*)\to R$ and the Fourier transform
$$
\Four_{\psi}: \D(Q^*(L))\to \D(Q^!(L^*))
$$
as above. Let $\bar S_{\psi}$ be the irreducible perverse sheaf on $Q^!(L^*)$ defined by 
$$
\bar S_{\psi}=\Four_{\psi}(\pi_{Q !}\Qlb[n])
$$
Clearly, $\bar S_{\psi}$ is $\GL(L)$-equivariant. We have canonically $\DD \bar S_{\psi}\,\iso\, \bar S_{\psi^{-1}}$. 
We write $_n{\bar S}_{\psi}$ if we need to express the dependence on $n$. One has a canonical isomorphism 
$$
\pr^*\bar S_{\psi}[n(n+1)/2-n]\,\iso\, S_{\psi},
$$
where $\pr: \Sym^{!2}(\tilde L^*)\to Q^!(L^*)$ is the projection.

 We say that \select{the Gauss sum} for $q\in Q^!(L^*)$ is the $*$-fibre of $\bar S_{\psi}$ at $q$. The Gauss sum for $\tilde\phi \in\Sym^{!2}(\tilde L^*)$ is the Gauss sum for its image in $Q^!(L^*)$.

\medskip\noindent
4.2 Let $Q^!_0(L^*)\subset Q^!(L^*)$ be the open subscheme of $q\in Q^!(L^*)$ whose image in $\Sym^{!2}(L^*)$ is a non degenerate symmetric bilinear form, that is, a symmetric isomorphism $L\,\iso\, L^*$ of $k$-vector spaces.

\begin{Lm} Over $Q^!_0(L^*)$ the sheaf $\bar S_{\psi}$ is a rank one local system placed in usual degree $-n-n(n+1)/2=-\dim Q^!(L^*)$. One has canonically over $Q^!_0(L^*)$
$$
S_{\psi}\otimes S_{\psi^{-1}}\,\iso\, \Qlb[-3n-n^2]
$$
\end{Lm}
\begin{Prf}
Let us explain the argument at the level of functions, its geometrization is straightforward. Let $\phi: \tilde L\to \tilde L^*$ be a symmetric isomorphism of $R$-modules.
For $l\in L$ write $\tilde l$ for any lifting of $l$ to an element of $\tilde L$. Let us calculate
\begin{equation}
\label{expr_sum_twice}
\sum_{l,u\in L} \psi(\<\tilde l, \phi(\tilde l)\>-\<\tilde u, \phi(\tilde u)\>)
\end{equation}
For $f(\tilde l, \tilde u):=\<\tilde l, \phi(\tilde l)\>-\<\tilde u, \phi(\tilde u)\>$ and $\tilde v\in\tilde L$ we have 
$$
f(\tilde l+\tilde v, \tilde u+\tilde v)=f(\tilde l, \tilde u)+2\<\phi(\tilde l-\tilde u), \tilde v\>
$$
First, summate along the fibres of the map $L\times L\to L$, $(l,u)\mapsto l-u$. The above formula shows that the result is supported by $\{0\}\subset L$. So, (\ref{expr_sum_twice}) equals
$$
\sum_{l\in L} \psi(\<\tilde l, \phi(\tilde l)\>-\<\tilde l, \phi(\tilde l)\>)=\sum_{l\in L} 1 
$$
\end{Prf}

\smallskip

 The group scheme $W_2^*$ of invertible elements with respect to the multiplication identifies canonically with $\Gm\times\Ga$. Let $\bar\cB_2$ be the $k$-stack classifying a rank one free $R$-module $\tilde W$ with a bilinear form $\tilde W\otimes\tilde W\to R$. This is the stack quotient $W_2/W_2^*$, where $b\in W_2^*$ acts on $a\in W_2$ as $b^2 a$. Let $\cB_2\subset \bar\cB_2$ be the open substack given by the condition that the bilinear form in non degenerate. Recall the group scheme $\mu_2$ from Section~3.1.
 
\begin{Lm} There is a canonical isomorphism of $k$-stacks $\cB_2\,\iso\, B(\mu_2\times\Ga)\times \A^1$.
\end{Lm}
\begin{Prf} An $S$-point of $\cB_2$ is a $W_2^*$-torsor  
on $S$ with trivialization of its tensor square. A $W_2^*$-torsor is a pair: a $\Gm$-torsor $\cF_1$ and a $\Ga$-torsor $\cF_2$. Since $k$ is of characteristic two, for any $\Ga$-torsor $\cF_2$, its tensor square is canonically trivialized, and the additional trivialization yields a point of $\A^1$. Our assertion follows.
\end{Prf}

 \medskip
 
\begin{Lm} There is a natural map $\Disc: Q^!(L^*)\to \bar\cB_2$, its restriction to $Q_0^!(L^*)\subset Q^!(L^*)$ factors as $Q_0^!(L^*)\toup{\Disc} \cB_2\hook{}\bar\cB_2$.
\end{Lm} 
\begin{Prf}
The map sending a symmetric bilinear form $\tilde \phi:\tilde L\to\tilde L^*$ to $\det\tilde \phi: \det\tilde L\to \det\tilde L^*$ factors through $\Sym^{!2}(\tilde L^*)\to Q^!(L^*)$. Indeed, fix a base in $\tilde L$, so view $\Sym^{!2}(\tilde L^*)$ as symmetric $n\times n$-matrices $B$ over $R$.  Then $B_a(L^*)$ becomes the set of zero-diagonal symmetric matrices $C$ with entries in $2R$. For such $B$ and $C$ we claim that $\det(B+C)=\det B$. To see this, one must prove
$$
\sum_{\sigma\in S_n} (\sign\sigma)\sum_{j=1}^n c_{j,\sigma(j)}\prod_{i=1, i\ne j}^n b_{i,\sigma(i)}=0
$$
in $R$. If $\sigma\ne \sigma^{-1}$ then the contributions of $\sigma$ and $\sigma^{-1}$ are the same, so there remains the sum over $\sigma\in S_n$ such that $\sigma^2=\id$. We claim that for $\sigma\in S_n$ with $\sigma^2=\id$ we have
$$
\sum_{j=1}^n c_{j,\sigma(j)}\prod_{i=1, i\ne j}^n b_{i,\sigma(i)}=0
$$
Indeed, the above sum is actually over those $j$ which are not fixed by $\sigma$. Such $j$ are divided into pairs $(j,\sigma(j))$, and the contribution of each pair vanishes. 
\end{Prf} 

\medskip

 We will refer to $\Disc: Q^!(L^*)\to \bar\cB_2$ as \select{the discriminant map}. Write $\disc: Q_0^!(L^*)\to \A^1$ for the composition $Q_0^!(L^*)\toup{\Disc} \cB_2\to\A^1$, where the second map is the projection.
 
\medskip\noindent
4.3  Remind the classical result of Albert (\cite{A}). If $n$ is odd then $\GL(L)$ acts transitively on the variety of non-degenerate symmetric bilinear forms $\phi$ on $L$. If $n$ is even then $\GL(L)$ has two orbits on this variety:  they are distinguished by the fact that the image of $\phi$ in $Q_a(L^*)$ vanishes or not. 

 For any $n$ let $\Sym^{!2}(L^*)_0\subset \Sym^{!2}(L^*)$ be the open subscheme given by the condition that the image of $\phi\in \Sym^{!2}(L^*)$ in $Q_a(L^*)$ does not vanish, and $\phi$ is non-degenerate. This is the open orbit of $\GL(L)$ on $\Sym^{!2}(L^*)$.

 Let $U(L^*)\subset Q^!(L^*)$ be the open subscheme classifying $q\in Q^!(L^*)$ such that its image $\phi\in \Sym^{!2}(L^*)$ lies in $\Sym^{!2}(L^*)_0$. For $u\in k$ denote by $Q^!_0(L^*)_u\subset Q^!_0(L^*)$ the closed subscheme given by the condition that the discriminant equals $u$. 
Write also $U(L^*)_u$ for the fibre of the disriminant map
$\disc: U(L^*)\to\A^1$ over a $k$-point $u$.  

 Let $\phi\in \Sym^{!2}(L^*)_0$ be a $k$-point and let $r(y)=\phi(y,y)$, $y\in L$ be the corresponding additive quadratic form. Let $
y^*\in L^*$ be such that $r(y)=\<y, y^*\>^2$. Write $L_1$ for the kernel of $y^*: L\to k$. 

 If $n$ is odd then the restriction of $\phi$ to $L_1$ is non-degenerate, so $\phi$ is a non-degenerate alternating form on $L_1$. Writing $L_2=L_1^{\perp}$ we get $L\,\iso\, L_1\oplus L_2$. The stabilizor $\OO(\phi)$ of $\phi$ in $\GL(L)$ is $\Sp(L_1,\phi)\times \mu_2$, it is  not reduced, and its reduced part identifies with $\Sp(L_1,\phi)$. 

 If $n$ is even then $L_2=L_1^{\perp}\subset L_1$. Again, the group $\OO(\phi)$ is not reduced (it acts on $L_2$ via a quotient $\mu_2$). Its reduced part is an extension of $\Sp(L_1/L_2,\phi)$ by a unipotent group $U_{\phi}$. The group $U_{\phi}$ fits into an exact sequence $1\to L_2^{\otimes 2}\to U_{\phi}\to \Hom(L_1/L_2, L_2)\to 1$, so $U_{\phi}$ looks like a Heisenberg group. Actually, the latter exact sequence splits, because $k$ is of characteristic two.

\begin{Lm} 
\label{Lm_stabilizors}
i) Let $\tilde\phi \in \Sym^{!2}(\tilde L^*)$ be a lifting of $\phi\in \Sym^{!2}(L^*)_0$. Let $q\in U(L^*)$ be the image of $\tilde\phi$. 
The stabilizor $\OO(q)$ of $q$ in $\GL(L)$ is as follows. 
\begin{itemize}
\item[1)] If $n$ is odd then this is $\OO(L_1,q)\times\mu_2\subset \Sp(L_1,\phi)\times\mu_2$. 
\item[2)] If $n$ is even then we have two cases. First, if  $q$ is nonzero on $L_2$, then $\OO(q)$ fits into an exact sequence $1\to \ZZ/2\ZZ\to \OO(q)\to \Sp_{n-2}\to 1$. Second, if $q$ vanishes on $L_2$ then $q$ yields a quadratic form $\bar q: L_1/L_2\to 2R$, and $\OO(q)$ fits into an exact sequence $1\to \Hom(L_1/L_2, L_2)\to \OO(q)\to \OO(L_1/L_2, \bar q)\to 1$. The bilinear form on $L_1/L_2$ associated to $\bar q$ is symplectic.
\end{itemize}
ii) If $n$ is odd then we have a transitive action of $\GL(L)\times \A^1$ on $U(L^*)$, and the map $\disc: U(L^*)\to\A^1$ is $\A^1$-equivariant, where $\A^1$ acts on itself by translations.   

Assume $n$ even. Then $q$ vanishes on $L_2$ iff the discriminant of $q$ equals 
\begin{equation}
\label{disc_critical_value}
\left\{
\begin{array}{cl}
0, & n=0\mod 4\\
1, & n=2\mod 4
\end{array}
\right.
\end{equation}
Besides, for any $u\in k$ and any $n>0$ the action of $\GL(L)$ on $U(L^*)_u$ is transitive.
\end{Lm} 
\begin{Prf}   
i) For $x,y\in L$ one has
\begin{equation}
\label{eq_comparing}
q(x+y)-q(x)-q(y)=\tilde\phi(\tilde x+\tilde y, \tilde x+\tilde y)-\tilde\phi(\tilde x,\tilde x)-\tilde\phi(\tilde y,\tilde y)=2\tilde\phi(\tilde x,\tilde y)=(0,\phi(x,y)^2),
\end{equation}
the last expression being written in Witt coordinates.\\
1) Let $n$ be odd. Then $\phi$ yields a decomposition $L=L_1\oplus L_2$ as above, and $\phi$ is alternating on $L_1$. The restriction of $q$ to $L_1$ is a map $q_1: L_1\to 2R$. By (\ref{eq_comparing}), we see that the symmetric bilinear form associated to $q$ is $\phi$. It follows that $\OO(L_1, q)\subset \Sp(L_1, \phi)$.    
The stabilizor of $(L_2, q)$ is $\mu_2$.

\medskip\noindent
2) Let $n$ be even. Then $\phi$ yields a decomposition $L_2\subset L_1\subset L$ as above.  The restriction of $q$ to $L_1$ is a map $q_1: L_1\to 2R$. First, assume $q$ nontrivial on $L_2$, then the intersection of the unipotent radical of $\OO(\phi)$ with $\OO(q)$ is isomorphic to $\ZZ/2\ZZ$ (realized as a closed subgroup of $\Ga$). So, $\dim\OO(q)\le \dim\Sp_{n-2}=(n-2)(n-1)/2$. It follows that $\dim(\OO(\phi)/\OO(q))\ge n-1$, because $\dim\OO(\phi)=n(n-1)/2$. Let $\cY_{\phi}$ be the preimage of $\phi$ under $Q^!(L^*)\to \Sym^{!2}(L^*)$. Let $\cY_{\phi, q}$ be the closed subscheme in $\cY_{\phi}$ given by the condition that the discriminant equals the discriminant of $q$. Since $\dim\cY_{\phi,q}=n-1$, we learn that each orbit of $\OO(q)$ on $\cY_{\phi,q}$ is open. So, one gets an exact sequence $1\to \ZZ/2\ZZ\to \OO(q)\to \Sp_{n-2}\to 1$.

 Now consider the case when $q$ vanishes on $L_2$. Then for $x\in L_2, y\in L_1$ we get from (\ref{eq_comparing}) that $q(x+y)-q(y)=0$. So, $q$ yields a map $\bar q: L_1/L_2\to 2R$ given by $\bar q(y \!\mod L_2)=q(y)$ for any $y\in L_1$. It is easy to see that $\OO(q)\cap L_2^{\otimes 2}=0$. Since $\dim \OO(\bar q)=(n-2)(n-3)/2$, we get 
$$
\dim\OO(q)\le \dim\Hom(L_1/L_2, L_2)+\dim\OO(\bar q)=n-2+(n-2)(n-3)/2
$$ 
So, $\dim GL(L)/\OO(q)\ge (n-1)+n(n+1)/2$. On the other hand, $GL(L)/\OO(q)$ is contained in the locus of $Q^!_0(L^*)$ with fixed discriminant, which is of dimension $\le (n-1)+n(n+1)/2$. Thus, $1\to \Hom(L_1/L_2, L_2)\to \OO(q)\to \OO(\bar q)\to 1$ is exact.

\medskip\noindent
ii) If $n$ is odd we let $a\in\A^1$ act on $q\in U(L^*)$ sending it to $(1,a)q\in U(L^*)$, here $(1,a)\in R$ is written in Witt coordinates. For any $q\in U(L^*)$ the $\GL(L)\times\A^1$-orbit on $U(L^*)$ through $q$ is of dimension equal to that of $U(L^*)$. Since $U(L^*)$ is an open subscheme of an affine space, it is irreducible. So, the action of $\GL(L)\times \A^1$ on $U(L^*)$ is transitive.

 Now let $n$ be even, let $\tilde\phi,\phi,q$ be as in i). We have the corresponding flag $L_2\subset L_1\subset L$ on $L$. Pick $l\in L-L_1$ with $q(l) \mod 2=1$. Adding to $l$ a suitable element of $L_2$ we may assume that $q(l)=1$.
Let $\tilde l\in \tilde L$ be a vector over $l$. Pick any $\tilde l_2\in \tilde L$ such that its image $l_2$ in $L$ lies in $L_2$ and $\tilde\phi(\tilde l_2, \tilde l)=1$. We have in Witt coordinates $\tilde\phi(\tilde l_2,\tilde l_2)=(0,\epsilon)$ for some $\epsilon\in k$. The bilinear form $\tilde\phi$ is non degenerate over $R\tilde l_2\oplus R\tilde l$, so we get an orthogonal (with respect to $\tilde\phi$) decomposition of $\tilde L$ into a direct sum of free $R$-submodules $(R\tilde l_2\oplus R\tilde l)\oplus \tilde W$. Let $W$ be the image of $\tilde W$ in $L$ then $L_1=L_2\oplus W$. So, $q: W\to 2R$. 
 
  The bilinear form $\phi$ on $W$ is symplectic, this is the bilinear form associated to $q$ on $W$. So, $q$ can be seen as a non-degenerate quadratic form on $W$ with values in $k$. 
 It follows that the determinant of $\tilde\phi\mid_W$ equals $(-1)^{(n-2)/2}$. Indeed, there is a base $e_i, e_{-i}$ of $W$ in which the form $q$ on $W$ writes in Witt coordinates 
$$
q(x)=(0, \sum_i x_i^2x_{-i}^2)\in 2R
$$ 
for $x=\sum_i (x_i e_i+x_{-i}e_{-i})$. Actually, we have proved that in a suitable base of $L$ the form $q$ writes as 
\begin{equation}
\label{normal_form_q}
q(y,y_2, x_i, x_{-i})=(0,\epsilon)\tilde y_2^2+2\tilde y\tilde y_2+\tilde y^2+2\!\!\!\!\sum_{i=1}^{(n-2)/2} \tilde x_i \tilde x_{-i}
\end{equation}
in $n$ variables $\{y,y_2, x_i, x_{-i}\}$ with $1\le i\le (n-2)/2$. 
Here $\tilde y, \tilde y_2, \tilde x_i, \tilde x_{-i}$ are any liftings of the corresponding variables with respect to $R\to k$. In particular, for each $u$ the action of $\GL(L)$ on $U(L^*)_u$ is transitive.
  
  Since the determinant of $\tilde\phi$ over $R\tilde l_2\oplus R\tilde l$ equals $(1,1+\epsilon)$ in Witt coordinates, we get that the discriminant of such $q$ equals
$$
\left\{
\begin{array}{cl}
\epsilon, & n=0\mod 4\\
\epsilon+1, & n=2\mod 4
\end{array}
\right.
$$  
Our assertion follows. 
\end{Prf}

\medskip

\begin{Rem}
\label{Rem_open_orbits_of_GL(L)}
i) For $n$ even the normal form (\ref{normal_form_q}) can be seen as a section $s: \A^1\to U(L^*)$ of $\disc: U(L^*)\to \A^1$.  \\
\noindent 
ii) For any $q\in U(L^*)$ the stabilizor of $q$ in $\GL(L)$ has two connected components (cf. also Section~4.4). So, for any $u\in k$ the shifted local system $\bar S_{\psi}^{\otimes 2}$ over $U(L^*)_u$ is constant. It follows that there exists a local system $\cW_n$ on $\A^1$ together with an isomorphism over $Q^!_0(L^*)$
$$
\bar S_{\psi}^{\otimes 2}[-3n-n^2]\,\iso\, \disc^* \cW_n,
$$ 
here $\disc: Q^!_0(L^*)\to\A^1$. One checks that for any $n$ the local system $\cW_n$ on $\A^1$ is of order two, it becomes constant on the covering $\A^1\to \A^1$, $b\mapsto b^2+b$.
 
  Actually, for $n$ even we can say more. The restriction $s^*\bar S_{\psi}$ under the section $s: \A^1\to U(L^*)$ is of order 4 and can be calculated using the formula (\ref{normal_form_q}). Let $i_R: \A^1\to W_2$ be the map sending $x$ to $(x,0)$ in Witt coordinates. There exists an isomorphism over $\A^1$
$$
s^*\bar S_{\psi}[-\dim U(L^*)]\;\iso\; i_R^*\cL_{\psi}
$$
This follows easily from the fact that, at the level of functions, for $\epsilon, y\in k$ the sum 
$$
\sum_{y_2\in k} \psi((0,\epsilon)\tilde y_2^2+2\tilde y\tilde y_2+\tilde y^2)
$$ 
vanishes unless $\epsilon=y^4$.
\end{Rem}

\begin{Def} For $n=1$ and $L=k$ one has canonically $Q^!(L^*)\,\iso\, R$. In this case write $\cE\in \D(\Spec k)$ for the fibre of $\bar S_{\phi}[-2]$ at $1\in Q^!(L^*)$. Then $\cE$ is a 1-dimensional $\Qlb$-vector space placed in usual degree zero.  
\end{Def} 

\begin{Pp}
\label{Pp_power_4}
For any $n$ we have canonically over $Q^!_0(L^*)$
$$
\bar S_{\psi}^{\otimes 4}\mid_{Q^!_0(L^*)}[-6n-2n^2]\,\iso\, \cE^{\otimes 4n}
$$ 
\end{Pp}
\begin{Prf}
1) Assume $n$ even. Set $u=0\in k$. By Remark~\ref{Rem_open_orbits_of_GL(L)}, $\bar S_{\psi}^{\otimes 2}[-3n-n^2]$ is a constant local system on $Q^!_0(L^*)_u$. Choosing a base in $\tilde L$, we get a $k$-point in $Q^!_0(L^*)$, the image of the symmetric bilinear form $\tilde\phi\in\Sym^{!2}(\tilde L^*)$ with the identity matrix. The fibre of $\bar S_{\psi}^{\otimes 2}[-3n-n^2]$ at this $k$-point is $\cE^{2n}$. This yields an isomorphism
$$
\bar S_{\psi}^{\otimes 2}[-3n-n^2]\,\iso\, \cE^{2n}
$$
over $Q^!_0(L^*)_u$. \\
2) Considering the map $Q^!_0(L^*)\to Q^!_0(L^*\oplus L^*)_u$ sending $\phi$ to $\phi\oplus\phi: \tilde L\oplus \tilde L\to \tilde L^*\oplus \tilde L^*$, we are reduced to 1).
\end{Prf}

\medskip

 If $n=1$ then the local system $\bar S_{\psi}^{\otimes 2}\mid_{Q^!_0(L^*)}$ is nonconstant. Indeed, consider the closed immersion $\kappa_R: \AA^1\to W_2^*$ sending $x$ to $(1,x)\in R^*$ in Witt coordinates. Then $\RG_c(\AA^1, \kappa_R^*\bar S_{\psi}^{\otimes 2})$ is easy to calculate, which shows that $\kappa_R^*\bar S_{\psi}^{\otimes 2}$ is nonconstant. Actually, $\kappa_R^*\bar S_{\psi}^{\otimes 2}$ trivializes on the covering $\AA^1\to \AA^1$, $z\mapsto z^2+z$. This is a consequence of the following lemma.
 
\begin{Lm} 
\label{Lm_covering_case_n=1} 
Let $\bar i_R: \A^1\to W_2$ be the map sending $x$ to $(x+1,0)\in R$ in Witt coordinates. For $n=1$ one has an isomorphism
$$
\bar i_R^*\cL_{\psi}\,\iso\, \kappa_R^*(_1\bar S_{\psi}[-2])
$$
\end{Lm}
\begin{Prf}
Consider the quadratic form $q(x,y)=(1,a)\tilde x^2+\tilde y^2$ with $\tilde x,\tilde y\in R$ over $x,y\in k$. Consider a new base $\{u, u_2\}$ of $k^2$ given by $u=(0,1)$ and $u_2=(1,1)$. The value of $q$ at $yu+y_2u_2\in k^2$ is $$
(0,a+1)(y_2^2,0)+(0, y^2y_2^2)+(y^2,0)=(0,a+1)\tilde y_2^2+2\tilde y\tilde y_2+\tilde y^2
$$ 
So, the assertion follows from Remark~\ref{Rem_open_orbits_of_GL(L)}.
\end{Prf}
 
\begin{Def} 
\label{Def_torsor_Cov}
Fix once and for all an isomorphism $\cE^4\,\iso\,\Qlb$. Then $\bar S_{\psi}[-n-n(n+1)/2]$ is a local system on $Q^!_0(L^*)$ whose 4th power is trivialized canonically by Proposition~\ref{Pp_power_4}. In view of the homomorphism $\psi: \ZZ/4\ZZ\to \Qlb^*$ we have fixed in Section~3.1, $\bar S_{\psi}$ yields a $\ZZ/4\ZZ$-torsor $\Cov(Q^!_0(L^*))\to Q^!_0(L^*)$. The restriction of $\bar S_{\psi}$ to $\Cov(Q^!_0(L^*))$ is equipped with a canonical trivialization.
\end{Def}

\begin{Def}
\label{Def_Q_0_and_Cov}
Let $Q_0(\tilde L)\subset \Sym^{!2}(\tilde L^*)$ be the open subscheme of $\phi: \tilde L\to \tilde L^*$, which are isomorphisms of $R$-modules. Write $\Cov(Q_0(\tilde L))\to Q_0(\tilde L)$ for the $\ZZ/4\ZZ$-torsor obtained from $\Cov(Q^!_0(L^*))$ by the base change $Q_0(\tilde L)\to Q^!_0(L^*)$. 
\end{Def} 

\medskip\noindent
4.4 {\scshape Arf invariant of quadratic forms} 
This subsection is independent of the rest of the paper. 

\medskip\noindent
4.4.1 Let $S$ be a scheme and $E$ a vector bundle on $S$. 
A quadratic form on $E$ is a morphism of sheaves $q: E\to \cO_S$ such that for local sections $s\in\cO_S, e\in E$ one has $q(se)=s^2q(e)$ and the map $b_q: E\times E\to\cO_S$, $(e_1,e_2)\mapsto q(e_1+e_2)-q(e_1)-q(e_2)$ is $\cO_S$-bilinear. Then the Clifford algebra $C(E)$ is a sheaf of $\cO_S$-algebras on $S$ defined as the quotient of the tensor algebra $T(E)$ of $E$ by the two-sided sheaf of ideals generated by local sections of the form $e\otimes e-q(e)$. This is a $\ZZ/2\ZZ$-graded sheaf of $\cO_S$-algebras, write $C(E)^+$ and $C(E)^-$ for even and odd parts respectively.

 For example, if $q=0$ then $C(E)$ is $\ZZ$-graded and $C(E)\,\iso\, \oplus_{n\ge 0} \wedge^n(E)$ (the $\ZZ/2\ZZ$-gradation is given by the parity of $n$).

  Write $Z[C(E)^+]$ for the center of the sheaf of $\cO_S$-algebras $C(E)^+$. By (\cite{B}, Theorem~3.6, p. 39), 
if $E$ is of rank $2n$ for some $n$ and $b_q$ is non-degenerate then $Z[C(E)^+]$ is a sheaf of quadratic separable algebras over $S$, that is, $\Spec( Z[C(E)^+])$ is an \'etale two-sheeted covering of $S$. 
This is \select{the Arf invariant} of the quadratic form $q: E\to \cO_S$. 
 
\medskip\noindent
4.4.2 Assume that $S=\Spec k$ with $k$ algebraically closed of characteristic two. Let $E$ be a $k$-vector space of dimension $2n$. A quadratic form $q: E\to k$ is called \select{non-degenerate} if the bilinear form $b_q(x,y)=q(x+y)-q(x)-q(y)$ on $E$ is non-degenerate. The group $\GL(E)$ acts transitively on the variety of non-degenerate quadratic forms on $E$ (\cite{B2}), and the bilinear form $b_q$ is symplectic. Let $q:E\to k$ be a non-degenerate bilinear form on $E$. 

\begin{Lm} The stabilizer $\OO(q)$ of $q$ in $\GL(E)$ has two connected components. We write $\SO(q)$ for the connected component of unity, although $\OO(q)\subset \Sp(b_q)\subset \SL(E)$. The group $\OO(q)$ acts naturally on $C(E)$, so on $Z[C(E)^+]$. The action of $\SO(q)$ on $Z[C(E)^+]$ is trivial, and the group $\OO(q)/\SO(q)$ induces a unique nontrivial $k$-automorphism of the $k$-algebra $Z[C(E)^+]$.
\end{Lm}
\begin{Prf}  This follows from (\cite{Bou}, exercice 9 on p. 155). Namely, if $\{e_i, e_{-i}\}$ ($1\le i\le n$) is a symplectic base in $E$ then $\{1,z\}$ is a $k$-base of $Z[C(E)^+]$, where $z=\sum e_i e_{-i}$. If $W$ is the Weyl group of $\OO(q)$ for the corresponding torus, we have an exact sequence $1\to S_n\to W\to (\ZZ/2\ZZ)^n\to 1$. The image of $\sigma\in S_n$ in $W$ sends $e_i$ to $e_{\sigma i}$ and $e_{-i}$ to $e_{-\sigma i}$ ($1\le i\le n$). The group $W$ contains also the elements $w_i$ permuting $e_i$ and $e_{-i}$. Let $W^+$ be the subgroup of those $w\in W$ whose image $a$ in  $(\ZZ/2\ZZ)^n$ satisfies $\sum_{i=1}^n a_i=0$. Let $W^-=W-W^+$. Then 
$$
\SO(q)=\cup_{w\in W^+} BwB\;\;\;\;\mbox{and}\;\;\;\;
\OO(q)-\SO(q)=\cup_{w\in W^-} BwB  
$$
Since any $w\in W^-$ acts on the base $\{1,z\}$ as $\{1, z+1\}$, we are done.
\end{Prf}

\medskip

 Now if $S$ is a $k$-scheme, let $q: \cE\to\cO_S$ be a vector bundle of rank $2n$ with a non-degenerate quadratic form. This is nothing but a torsor for $\OO(q_0)$, where $q_0$ is the split quadratic form of rank $2n$. So it induces via extension of scalars $\OO(q_0)\to \OO(q_0)/\SO(q_0)$ a $\ZZ/2\ZZ$-torsor, which identifies canonically with $\Spec Z[C(\cE)^+]$. Another way is to say that $\GL(E)/\OO(q_0)$ is the stack classifying non-degenerate quadratic forms of dimension $2n$ , and $(\cE,q)$ is a datum of a map $\tau: S\to \GL(E)/\OO(q_0)$. Then the Arf invariant of $(\cE,q)$ is the restriction under $\tau$ of the $\ZZ/2\ZZ$-torsor $\GL(E)/\SO(q_0)\to \GL(E)/\OO(q_0)$. 

\medskip\noindent
4.5 {\scshape Description of $\Cov(Q^!_0(L^*))$ generically}

\medskip\noindent
4.5.1 Use the notations of Section 4.3. Let $q\in Q^!(L^*)$, so $q: L\to R$. Say that a 1-dimensional $k$-subspace $E\subset L$ is \select{multiplicative} if for any nonzero $e\in E$ we have in Witt coordinates $q(e)=(a,0)$ for some $a\in k$. This condition does not depend on a choice of $e$. 

 For a point $q\in U(L^*)$ let $r$ be the image of $\Fr(q)$  in $Q_a(L^*)$, where $\Fr: U(L^*)\to U(L^*)$ is the Frobenius map. Let the hyperplane $L_1\subset L$ be the kernel of $r$.  As $q$ varies these hyperplanes form a vector bundle over $U(L^*)$.   
 
 If $n$ is odd, we consider the vector bundle $\cE_1$ over $U(L^*)$ whose fibre at $q$ is $L_1$ for the corresponding $r$. It is equipped with the quadratic form $q: L_1\to 2R$, which organize into a quadratic form on $\cE_1$. As in Lemma~\ref{Lm_stabilizors}, the bilinear form $(x,y)\mapsto q(x+y)-q(x)-q(y)$ on $L_1$ is symplectic. The corresponding Clifford algebra over $U(L^*)$ is denoted $C(\cE_1)$. Write $\tilde U(L^*)$ for $\Spec Z[C(\cE_1)]$, this is a $\ZZ/2\ZZ$-Galois \'etale covering over $U(L^*)$. The group $\GL(L)$ acts naturally on $\tilde U(L^*)$, the projection $\eta: \tilde U(L^*)\to U(L^*)$ is $\GL(L)$-equivariant.

 By Remark~\ref{Rem_open_orbits_of_GL(L)}, the local system $\eta^*\bar S_{\psi}[-\dim U(L^*)]$ descends with respect to the composition 
$
\tilde U(L^*)\toup{\eta} U(L^*)\toup{\disc}\A^1
$.  
More precisely, there exists an isomorphism over $\tilde U(L^*)$
$$
\eta^*\bar S_{\psi}[-\dim U(L^*)]\,\iso\, \eta^*\disc^*\kappa_R^*(_{1}\bar S_{\psi}[-2]),
$$
where $\kappa_R$ is the map from Lemma~\ref{Lm_covering_case_n=1}.
     
\medskip\noindent
4.5.2 Now consider the case of $n$ even. For $q\in U(L^*)$ write $r\in Q_a(L^*)$ and $\phi\in \Sym^{!2}(L^*)_0$ for its images in $Q_a(L^*)$ and $\Sym^{!2}(L^*)_0$ respectively. We also have the corresponding flag $L_2\subset L_1\subset L$, where $L_1$ is the kernel of $r$, and $L_2=L_1^{\perp}$ with respect to $\phi$.

 Assume that $q$ is nonzero on $L_2$, then there is a unique subspace $W\subset L_1$ such that $L_2\oplus W=L_1$ and $W$ is preserved by the stabilizer $\OO(q)$ of $q$ in $\GL(L)$. The bilinear form $\phi\mid_W$ is symplectic. Let $W^{\perp}\subset L$ be its orthogonal complement in $L$, so $L_2\subset W^{\perp}$ and $W^{\perp}/L_2\,\iso\, L/L_1$ naturally.    
     
  Then the form $q\mid_{W^{\perp}}$ lies in $U((W^{\perp})^*)$, that is, it has the same properties as $q$ itself but now for dimension two. One checks that there are exactly two 1-dimensional multiplicative subspaces in $W^{\perp}$ for $q$. The stabilizer $\OO(q)\subset \GL(L)$ of $q$ in $\GL(L)$ acts transitively on these two multiplicative subspaces. 
   
    Let $(\A^1)^0\subset \A^1$ be the open subscheme given by the condition that $a\in \A^1$ is not the critical value given by (\ref{disc_critical_value}).  Let $U(L^*)^0\subset U(L^*)$ be the preimage of $(\A^1)^0$ under $\disc: U(L^*)\to \A^1$. 
   
   We get a 2-sheeted \'etale covering 
$$
\eta: \tilde U(L^*)^0\to U(L^*)^0,
$$
which classifies a point $q\in U(L^*)^0$ together with a one-dimensional multiplicative subspace $E\subset W^{\perp}$. The section $s: \A^1\to U(L^*)$ from Remark~\ref{Rem_open_orbits_of_GL(L)} extends to a section $\tilde s: (\A^1)^0\to \tilde U(L^*)^0$, namely, in the notation of formula (\ref{normal_form_q}) we have a distinguished multiplicative subspace given by $y_2=x_i=x_{-i}=0$.

 By Lemma~\ref{Lm_stabilizors}, $\eta^*\bar S_{\psi}$ descends with respect to the composition 
$$
\tilde U(L^*)^0\toup{\eta} U(L^*)^0\toup{\disc}\A^1
$$  
More precisely, there is an isomorphism over $\tilde U(L^*)^0$
$$
\eta^*\bar S_{\psi}\,\iso\, \eta^*\disc^* s^*(\bar S_{\psi})
$$
Remind that the shifted local system $s^*(\bar S_{\psi})$ is of rank one and order four, it is described in Remark~\ref{Rem_open_orbits_of_GL(L)}.

\medskip\noindent
4.6 {\scshape Description of $\bar S_{\psi}$ on strata}

\medskip\noindent
 Write $\Sym^{!2}_i(L^*)\subset \Sym^{!2}(L^*)$ for the locus of symmetric bilinear forms on $L$ whose kernel is of dimension $i$. Let $Q^!_i(L^*)$ be the preimage of $\Sym^{!2}_i(L^*)$ under the projection $Q^!(L^*)\to \Sym^{!2}(L^*)$. Write $'Q_i(L^*)\subset Q^!_i(L^*)$ for the closed subscheme classifying those $q\in Q^!_i(L^*)$ which vanish on the kernel of $\phi$, where $\phi\in \Sym^{!2}_i(L^*)$ is the image of $q$. Then $'Q_i(L^*)$ can be seen as a scheme classifying an $i$-dimensional subspace $M\subset L$ and $q\in Q^!_0((L/M)^*)$. The  corresponding map $L\to R$ is then the composition $L\to L/M\toup{q} R$.

\begin{Pp} 
\label{Pp_restriction_strata}
For any $i$ the $*$-restriction of $\bar S_{\psi}$ to $Q^!_i(L^*)$ is the extension by zero from $'Q_i(L^*)$. The corresponding sheaf 
$$
\bar S_{\psi}\mid_{'Q_i(L^*)} [i-n-\frac{n(n+1)}{2}]
$$
is a local system of rank one and order four (except for the last stratum $'Q_n(L^*)=\Spec k$, over which it is a trivial local system of rank one). If $M\subset L$ is a subspace and $q\in Q^!_0((L/M)^*)$ then the fibre of $\bar S_{\psi}$ at $L\to L/M\toup{q}R$ is the shifted Gauss sum for $q$.
\end{Pp}
\begin{Prf}
Let $q\in Q^!_i(L^*)$ and $\phi$ be its image in $\Sym^{!2}_i(L^*)$. Let $M\subset L$ be the kernel of $\phi$. Then for $x\in L, y\in M$ we have $q(x+y)=q(x)+q(y)$. At the level of functions, for any $x\in L$ the sum $\sum_{y\in M} \psi(q(x+y))$ will vanish unless $q\in {'Q_i(L^*)}$. The geometrization is straightforward, our assertion easily follows.
\end{Prf}

\bigskip

\centerline{\scshape 5. Maslov index in characteristic two}

\bigskip\noindent
In Sections~5.1-5.3 we generalize the results of \cite{TT} to the case of characteristic two. 

\medskip\noindent
5.1 Let $\cO$ be a discrete valuation ring, 
$I\subset \cO$ be a non zero ideal and $\bar\cO=\cO/I$.
Let $\bar V$ be a free $\bar\cO$-module of rank $2n$ with a symplectic form $\omega: \bar V\otimes\bar V\to\bar\cO$. For a submodule $M\subset \bar V$ let $M^{\perp}=\{x\in \bar V\mid \omega(x,m)=0\; \mbox{for all}\; m\in M\}$. One checks that $(M^{\perp})^{\perp}=M$. 

 For a $\bar \cO$-module $M$ set $M^*=\Hom_{\bar\cO}(M,\bar \cO)$ in the category of $\bar\cO$-modules. The functor $M\mapsto M^*$ is exact on the category of finite type $\bar\cO$-modules. For any $\bar\cO$-module of finite type $M$ one has $(M^*)^*\,\iso\, M$ canonically. If $K$ is a bounded complex of free $\bar\cO$-modules of finite type then $\H^0(\Hom_{\bar\cO}(K,\bar\cO))\,\iso\, \Hom_{\bar\cO}(\H^0(K), \bar\cO)$ canonically.
 
 For any submodule $M\subset \bar V$ we have an isomorphism of exact sequences of $\bar\cO$-modules
$$
\begin{array}{ccccccc}
0\to & M & \to & \bar V& \to &\bar V/M&\to 0\\ 
& \downarrow && \downarrow&& \downarrow\\
0\to &(\bar V/M^{\perp})^* & \to & \bar V^* & \to & (M^{\perp})^* & \to 0
\end{array}
$$  

\noindent
5.2 Keep notation of Section~3. Let $\tilde L$ be a free $R$-module of finite type with symmetric bilinear form $\phi: \tilde L\to\tilde L^*$. We say that $\phi$ is \select{hyperbolic} if there is a base $\{e_i, e_{-i}\}$ of $\tilde L$ in which the form is given by $\phi(e_i, e_j)=0$ unless $i=-j$, and $\phi(e_i, e_{-i})=
\phi(e_{-i}, e_i)=1$. 

\begin{Lm} If $\phi$ is an isomorphism then   \\
1) for any free submodule $\tilde N\subset\tilde L$ its orthogonal complement $\tilde N^{\perp}\subset\tilde L$ is free. \\
2) If $\tilde N$ is a free isotropic submodule in $\tilde L$ then the induced symmetric bilinear form $\bar\phi$ on $\tilde N^{\perp}/\tilde N$ is non degenerate, and we have an isometry $(\tilde N^{\perp}/\tilde N)\oplus (\tilde N\oplus\tilde M)\,\iso\, \tilde L$ for any free submodule $\tilde M\subset \tilde L$ such that $\tilde N^{\perp}\oplus \tilde M=\tilde L$. Here $\tilde N\oplus\tilde M$ is equipped with the form induced by $\phi$ (and $\tilde N\oplus\tilde M$ is not necessarily hyperbolic).
\end{Lm}
\begin{Prf}
2) We have an orthogonal decompositon $(\tilde N\oplus\tilde M)\oplus (\tilde N\oplus\tilde M)^{\perp}\,\iso\,\tilde L$, and $(\tilde N\oplus\tilde M)^{\perp}\subset \tilde N^{\perp}$ maps isometrically onto $\tilde N^{\perp}/\tilde N^{\perp}$.
\end{Prf}   

\medskip

\begin{Rem} 
\label{Rem_free_isotropic_submodule_and_Gauss}
If $\phi: \tilde L\to\tilde L^*$ is non degenerate of rank $2r$, and $\tilde N\subset\tilde L$ is a free isotropic submodule of rank $r$ (it is automatically a direct summand), then $(\tilde L,\phi)$ is not always hyperbolic. At the level of matrices the reason is that one can not in general present a given symmetric matrix $B\in\Mat_r(R)$ as $B=A+{^tA}$ for another matrix $A\in\Mat_r(R)$. However, the Gauss sum of such $(\tilde L,\phi)$ is canonically trivialized.
\end{Rem}

\medskip\noindent
5.3 Let $\tilde V$ be a free $R$-module of rank $2n$ with symplectic form $\tilde\omega: \tilde V\times\tilde V\to R$. Write $\cL(\tilde V)$ for the variety of free lagrangian $R$-submodules in $\tilde V$. We view it as a $k$-scheme via its Greenberg realization.  

  For any submodule $\tilde M\subset\tilde V$ we have $(\tilde M^{\perp})^{\perp}=\tilde M$. If $\tilde M\subset\tilde V$ is a free isotropic submodule then $\tilde M^{\perp}$ is free, and $\tilde M^{\perp}/\tilde M$ is equipped with a symplectic form with values in $R$. For a family of submodules $\tilde M_i\subset \tilde V$ we have $\cap_i (\tilde M_i^{\perp})=(\sum_i \tilde M_i)^{\perp}$.
  
  Let $m\ge 2$. As in \cite{TT}, think of $\ZZ/m\ZZ$ as the vertices of a graph whose set $\EE$ of edges is the set of pairs of consecutive numbers $\{i,i+1\}$, $i\in\ZZ/m\ZZ$. For 
$$
v=(v_{\{i,i+1\}})\in \mathop{\oplus}\limits_{\{i,i+1\}\in\EE} \tilde V
$$
its \select{derivative} is 
$$
\partial v=(\partial v_i)\in \mathop{\oplus}\limits_{i\in \ZZ/m\ZZ} V,\;\;\;\;\;\;\;\;\;\;\;\;\; \partial v_i=v_{\{i,i+1\}}-v_{\{i-1,i\}}
$$
Conversely, for $w=(w_i)\in \mathop{\oplus}\limits_{i\in\ZZ/m\ZZ} V$ its \select{antiderivative} is
$$
\hat w=(\hat w_{\{i, i+1\}})\in\mathop{\oplus}\limits_{\{i,i+1\}\in\EE} V\;\;\;\;\mbox{such that}\;\;\;\;\; \partial(\hat w)=w
$$
An anti-derivative exists if $\sum_{i\in\ZZ/m\ZZ} w_i=0$ and is unique up adding a constant function.
 
\begin{Def} 
\label{Def_K_many}
For $\tilde L_1,\ldots,\tilde L_m\in\cL(\tilde V)$ let $K_{1,2,\ldots,m}$ be the kernel of
$$
\mathop{\oplus}\limits_{i\in\ZZ/m\ZZ} \tilde L_i\;\toup{\sum}\;\tilde V
$$
The surjectivity of the above $\Sigma$ is equivalent to requiring that $\cap_i \tilde L_i=0$, and in this case $K_{1,\ldots,m}$ is a free $R$-module.
\end{Def}
 
 As in \cite{TT}, $K_{1,\ldots,m}$ is equipped with the $R$-valued symmetric bilinear form 
$$ 
q_{1,2,\ldots,m}(v,w)=\sum_{i\in\ZZ/m\ZZ} \tilde\omega(v_i, \hat\omega_{\{i,i+1\}})
$$
for any anti-derivative $\hat w\in \oplus_{\{i,i+1\}\in\EE} V$ of $w$. The derivative restricts to a map
$$
\partial: \mathop{\oplus}_{\{i,i+1\}\in\EE} \tilde L_i\cap \tilde L_{i+1}\to K_{1,\ldots,m},
$$
whose kernel is $\cap_i \tilde L_i$, and its image is contained in the kernel of $q_{1,2,\ldots,m}$. Set 
$$
T_{1,2,\ldots,m}=K_{1,\ldots,m}/\Im\partial
$$
This is an $R$-module, which is not necessarily free, it is equipped with the induced symmetric bilinear form $q_{1,2,\ldots,m}$. For a $R$-module of finite type $M$ write $M^*=\Hom_R(M,R)$.
 
\begin{Lm} The form $q_{1,\ldots,m}$ induces an isomorphism of $R$-modules $T_{1,\ldots,m}\,\iso\, T_{1,\ldots,m}^*$.
\end{Lm}
\begin{Prf}
Consider the complex placed in degrees $-1,0,1$
\begin{equation}
\label{complex_Thomas_dual}
V^*\toup{\diag} \mathop{\oplus}_{i\in \ZZ/n\ZZ} \tilde L_i^*\toup{\partial^*} \mathop{\oplus}_{\{i,i+1\}\in\EE} (\tilde L_i\cap \tilde L_{i+1})^*,
\end{equation}
where $\partial^*$ is the transpose to $\partial: \oplus_{\{i,i+1\}\in\EE} (\tilde L_i\cap \tilde L_{i+1})\to \oplus_{i\in\ZZ/m\ZZ} \tilde L_i$. The proof of (\cite{TT}, Proposition~3) goes through in this situation and yields in an isomorphism between $K_{1,\ldots,m}/\Im \partial$ and $\H^0$ of (\ref{complex_Thomas_dual}). The cokernel of $\partial^*$ is $(\cap_i \tilde L_i)^*$. To finish the proof use 5.1 for $\bar\cO=R$.
\end{Prf} 
 
\medskip

 Assume that $\cap_i \tilde L_i=0$. Then $T_{1,\ldots,m}$ is free iff 
for any $\{i,i+1\}\in\EE$, $\tilde L_i\cap\tilde L_{i+1}$ is free (this is also equivalent to $\tilde L_i+\tilde L_{i+1}$ being free). A straightforward analog of (\cite{TT}, Proposition~5) yields the following.

\begin{Lm}
\label{Lm_dihedral}
 For any $\tilde L_i\in \cL(\tilde V)$ there are 
canonical isometries $T_{1,\ldots, m}\,\iso\, T_{2,3,\ldots,m,1}$ and 
$$
(T_{1,\ldots, m}, q_{1,\ldots,m})\,\iso\, (T_{m,\ldots, 1}, -q_{m,\ldots,1})
$$
One can also replace $T$ by $K$ in the both above isomorphisms. \QED
\end{Lm} 

If $N$ is an $R$-module of finite type with a non degenerate symmetric bilinear form $\phi: N\,\iso\, N^*$ and $M\subset N$ is an isotropic submodule, let $M^{\perp}=\{n\in N\mid \<n,\phi(m)\>=0\;\mbox{for all}\; m\in M\}$. Then the $R$-module $M^{\perp}/M$ is equipped with a non degenerate symmetric bilinear form, and we say that $M^{\perp}/M$ is \select{a quadratic subquotient} of $N$. 
The proof of (\cite{TT}, Proposition~6) applies without changes in our situation and yields the following.

\begin{Pp} 
\label{Pp_chain_condition}
Let $\tilde L_1,\ldots,\tilde L_m\in \cL(\tilde V)$ and $k\in \{2,\ldots,m\}$. \\  
1) If $\tilde L_1\cap\tilde L_k=0$ then there is a canonical isometry $K_{1,2,\ldots,k}\oplus K_{1,k,\ldots,m} \,\iso\, K_{1,2,\ldots,m}$ and also 
$$
T_{1,2,\ldots,k}\oplus T_{1,k,\ldots,m} \,\iso\, T_{1,2,\ldots,m}
$$ 
2) Without conditions, $T_{1,2,\ldots,k}\oplus T_{1,k,\ldots,m}$ is a quadratic subquotient of $T_{1,2,\ldots,m}$.
\QED
\end{Pp} 

\noindent
5.4 {\scshape Maslov index in families}

\medskip\noindent
Let $\tilde L$ be a free $R$-module of rank $n$. Let $U_{23}$ be the variety of triples $\tilde L_1, \tilde L_2, \tilde L_3\in \cL(\tilde V)$ such that $\tilde L_1\cap \tilde L_2=\tilde L_1\cap\tilde L_3=0$. Define a morphism of $k$-stacks
\begin{equation}
\label{map_pi_23}
\pi_U: U_{23}\to \Sym^{!2}(\tilde L^*)/\GL(\tilde L)
\end{equation}
as follows (it is understood that we first take the Greenberg realizations of $\Sym^{!2}(\tilde L^*)$ and of $\GL(\tilde L)$ and then the stack quotient of one by the other).

 Given $(\tilde L_1,\tilde L_2,\tilde L_3)\in U_{23}$
let $r: \tilde L_2\to \tilde L_1$ be the $R$-linear map such that 
$$
\tilde L_3=\{r(x)-x\in \tilde L_1\oplus \tilde L_2\mid x\in \tilde L_2\}
$$ 
Let $\tilde\phi\in \Sym^{!2}(\tilde L_2^*)$ be given by $\tilde\phi(\tilde x, \tilde y)=\tilde\omega(r(\tilde x), \tilde y)$ for $\tilde x,\tilde y\in \tilde L_2$. 

\begin{Lm} For $(\tilde L_1,\tilde L_2,\tilde L_3)\in U_{23}$ there is a canonical isometry $(\tilde L_2,\tilde\phi)\,\iso\, K_{1,2,3}$. It identifies $\tilde L_2\cap\tilde L_3$ with the kernel of $(K_{1,2,3}, \; q_{1,2,3})$.
\end{Lm}
\begin{Prf}
The map $\tilde L_2\to K_{1,2,3}$ given by $x\mapsto (-r(x), x, r(x)-x)$ is the desired isometry.
\end{Prf}

\medskip

We let $\pi_U$ send $(\tilde L_1, \tilde L_2, \tilde L_3)$ to $(\tilde L_2, \tilde\phi)$. 
Note that the variety $\{\tilde N\in \cL(\tilde V)\mid \tilde N\cap\tilde L_1=0\}$ admits a free transitive action of the $R$-module $\Sym^{!2}(\tilde L_1)$.
 
 Write $U_3\subset U_{23}$ for the open subscheme of triples $(\tilde L_1,\tilde L_2,\tilde L_3)\in U_{23}$ such that $\tilde L_i\cap \tilde L_j=0$ for $i\ne j$. We also denote by the same symbol
\begin{equation}
\label{map_pi_3} 
 \pi_U: U_3\to Q_0(\tilde L)/\GL(\tilde L)
\end{equation} 
the restriction of $\pi_U$. Let $\Cov(U_3)\to U_3$ be the $\ZZ/4\ZZ$-torsor obtained from 
$$
\Cov(Q_0(\tilde L))/\GL(\tilde L)\to Q_0(\tilde L)/\GL(\tilde L)
$$ 
by the base change $\pi_U$. 
 
 Given a $k$-scheme $S$ and a $S$-point $(\tilde L_1,\tilde L_2,\tilde L_3)\in U_3$, write $\cC(\tilde L_1,\tilde L_2,\tilde L_3)$ for the $\ZZ/4\ZZ$-torsor over $S$ obtained by restricting $\Cov(U_3)$ under $S\to U_3$. This notation agrees with (\ref{Gauss_sum_first_classical}) and \cite{GH}.
 
\begin{Pp} 
\label{Con_Maslov_is_cocycle}
There is a canonical isomorphism
$$
\cC(\tilde L_1,\tilde L_2,\tilde L_3) \otimes \cC(\tilde L_0, \tilde L_1, \tilde L_3)=\cC(\tilde L_0,\tilde L_2,\tilde L_3)\otimes \cC(\tilde L_0,\tilde L_1,\tilde L_2)
$$
of $\ZZ/4\ZZ$-torsors on the variety $U_4$ classifying $\tilde L_0,\ldots,\tilde L_3\in \cL(\tilde V)$, which are paiwise transverse.
\end{Pp}
\begin{Prf}
The Gauss sum takes a direct sum to the tensor product of $\ZZ/4\ZZ$-torsors. So, it suffices to show that, given $\tilde L_0,\ldots,\tilde L_3\in U_4$, there is a canonical isometry of $R$-modules with the corresponding symmetric bilinear forms
$$
T_{1,2,3}\oplus T_{0,1,3}\,\iso\, T_{0,1,2}\oplus T_{0,2,3}
$$
And these isometries naturally organize into a family over $U_4$. Indeed, by Lemma~\ref{Lm_dihedral} and Proposition~\ref{Pp_chain_condition}, one has canonical isometries
$$
T_{1,2,3}\oplus T_{1,3,0}\,\iso\, T_{1,2,3,0}\,\iso\, T_{0,1,2,3}\,\iso\, T_{0,1,2}\oplus T_{0,2,3}
$$
To finish, use the dihedral isometry $T_{1,3,0}\,\iso\, T_{0,1,3}$ of Lemma~\ref{Lm_dihedral}.
\end{Prf} 
 
\medskip
\begin{Rem} 
\label{Rem_comparison_two_q_forms}
Let $(\tilde L_1,\tilde L_2, \tilde L_3)\in U_{23}$ and $L_i=\epsilon(\tilde L_i)$ be the corresponding enhanced lagrangians. Then the image of $\pi_U(\tilde L_1,\tilde L_2, \tilde L_3)$ in $Q^!(L_2^*)$ is the form $Q_{L_1,L_2,L_3}: L_2\to R$ defined by (\ref{eq_def_for_Q_123}).
\end{Rem}  

\medskip\noindent
5.5 {\scshape New isometries for the Maslov index}

\medskip\noindent
The following result will be used in Section~6.3.2.

\begin{Lm} 
\label{Lm_n-1_inclusion}
Let $\tilde L_1, \ldots,\tilde L_m\in \cL(\tilde V)$. There is a canonical isometric inclusion $K_{1,\ldots, m-1}\hook{} K_{1,\ldots,m}$.
\end{Lm}
\begin{Prf} Write $v\in K_{1,\ldots, m}$ as 
$v=(v_1,\ldots, v_m)$ with $\sum_i v_i=0$, $v_i\in \tilde L_i$. Then the desired isometry is given by the map $(v_1,\ldots, v_{m-1})\mapsto (v_1,\ldots, v_{m-1},0)$.
\end{Prf} 
 
\begin{Pp} 
\label{Pp_strange_isometry}
Let $(\tilde N, \tilde N', \tilde N'')$ and $(\tilde L,\tilde M)$ be two collections of lagrangians  in $\cL(\tilde V)$. Assume that each lagrangian from the first collection is transverse to each lagrangian from the second one. Then 
\begin{equation}
\label{isometry_new}
K_{\tilde N'', \tilde M, \tilde N', \tilde L}\oplus K_{\tilde N', \tilde M, \tilde N, \tilde L}
\end{equation}
admits a free isotropic $R$-submodule $D$ such that one has a canonical isometry $D^{\perp}/D\,\iso\, K_{\tilde N'', \tilde M, \tilde N, \tilde L}$.
\end{Pp}
\begin{Prf}
By Lemmas~\ref{Lm_dihedral} and \ref{Lm_n-1_inclusion}, (\ref{isometry_new}) is canonically isometric to $K_{\tilde M, \tilde N', \tilde L, \tilde N''}\oplus K_{\tilde L, \tilde N', \tilde M, \tilde N}$ and admits canonically the isometric
subspace $K_{\tilde M, \tilde N', \tilde L}\oplus K_{\tilde L, \tilde N', \tilde M}$. By Lemma~\ref{Lm_dihedral}, one has an isometry
$$
\sigma: (K_{\tilde M, \tilde N', \tilde L}, q_{\tilde M, \tilde N',\tilde L})\,\iso\, (K_{\tilde L, \tilde N', \tilde M}, -q_{\tilde L, \tilde N',\tilde M})
$$ 
Then $D=\{x+\sigma(x)\in K_{\tilde M, \tilde N', \tilde L}\oplus K_{\tilde L, \tilde N', \tilde M}\mid x\in K_{\tilde M, \tilde N', \tilde L}\}$ is the desired free isotropic $R$-module of rank $n$. Indeed, 
one checks that $D^{\perp}\subset K_{\tilde M, \tilde N', \tilde L, \tilde N''}\oplus K_{\tilde L, \tilde N', \tilde M, \tilde N}$ is the submodule
\begin{equation}
\label{eq_D_perp}
\{(m_1, n'_1, l_1, n''_1)\in K_{\tilde M, \tilde N', \tilde L, \tilde N''}, \; (l_2, n'_2, m_2, n_2)\in  K_{\tilde L, \tilde N', \tilde M, \tilde N}\;\mid n'_1=n'_2\}
\end{equation}
and $D=\{(m_1, n'_1, l_1, 0), (l_1, n'_1, m_1, 0)\;\mid (m_1, n'_1, l_1)\in K_{\tilde M, \tilde N', \tilde L}\}\subset D^{\perp}$.
The map sending a collection (\ref{eq_D_perp}) to 
$$
(n''_1, m_1-m_2, -n_2, l_1-l_2)\in K_{\tilde N'', \tilde M, \tilde N, \tilde L}
$$ 
yields the desired isometry.
\end{Prf}

\bigskip

\centerline{\scshape 6. Construction of the gerb $\hat\cL(\tilde V)$ over $\cL(\tilde V)$}
 
\bigskip\noindent
6.1  Use notation of Section~5. Our aim now is to construct a $\ZZ/4\ZZ$-gerb 
$\hat\cL(\tilde V)\to \cL(\tilde V)$ (cf. Appendix~A for a definition of a $\ZZ/4\ZZ$-gerb). 

 For $\tilde L\in \cL(\tilde V)$ write $\cU(\tilde L)=\{\tilde M\in \cL(\tilde V)\mid \tilde L\cap\tilde M=0\}$. Note that $\cU(\tilde L)$ depends only on the image of $\tilde L$ under $\cL(\tilde V)\to \cL(V)$. We want to construct $\hat\cL(\tilde V)$ by gluing the trivial gerbs $\cU(\tilde L)\times B(\ZZ/4\ZZ)$ over the open subschemes $\cU(\tilde L)$.

\begin{Def} For $m\ge 1$ say that a finite family of lagrangians $\tilde L_i\in \cL(\tilde V)$, $i\in I$ is \select{$m$-sweeping} if they are pairwise transverse and 
$$
\cup_{i\in I} \; \cU(\tilde L_i)^m=\cL(\tilde V)^m,
$$ 
here $Y^m$ denotes the $m$-th cartesian power of a scheme $Y$ over $k$. Since $\cL(\tilde V)^m$ is quasi-compact, for any $m$ a $m$-sweeping family exists.
\end{Def}
 
 Pick a 1-sweeping family $\tilde L_i\in \cL(\tilde V)$, $i\in I$.
The $\ZZ/4\ZZ$-torsors $\cC(\tilde L_i,\tilde L_j, \tilde L_k)$ (with $i,j,k\in I$) over $\Spec k$ satisfy the cocycle condition given by Proposition~\ref{Con_Maslov_is_cocycle}. So, we may and pick $\ZZ/4\ZZ$-torsors $\cB_{ij}$ over $\Spec k$ together with isomorphisms $\cB_{ji}\,\iso\, \cB_{ij}^{-1}$ and
\begin{equation}
\label{coboundary_for_cC}
\cC(\tilde L_i, \tilde L_j,\tilde L_k)^{-1}\,\iso\, \cB_{ij}\otimes\cB_{jk}\otimes\cB_{ki}
\end{equation}
for all $i,j,k$. We call \select{an enriched sweeping family} a collection $\{\tilde L_i\}_{i\in I}$ together with $\cB_{ij}$ and isomorphisms (\ref{coboundary_for_cC}).

 Set $\cU_i=\cU(\tilde L_i)$ and 
$$
\cU_{ij}=\cU_i\cap \cU_j, \;\;\;\; \cU_{ijk}=\cU_i\cap\cU_j\cap \cU_k
$$ 
 Let $\Cov(\cU_{ij})$ be the $\ZZ/4\ZZ$-torsor over $\cU_{ij}$ obtained from $\Cov(U_3)$ by the base change $\cU_{ij}\to U_3$, $\tilde L\mapsto (\tilde L, \tilde L_i,\tilde L_j)$. The gluing data consists of
\begin{itemize}
\item the $\ZZ/4\ZZ$-torsor $\Cov(\cU_{ij})\otimes\cB_{ij}$ over $\cU_{ij}$ giving rise to the automorphism
$$
\sigma_{ji}: \cU_{ij}\times B(\ZZ/4\ZZ)\,\iso\, \cU_{ij}\times B(\ZZ/4\ZZ)
$$
sending $(\tilde L, \cF)$ to $(\tilde L, \cF\otimes\Cov(\cU_{ij})\otimes\cB_{ij})$. To be presice, for a scheme $Z$ it sends a $Z$-point $(f,\cF)$ on the left to the $Z$-point $(f, \cF\otimes f^*\Cov(\cU_{ij})\otimes\cB_{ij})$ on the right. Here $f: Z\to \cU_{ij}$ is a morphism, and $\cF$ is a $\ZZ/4\ZZ$-torsor on $Z$. 
\item a 2-morphism $\delta_{kji}: \sigma_{kj}\comp \sigma_{ji}\,\iso\, \sigma_{ki}$ making the following diagram 2-commutative
$$
\begin{array}{ccc}
\cU_{ijk}\times B(\ZZ/4\ZZ) & \toup{\sigma_{ji}} & \cU_{ijk}\times B(\ZZ/4\ZZ)\\
& \searrow\lefteqn{\scriptstyle \sigma_{ki}} & \downarrow\lefteqn{\scriptstyle\sigma_{kj}}\\
&&\cU_{ijk}\times B(\ZZ/4\ZZ) 
\end{array}
$$
It is given by the $\ZZ/4\ZZ$-torsor $\cC(\tilde L_i, \tilde L_j,\tilde L_k)$ over $\Spec k$ together with the isomorphism of $\ZZ/4\ZZ$-torsors over $\cU_{ijk}$
\begin{equation}
\label{eq_for_cC_cocycle}
\cC(\tilde L, \tilde L_i,\tilde L_j)\otimes \cC(\tilde L, \tilde L_j,\tilde L_k)\,\iso\, \cC(\tilde L, \tilde L_i,\tilde L_k)\otimes \cC(\tilde L_i, \tilde L_j,\tilde L_k)
\end{equation}
(when $\tilde L$ runs through $\cU_{ijk}$ the above isomorphisms of $\ZZ/4\ZZ$-torsors over $\Spec k$ organize into an isomorphism of $\ZZ/4\ZZ$-torsors over $\cU_{ijk}$). So, $\delta_{kji}$ can be seen as a trivialization of the $\ZZ/4\ZZ$-torsor over $\cU_{ijk}$, whose fibre at $\tilde L$ is
$$
\frac{\cC(\tilde L, \tilde L_i,\tilde L_j)\otimes\cB_{ij}\otimes \cC(\tilde L, \tilde L_j,\tilde L_k)\otimes\cB_{jk}}
{\cC(\tilde L, \tilde L_i,\tilde L_k)\otimes\cB_{ik}}
$$
\end{itemize}

  It is understood that (\ref{eq_for_cC_cocycle}) is the isomorphism given by Proposition~\ref{Con_Maslov_is_cocycle}.
Finally, the 2-morphisms $\delta$ satisfy the following compatibility property: over any $\cU_{ijks}$ the diagram of 2-morphisms commutes
$$
\begin{array}{ccc}
\sigma_{sk}\comp\sigma_{kj}\comp\sigma_{ji} & \toup{\delta_{kji}} & \sigma_{sk}\comp\sigma_{ki}\\
\downarrow\lefteqn{\scriptstyle \delta_{skj}} && \downarrow\lefteqn{\scriptstyle \delta_{ski}}\\
\sigma_{sj}\comp \sigma_{ji} & \toup{\delta_{sji}} & \sigma_{si}
\end{array}
$$
  
  Though this is not reflected in our notation, the gerb $\hat\cL(\tilde V)$ depends on the family $\{\tilde L_i\}, i\in I$.
  
\begin{Rem} 
\label{Rem_elem_transforms}
Call an \select{elementary transformation} of  an enriched 1-sweeping family the procedure of adding or throwing away one lagrangian $\tilde L\in \cL(\tilde V)$ such that the obtained family is still 1-sweeping (together with a compatible change of the family $\cB_{ij}$). Clearly, the gerb associated to the enriched 1-sweeping family obtained from the original one by an elementary transformation is isomorphic to the initial gerb $\hat\cL(\tilde V)$ over $\cL(\tilde V)$. One can pass from one 1-sweeping family to another by a finite number of elementary transformations. So, the isomorphism class of the gerb $\hat\cL(\tilde V)$ does not depend on a choice of $\{\tilde L_i\}, i\in I$ and of $\cB_{ij}$. 

 This also shows that the corresponding element of $\H^2(\cL(\tilde V), \ZZ/4\ZZ)$ is invariant under the action of (the Greenberg realization of) $\Sp(\tilde V)$.
\end{Rem} 
 
\begin{Rem} 
\label{Rem_ELag_so_on}
i) Remind the variety $ELag(V)$ from Section~2.4.1. A similar $\ZZ/4\ZZ$-gerb can be defined over $ELag(V)$. It is not used in this paper.\\
ii) Consider the $W_2^*$-torsor over $\cL(\tilde V)$ whose fibre over $\tilde L$ is the set of generators of $\det\tilde L$. Since $W_2^*\,\iso\, \Gm\times \A^1$ canonically, it can be seen as a pair: a $\Gm$-torsor and a $\A^1$-torsor over $\cL(\tilde V)$. Denote the corresponding $\A^1$-torsor by $\cL(\tilde V)_{ex}\to \cL(\tilde V)$. Let $\hat\cL(\tilde V)_{ex}$ be the $\ZZ/2\ZZ$-gerb over $\cL(\tilde V)$ obtained from $\cL(\tilde V)_{ex}$ via the exact sequence 
$$
0\to \ZZ/2\ZZ\to \A^1\,\toup{x\mapsto x^2+x}\, \A^1\to 0
$$ 
One may show that $\hat\cL(\tilde V)_{ex}$ identifies with the extension of $\hat\cL(\tilde V)$ under the surjective homomorphism of structure groups $\ZZ/4\ZZ\to \ZZ/2\ZZ$. In particular, $\hat\cL(\tilde V)_{ex}$ admits a $\Sp(\tilde V)$-equivariant structure. This will not be used in the present paper.
\end{Rem} 
 
\medskip\noindent
6.2 Set $Y=\cL(\tilde V)\times \cL(\tilde V)$, the product being taken over $k$. Let $\hat Y$ be the $\ZZ/4\ZZ$-gerb over $Y$ obtained from $\hat\cL(\tilde V)\times\hat\cL(\tilde V)$ via the extension of the structure group 
$$
\ZZ/4\ZZ\times \ZZ/4\ZZ\to \ZZ/4\ZZ, \; (a,b)\mapsto b-a
$$ 
Let $Y_0\subset Y$ be the open subscheme classifying $(\tilde L_1,\tilde L_2)\in Y$ such that $\tilde L_1\cap\tilde L_2=0$. 

\begin{Lm} There exists a canonical section $\gs$ of the gerb $\hat Y\to Y$ over $Y_0$.
\end{Lm}
\begin{Prf}
Assume that the collection $(\{\tilde L_i\}_{i\in I}$, $\cB_{ij})$ giving rise to $\hat\cL(\tilde V)$ is such that $\{\tilde L_i\}_{i\in I}$ is 2-sweeping. 

 Define $\gs$ by gluing. Over the open subscheme $Y_0(\tilde L_i)=\{(\tilde L,\tilde M)\in Y_0\mid \tilde L, \tilde M\in \cU(\tilde L_i)\}$ we consider the $\ZZ/4\ZZ$-torsor, whose fibre at $(\tilde L, \tilde M)$ is $\cC(\tilde M, \tilde L, \tilde L_i)$. For any $i,j$ over 
$Y_0(\tilde L_i)\cap Y_0(\tilde L_j)$ we glue the corresponding sections of our gerb via the isomorphism (given by Proposition~\ref{Con_Maslov_is_cocycle})
$$
\cC(\tilde M, \tilde L, \tilde L_i)\otimes \frac{\cC(\tilde M, \tilde L_i, \tilde L_j)}{\cC(\tilde L, \tilde L_i, \tilde L_j)}\,\iso\, 
\cC(\tilde M, \tilde L, \tilde L_j)
$$ 
The fraction in the above formula is, according to 6.1, exactly the gluing data for the gerb $\hat Y_d$. It is understood that we have chosen the same $\cB_{ij}$ in the numerator and the denominator, so they have disappeared.
\end{Prf} 
 
\begin{Def} 
\label{Def_sheaf_S_V_psi}
Define a perverse sheaf $S_{\tilde V, \psi}$ on $\hat Y$ as follows. The section $\gs$ yields an isomorphism $\hat Y\mid_{Y_0}\,\iso\, Y_0\times B(\ZZ/4\ZZ)$. Let $W_{\psi}$ be the rank one local system on $B(\ZZ/4\ZZ)$ corresponding to the representation $\psi: \ZZ/4\ZZ\to \Qlb^*$. If $p: \Spec k\to B(\ZZ/4\ZZ)$ is the quotient map then $W_{\psi}$ is the direct summand in $p_!\Qlb$ on which $\ZZ/4\ZZ$ acts via $\psi$.
Let $S_{\tilde V, \psi}$ be the intermediate extension of $\Qlb\boxtimes W_{\psi}$ from $\hat Y\mid_{Y_0}$ to $\hat Y$. We refer to $S_{\tilde V, \psi}$ as \select{the finite-dimensional theta-sheaf}.
\end{Def}

\medskip\noindent 
6.3.1 Remind the variety $U_{23}$ defined in Section~5.4, it classifies triples $(\tilde N, \tilde L, \tilde M)$ in $\cL(\tilde V)$ such that $\tilde N\cap \tilde L=\tilde N\cap \tilde M=0$. Our purpose now is to define a morphism of stacks $\hat\nu_{23}: U_{23}\to \hat Y$ extending the projection $\nu_{23}: (\tilde N, \tilde L, \tilde M)\to (\tilde L, \tilde M)$. 
 
 Assume that the collection $(\{\tilde L_i\}_{i\in I}, \cB_{ij})$ is such that $\{\tilde L_i\}_{i\in I}$ is 3-sweeping. Define the section $\hat\nu_{23}$ of the gerb $\hat Y$ over $U_{23}$ by gluing sections on the open subschemes
$$
U_{23}(\tilde L_i):=\{(\tilde N,\tilde L,\tilde M)\in \tilde U_{23}\mid \tilde N, \tilde L,\tilde M\in \cU(\tilde L_i)\}
$$ 
Namely, over the open subscheme $U_{23}(\tilde L_i)$ we consider the $\ZZ/4\ZZ$-torsor whose fibre at $(\tilde N,\tilde L,\tilde M)$ is
$$
\frac{\cC(\tilde M, \tilde N, \tilde L_i)}{\cC(\tilde L, \tilde N, \tilde L_i)}
$$
Over the intersection $U_{23}(\tilde L_i)\cap U_{23}(\tilde L_j)$ we glue the above sections via the isomorphism
\begin{equation}
\label{iso_three_fractions}
\frac{\cC(\tilde M, \tilde N, \tilde L_i)}{\cC(\tilde L, \tilde N, \tilde L_i)}\otimes \frac{\cC(\tilde M, \tilde L_i, \tilde L_j)}{\cC(\tilde L, \tilde L_i, \tilde L_j)}\,\iso\, 
\frac{\cC(\tilde M, \tilde N, \tilde L_j)}{\cC(\tilde L, \tilde N, \tilde L_j)}
\end{equation}
obtained from the two isomorphisms (of Proposition~\ref{Con_Maslov_is_cocycle})
$$
\cC(\cdot, \tilde N, \tilde L_i)\otimes \cC(\cdot, \tilde L_i, \tilde L_j)\,\iso\, \cC(\tilde N, \tilde L_i, \tilde L_j)\otimes \cC(\cdot, \tilde N, \tilde L_j),
$$
where $\cdot$ is $\tilde M$ or $\tilde L$. This completes the definition of $\hat\nu_{23}: U_{23}\to\hat Y$.

\begin{Pp} 
\label{Pp_two_coverings_isomorphic}
The $\ZZ/4\ZZ$-coverings of $U_3$ defined by the following two cartesian squares 
$$
\begin{array}{ccc}
Y_0 & \gets & Y_0\times_{\hat Y} U_3\\
\downarrow\lefteqn{\scriptstyle \gs}&& \downarrow\\
\hat Y & \getsup{\hat\nu_{23}} & U_3
\end{array}
\;\;\;\;\;\;
\begin{array}{ccc}
\Cov(U_3) & \to & \Cov(Q_0(\tilde L))/\GL(\tilde L)\\
\downarrow && \downarrow\\
U_3 & \to & Q_0(\tilde L)/\GL(\tilde L)
\end{array}
$$
are canonically isomorphic.
\end{Pp}
\begin{Prf}
Over $Y_0$ we trivialize the $\ZZ/4\ZZ$-gerb $\hat Y\to Y$ via the section $\gs$. Then the map $\hat\nu_{23}: U_3\to Y_0\times B(\ZZ/4\ZZ)$ is given by some $\ZZ/4\ZZ$-torsor, say $\cF$, over $U_3$. Note that $Y_0\times_{\hat Y} U_3\to U_3$ is the total space of $\cF$. 

 We identify $\cF$ with the $\ZZ/4\ZZ$-torsor $\Cov(U_3)\to U_3$ by gluing. First, $\cF$ itself is constructed by gluing. Over the open subscheme 
$$
U_3(\tilde L_i):=\{(\tilde N,\tilde M,\tilde L)\in U_3\mid \tilde N,\tilde M,\tilde L\in \cU(\tilde L_i)\},
$$
$\cF$ is given by the $\ZZ/4\ZZ$-torsor whose fibre at $(\tilde N, \tilde L, \tilde M)$ is
$$
\frac{\cC(\tilde M, \tilde N, \tilde L_i)}{\cC(\tilde L, \tilde N, \tilde L_i)\otimes \cC(\tilde M, \tilde L, \tilde L_i)}
$$
The fibre of $\Cov(U_3)\to U_3$ over $(\tilde N, \tilde L, \tilde M)$ is $\cC(\tilde N, \tilde L, \tilde M)$. The desired isomorphism of $\ZZ/4\ZZ$-torsors over $U_3(\tilde L_i)$ is the isomorphism
$$
\frac{\cC(\tilde M, \tilde N, \tilde L_i)}{\cC(\tilde L, \tilde N, \tilde L_i)\otimes \cC(\tilde M, \tilde L, \tilde L_i)}\;\iso\;
\cC(\tilde N, \tilde L, \tilde M)
$$
given by Proposition~\ref{Con_Maslov_is_cocycle}. These isomorphisms are compatible with the gluing data, so they yield the desired isomorphism over $U_3=\cup_{i\in I} \, U_3(\tilde L_i)$.
\end{Prf}

\medskip

 For any $\tilde L\in \cL(\tilde V)$ there is an isomorphism
$Y_0\,\iso\, \Sp(\tilde V)/\GL(\tilde L)$, so $Y_0$ is affine.
Note that $U_{23}$ can be seen as a variety of triples $(\tilde N, \tilde L, \tilde\phi)$, where $(\tilde N, \tilde L)\in Y_0$ and $\tilde\phi\in \Sym^{!2}(\tilde L^*)$. Therefore, 
the map (\ref{map_pi_23}) is smooth, and $U_{23}$ is affine. 
Proposition~\ref{Pp_two_coverings_isomorphic}  immediately yields the following.
 
\begin{Cor} 
\label{Cor_theta_descent}
There is a canonical isomorphism
$$
\pi_U^*\SSS_{\psi}[\dimrel(\pi_U)]\,\iso\,
\hat\nu_{23}^*S_{\tilde V,\psi}[\dimrel(\nu_{23})]
$$
of perverse sheaves on $U_{23}$. \QED 
\end{Cor}

\noindent
6.3.2 {\scshape An alternative description of $\hat Y$}

\medskip\noindent
For the map $\nu_{23}: U_{23}\to Y$ consider the scheme $U_Y:=U_{23}\times_Y U_{23}$, it classifies collections 
$$
(\tilde N, \tilde L, \tilde M)\in U_{23}, (\tilde N',\tilde L, \tilde M)\in U_{23}
$$ 
For such a point of  $U_Y$ consider the ordered collection of lagrangians $(\tilde N', \tilde M, \tilde N, \tilde L)$. Let $K_{\tilde N', \tilde M, \tilde N, \tilde L}$ be the corresponding free $R$-module with the symmetric bilinear form given by Definition~\ref{Def_K_many}. It is crucial that in the collection $\tilde N', \tilde M, \tilde N, \tilde L$ indexed by $\ZZ/4\ZZ$ each lagrangian is transversal to the next one, so the bilinear form on $K_{\tilde N', \tilde M, \tilde N, \tilde L}$ is non degenerate. Let $\tilde K$ be a free $R$-module of rank $2n$. This gives a morphism
$$
\pi_{UY}: U_Y\to Q_0(\tilde K)/\GL(\tilde K)
$$
sending $(\tilde N', \tilde M, \tilde N, \tilde L)$ to $K_{\tilde N', \tilde M, \tilde N, \tilde L}$ with the corresponding bilinear form.

 Let $\Cov(U_Y)$ be the $\ZZ/4\ZZ$-torsor over $U_Y$ obtained from $\Cov(Q_0(\tilde K))/\GL(\tilde K)$ by the base change $\pi_{UY}$. For a scheme $S$ and an $S$-point $(\tilde N', \tilde M, \tilde N, \tilde L)$ of $U_Y$ write $\cC(\tilde N', \tilde M, \tilde N, \tilde L)$ for the $\ZZ/4\ZZ$-torsor over $S$ obtained from $\Cov(U_Y)$ by the base change $S\to U_Y$.
 
\begin{Lm} 
\label{Lm_new_isometries}
Let $S$ be a scheme and 
$$
(\tilde N, \tilde L, \tilde M)\in U_{23}, (\tilde N',\tilde L, \tilde M)\in U_{23}, (\tilde N'', \tilde L, \tilde M)\in U_{23}
$$
be an $S$-point of $U_{23}\times_Y U_{23}\times_Y U_{23}$. One has a canonical isomorphism of $\ZZ/4\ZZ$-torsors on $S$
\begin{equation}
\label{iso_three_U_Y}
\cC(\tilde N'', \tilde M, \tilde N', \tilde L)\otimes \cC(\tilde N', \tilde M, \tilde N, \tilde L)\,\iso\, \cC(\tilde N'', \tilde M, \tilde N, \tilde L)
\end{equation}
In particular, for $(\tilde N, \tilde L, \tilde M)\in U_{23}$ it yields a trivialization of the $\ZZ/4\ZZ$-torsor $\cC(\tilde N, \tilde M, \tilde N, \tilde L)$.
\end{Lm} 
\begin{Prf} Combine Proposition~\ref{Pp_strange_isometry} and Remark~\ref{Rem_free_isotropic_submodule_and_Gauss}.
\end{Prf}

\medskip

Let $\hat Y_d$ be the $\ZZ/4\ZZ$-gerb over $Y$ obtained as the descent of the trivial gerb $U_{23}\times B(\ZZ/4\ZZ)$ with respect to the morphism $\nu_{23}:U_{23}\to Y$ for the descent data given by 
\begin{itemize}
\item the isomorphism
$$
\sigma_{UY}: U_Y\times B(\ZZ/4\ZZ)\to U_Y\times B(\ZZ/4\ZZ)
$$
over $U_Y$ sending a given $\ZZ/4\ZZ$-torsor $\cF$ to $\cF\otimes \Cov(U_Y)$;
\item the isomorphism (\ref{iso_three_U_Y}) of $\ZZ/4\ZZ$-torsors over $U_{23}\times_Y U_{23}\times_Y U_{23}$. 
\end{itemize}

 The gerb $\hat Y_d$ is naturally $\Sp(\tilde V)$-equivariant, thus it can be seen as a $\ZZ/4\ZZ$-gerb over the stack quotient $Y/\Sp(\tilde V)$. It is understood that $\Sp(\tilde V)$ acts on $Y$ diagonally. 

\begin{Lm} 1) There exists a canonical isomorphism $\hat Y_d\,\iso\, \hat Y$ of $\ZZ/4\ZZ$-gerbs over $Y$. \\
2) The restriction of $\hat Y_d$ to the diagonal $\cL(\tilde V)\hook{} Y$ admits a canonical $\Sp(\tilde V)$-equivariant section $\gt: \cL(\tilde V)\to \hat Y_d$.
\end{Lm}
\begin{Prf}
Remind the original definition of $\hat Y$. One first picks a 2-sweeping family $(\{\tilde L_i\}_{i\in I})$. For $i\in I$ let $Y(\tilde L_i)=\{(\tilde L, \tilde M\in Y\mid \tilde L, \tilde M\in \cU(\tilde L_i)\}$. Then $\hat Y$ is obtained by gluing the trivial $\ZZ/4\ZZ$-gerbs $Y(\tilde L_i)\times B(\ZZ/4\ZZ)$ over the open subschemes $Y(\tilde L_i)\cap Y(\tilde L_j)$ via the isomorphism sending a $\ZZ/4\ZZ$-torsor $\cF_i$ over $Y(\tilde L_i)\cap Y(\tilde L_j)$ to the $\ZZ/4\ZZ$-torsor whose fibre at $(\tilde L, \tilde M)\in Y(\tilde L_i)\cap Y(\tilde L_j)$ is
$$
(\cF_i)_{\tilde L, \tilde M}\otimes \frac{\cC(\tilde M, \tilde L_i, \tilde L_j)}{\cC(\tilde L, \tilde L_i, \tilde L_j)}
$$
From Proposition~\ref{Pp_chain_condition} and Lemma~\ref{Lm_dihedral} one gets an isomorphism 
$$
\frac{\cC(\tilde M, \tilde L_i, \tilde L_j)}{\cC(\tilde L, \tilde L_i, \tilde L_j)}\,\iso\, \cC(L_j, M, L_i, L)
$$
2) For $(\tilde N, \tilde L, \tilde M)\in U_{23}$ the trivialization of $\cC(\tilde N, \tilde M, \tilde N, \tilde L)$ given by Lemma~\ref{Lm_new_isometries} shows the following. After the diagonal base change $\cL(\tilde V)\to Y$, the tautological  section of the gerb $(U_{23}\times_Y \cL(\tilde V))\times B(\ZZ/4\ZZ)$ is compatible with the descent data with respect to $\nu_{23}\times\id: U_{23}\times_Y \cL(\tilde V)\to \cL(\tilde V)$. It yields the desired section.
\end{Prf}

\begin{Rem} 
\label{Rem_theta_sheaf_is_hat_G-equiv}
The section $\gs: Y_0\to \hat Y_0$ can now be described as follows. Consider the $\ZZ/4\ZZ$-torsor over $U_3$ whose fibre at $(\tilde N, \tilde L, \tilde M)\in U_3$ is $\cC(\tilde M, \tilde L, \tilde N)$. It is compatible with the descent data for $\hat Y_d$ with respect to the morphism $\nu_{23}: U_3\to Y_0$. Thus $\gs: Y_0\to \hat Y_0$ is $\Sp(\tilde V)$-equivariant. This implies that $S_{\tilde V, \psi}$ is also $\Sp(\tilde V)$-equivariant.
\end{Rem}

\medskip 
\noindent
6.4 {\scshape The metaplectic group} Let $G$ be the Greenberg realization of $\Sp(\tilde V)$. It acts naturally on $\cL(\tilde V)$. According to Appendix~A.2, the gerb $\hat\cL(\tilde V)$ yields a group stack $\hat G$ over $G$. From Remark~\ref{Rem_elem_transforms} we conclude that the homomorphism $\hat G\to G$ is surjective. We refer to $\hat G$ as \select{the metaplectic group}. 

 In the rest of Section~6.4 we prove the following 
\begin{Pp} 
\label{Pp_Z/2Z-torsor_is_trivial}
Any $\ZZ/2\ZZ$-torsor over $\cL(\tilde V)$ is trivial.
\end{Pp}   
By A.3, this implies that $\hat G$ fits into an exact sequence $1\to B(\ZZ/4\ZZ)\to \hat G\to G\to 1$ and is algebraic. Besides, $\hat G$ acts naturally on $\hat\cL(\tilde V)$, and the projection $
\hat\cL(\tilde V)\to \cL(\tilde V)$ is equivariant with respect to $\hat G\to G$. 
 
\begin{Lm} 
\label{Lm_functions_on_vector_bunldes}
Let $Z$ be a $k$-scheme, $q: \WW\to Z$ be the total space of a vector bundle $W$ on $Z$. One has a canonical isomorphism of sheaves of $\cO_Z$-algebras $q_*\cO\,\iso\, \oplus_{d\ge 0} \Sym^{*d}(W^*)$. \QED
\end{Lm}

\begin{Lm} 
\label{Lm_Z/2Z_torsors_for_P1}
Let $W$ be a vector bundle on $\PP^1_k$ of the form $W=\oplus_{i=1}^r \cO(n_i)\otimes W_i$, where $n_i>0$ and $W_i$ are fixed finite-dimensional $k$-vector spaces.  Let $\WW$ be the total space of $W$. Then any $\ZZ/2\ZZ$-torsor on $\WW$ is trivial.
\end{Lm}
\begin{Prf} Let $q: \WW\to \PP^1$ be the projection. By Lemma~\ref{Lm_functions_on_vector_bunldes}, $\H^0(\WW, \cO)=k$. One has the exact sequence of groups $1\to \ZZ/2\ZZ\to \A^1\toup{a}\A^1\to 0$ over $\Spec k$, where $a(x)=x^2+x$. So, it suffices to show that $\H^1(\WW, \cO)\toup{a_*} \H^1(\WW,\cO)$ is injective. The space $\H^1(\WW,\cO)=\oplus_{d\ge 0} \H^1(\PP^1,  \Sym^{*d} (W^*))$ is graded by $d\ge 0$. Let $x\in \H^1(\WW, \cO)$ non zero, let $x_d$ be its non zero component of the biggest degree $d\ge 0$. It suffices to show that $x_d^2\in \H^1(\PP^1, \Sym^{*2d}(W^*))$ does not vanish. One has 
$$
\Sym^{*d}(W^*)\,\iso\, \sum_{d_1+\ldots+d_r=d} \cO(-\sum_i d_in_i)\otimes \Sym^{*d_1}(W_1^*)\otimes\ldots\otimes\Sym^{*d_r}(W_r^*)
$$ 
So, our assertion follows from the fact that for $m\ge 0$ the map $\H^1(\PP^1, \cO(-m))\toup{v_*} \H^1(\PP^1, \cO(-2m))$ is injective, where $v: \cO(-m)\to \cO(-2m)$, $x\mapsto x^2$ is a homomorphism of sheaves of abelian groups on $\PP^1_k$. We also used the property that for a finite-dimensional $k$-vector space $U$ and $u\in \Sym^{*d}U$ the condition $u^2=0$ in $\Sym^{*2d}U$ implies $u=0$.
\end{Prf}
 
\medskip
\begin{Prf}\select{of Proposition~\ref{Pp_Z/2Z-torsor_is_trivial}}

\medskip\noindent
Let $\cE$ be the vector bundle over $\cL(V)$ whose fibre at $L\in \cL(V)$ is $\Sym^{!2}(L^*)$. The projection $\cL(\tilde V)\to\cL(V)$ is a torsor under $F^*\cE$, the inverse image of $\cE$ by the Frobenius map $F$. For $n=1$ one has $\cL(V)\,\iso\, \PP^1_k$, and $F^*\cE$ is isomorphic to the line bundle $\cO(4)$  on $\PP^1_k$.  

 Let $\cF$ be a $\ZZ/2\ZZ$-torsor over $\cL(\tilde V)$. It suffices to show that $\cF$ is constant along the fibres of $\cL(\tilde V)\to\cL(V)$.

\Step 1 Pick an isometry $V\,\iso\, V_1\oplus V_2$, where $V_i$ are symplectic vector spaces, $\dim V_1=2$, $\dim V_2=2n-2$. Pick a $k$-point $L_2\in \cL(V_2)$. It yields a closed immersion $\cL(V_1)\to \cL(V)$, $L_1\mapsto L_1\oplus L_2$.  Let $Z_1=\cL(V_1)\times_{\cL(V)} \cL(\tilde V)$. 
 
 Let $\cE_1$ be the vector bundle on $\cL(V_1)$ whose fibre at $L_1$ is $\Sym^{!2}(L_1^*)\oplus (L_1\otimes L_2)^*$. Then $Z_1\to \cL(V_1)$ is a torsor under the vector bundle $F^*(\cE_1\oplus \Sym^{!2}(L_2^*))$. Any such torsor is trivial. Let $Z$ be the total space of $F^*\cE_1$ then $Z_1\,\iso\, Z\times U$, where $U=F^*\Sym^{!2}(L_2^*)$ is an affine space over $k$. By Lemma~\ref{Lm_Z/2Z_torsors_for_P1}, the $*$-restriction $\cF\mid_{Z_1}$ descends under the projection $Z\times U\to U$ to a local system on $U$.

\medskip

\Step 2 Pick a $k$-point $L\in \cL(V)$. Let $T$ be the fibre of $\cL(\tilde V)\to\cL(V)$ over $L$. For any decomposition $L\,\iso\, L_1\oplus L_2$ into a direct sum of vector subspaces with $\dim L_1=1$, the torsor $\cF\mid_T$ is equivariant under the action of the vector space $F^*(\Sym^{!2}(L_1^*)\oplus (L_1\otimes L_2)^*)$. Indeed, this follows from Step 1. Since the decomposition $L\,\iso\, L_1\oplus L_2$ was arbitrary, the $\ZZ/2\ZZ$-torsor $\cF\mid_T$ is trivial.
\end{Prf} 
  
\medskip 
\noindent
6.5 {\scshape CASE $n=1$.} 

\medskip\noindent
6.5.1 Let us give some explicit formulas for $\hat\cL(\tilde V)$ in the simplest case $n=1$. Take $\tilde V$ to be the free $R$-module with a symplectic base $e_1,e_2$ such that $\tilde\omega(e_1, e_2)=1$. Then $\cL(\tilde V)$ is the Greenberg scheme of $\PP^1_R$. 

 Let $\tilde L_i\subset\tilde V$ be the $R$-submodule generated by $e_i$. Set $\cU_i=\cU(\tilde L_i)$ then $\cU_1\cup\cU_2=\cL(\tilde V)$. We have the isomorphism $R\,\iso\,\cU_1$ sending $a\in R$ to the $R$-submodule in $\tilde V$ generated by $ae_1+e_2$. We have the isomorphism $R\,\iso\,\cU_2$ sending $b$ to the $R$-submodule generated by $e_1+be_2$. So, $\cU_1\cap\,\cU_2\,\iso\, R^*$ and the corresponding identification is given by $b=a^{-1}$.

 View $\tilde L\in \cU_{12}$ as a $R$-submodule in $\tilde V$ generated by $e_1+be_2$ with $b\in R^*$. In this notation the $\ZZ/4\ZZ$-torsor $\Cov(\cU_{12})\to \cU_{12}$ becomes the $\ZZ/4\ZZ$-torsor over $R^*$ whose fibre over $b=(b_0,b_1)$ is $\{z\in R\mid Fz-z=(b_1b_0^{-2},0)\}$, here $b$ is written in Witt coordinates. Indeed, this follows from Lemma~\ref{Lm_covering_case_n=1} and the fact that $C(\tilde L,\tilde L_1,\tilde L_2)$ is the Gauss sum for the quadratic form $x\mapsto -b\tilde x^2$ (here $\tilde x\in R$ is any lifting of $x\in k$). 

 Consider the $R$-torsor $\cT\cL(\tilde V)\to\cL(\tilde V)$ defined as the gluing of the trivial $R$-torsors over $\cU_i$ by the 1-cocycle $h: \cU_{12}\to R$ sending $b=(b_0,b_1)$ to $(b_1b_0^{-2},0)$. The $R$-torsor $\cT\cL(\tilde V)$ over $\cL(\tilde V)$ is nontrivial (even its extension of scalars via $R\to k$ is a nontrivial $\A^1$-torsor over $\cL(\tilde V)$). The gerb $\hat\cL(\tilde V)$ is obtained from this $R$-torsor via the exact sequence $0\to \ZZ/4\ZZ\to R\,\toup{z\mapsto Fz-z}\, R\to 0$.

 Set $\Sp^1(V)=\{g\in\Sp(\tilde V)\mid g=\id\!\mod 2\}$. It is easy to check that $\cT\cL(\tilde V)$ admits a $\Sp^1(V)$-equivariant structure.
 
\medskip\noindent
6.5.2 Following \cite{GH}, consider the variety $\cL^0(\tilde V)$ classifying $\tilde L\in\cL(\tilde V)$ together with a generator $o_{\tilde L}$ of the $R$-module $\det L$. A point $\tilde L^0=(\tilde L, o_{\tilde L})$ of $\cL^0(\tilde V)$ is called \select{an oriented lagrangian}.

 Let $\cT\cL^0(\tilde \cV)\to \cL^0(\tilde V)$ be the $R$-torsor obtained from $\cT\cL(\tilde V)\to\cL(\tilde V)$ by the base change $\cL^0(\tilde V)\to \cL(\tilde V)$. 

  Let $\cV_i\subset \cL^0(\tilde V)$ be the open subscheme classifying oriented lagrangians $\tilde L^0$ generated by $o_{\tilde L}=(a,b)\in R\times R$ such that $b\in R^*$ (resp., $a\in R^*$) for $i=1$ (resp., for $i=2$). 
  
  Let $\cT^0(\tilde V)\to \cL^0(\tilde V)$ be the $2R$-torsor defined as the gluing of the trivial $2R$-torsors over $\cV_i$ via
the cocycle over $\cV_{12}=\cV_1\cap\cV_2$ sending $(a,b)$ to 
$$
(0, \frac{a_1b_1}{a_0^2b_0^2})\in 2R
$$
The $R$-torsor $\cT\cL^0(\tilde \cV)\to \cL^0(\tilde V)$ is isomorphic to the extension of scalars of $\cT^0(\tilde V)$ under the inclusion $2R\to R$.

  One checks that the $\Sp(\tilde V)$-orbit in $\H^1(\cL^0(\tilde V), \cO)$ passing through the $2R$-torsor $\cT^0(\tilde V)$ is not a point. For example, take $v\in R$ and $g\in\Sp(\tilde V)$ given by $g(e_1)=e_1$ and $g(e_2)=e_2+ve_1$. The $2R$-torsors $g^*\cT^0(\tilde V)$ and $\cT^0(\tilde V)$ are not isomorphic over $\cL^0(\tilde V)$ unless $v_0=0$, here $v=(v_0,v_1)$ is written in Witt coordinates.
 
 Actually, there exists a $2R$-torsor, say $\cT(\tilde V)\to \cL(\tilde V)$, whose restriction to $\cL^0(\tilde V)$ is isomorphic to $\cT^0(\tilde V)\to \cL^0(\tilde V)$.

\bigskip

\centerline{\scshape 7. Canonical interwining operators}

\medskip\noindent
7.1.1 Our purpose now is to generalize the theory of canonical interwining operators (\cite{L2}, Theorem~1) to the case of characteristic two. 

 Remind that $\tilde V$ is a free $R$-module of rank $2n$ with symplectic form $\tilde\omega: \tilde V\otimes\tilde V\to R$, and $V=\tilde V\otimes_R k$. Pick a bilinear form $\tilde\beta: \tilde V\times \tilde V\to R$ satisfying (\ref{def_tilde_beta}). 
 
  Let $\beta: V\times V\to R$ be defined as in 2.2. It gives rise to the Heisenberg group $H=H(V)=V\times R$ with operation (\ref{def_Heis_group}). We view it as an an algebraic group over $k$ (the Greenberg realization of the corresponding $R$-scheme).
The center of $H$ is $Z(H)=\{(0,z)\in H(V)\mid z\in R\}$.   

  The affine symplectic group $\ASp(V)$ is defined as in 2.2 (it is understood that for $(g,\alpha)\in \ASp(V)$ the map $\alpha: V\to R$ must be a morphism of $k$-schemes), it is an algebraic group over $k$ acting on $H$ by automorphisms. 
  
  Let $G$ be the Greenberg realization of $\Sp(\tilde V)$ over $k$.  As in 2.2, one defines a homomorphism $\xi: G\to \ASp(V)$ of algebraic groups over $k$. Let $ELag(V)$ be defined as in Section~2.4.1. The map $\epsilon: \cL(\tilde V)\to ELag(V)$ defined as in 2.4.2 is a morphism of schemes  
over $\cL(V)$, the action of $\ASp(V)$ on $ELag(V)$ is algebraic, and the diagram (\ref{diag_action_of_ASp(V)_and_G}) commutes.   
  
 Given a $k$-point $(L,\tau)\in ELag(V)$, write $\cH_L$ for the category of $\Qlb$-perverse sheaves on $H$, which are equvariant with respect to $L$ acting on $H$ by the left multiplication via $\tau$, and also $(Z(H), \cL_{\psi})$-equivariant. This is a full subcategory in $\P(H)$. Write $\D\cH_L\subset \D(H)$ for the full subcategory of objects whose all perverse cohomologies lie in $\cH_L$. 

 For a $k$-point $\tilde L\in \cL(\tilde V)$ set $\cH_{\tilde L}=\cH_L$ for the enhanced lagrangian $L=\epsilon(\tilde L)$. For a $k$-point $g\in G$ the inverse image under the map $H\to H$, $h\mapsto g^{-1}h$ yields an equivalence $g: \cH_{\tilde L}\,\iso\, \cH_{g\tilde L}$, which we denote by $g$ by some abuse of notation.

\medskip\noindent
7.1.2 Write $\hat L$ for a point of $\hat\cL(\tilde V)$ over $\tilde L\in \cL(\tilde V)$. For $g\in\hat G$, $\hat L\in\hat\cL(\tilde V)$ we write $g\hat L$ for the image of $(g, \hat L)$ under the action map (defined in A.2)
$$
\hat G\times \hat\cL(\tilde V)\to \hat\cL(\tilde V)
$$ 
Write $(\hat L: \hat M)$ for the image of a pair $(\hat L, \hat M)\in \hat\cL(\tilde V)\times\hat\cL(\tilde V)$ in $\hat Y$. 

 For a pair of points $\hat L, \hat M\in \hat\cL(\tilde V)$ we will define a canonical interwining operator 
\begin{equation}
\label{functor_cF_CIO}
\cF_{\hat L, \hat M}: \cH_{\tilde M}\to \cH_{\tilde L}
\end{equation}
depending only on the image $(\hat L: \hat M)\in\hat Y$. 
They will be equipped with isomorphisms
\begin{itemize}
\item $\cF_{\hat L, \hat L}\;\iso\;\id$
\item $\cF_{\hat N, \hat L}\comp \cF_{\hat L, \hat M}\;\iso\; \cF_{\hat N, \hat M}$
\item for $g\in \hat G$ one has $g\comp \cF_{\hat L, \hat M}\comp g^{-1}\;\iso\; \cF_{g\hat L, g\hat M}$
\end{itemize}
satisfying natural compatibility properties. As in \cite{L2}, the functor $\cF_{\hat L, \hat M}$ will be defined as a convolution with a suitable complex on $H$, and as $\hat L, \hat M$ vary, these complexes will organize into a perverse sheaf $F$ on $\hat\cL(\tilde V)\times\hat\cL(\tilde V)\times H$. 

  Denote by $C\to \cL(\tilde V)$ the vector bundle whose fibre at $\tilde L\in \cL(\tilde V)$ is $L=\tilde L\otimes_R k$. Consider the maps
$$
\pr,\act_{lr}: C\times C\times H\to \cL(\tilde V)\times\cL(\tilde V)\times H,
$$ 
where $\act_{lr}$ sends $(\tilde L, l\in L, \tilde M, m\in M, h)$ to $(\tilde L, \tilde M, \tau_{\tilde L}(l)h\tau_{\tilde M}(m))$, and $\pr$ sends the above point to $(\tilde L,\tilde M, h)$. Say that a perverse sheaf $K$ on $\cL(\tilde V)\times\cL(\tilde V)\times H$ is \select{$\act_{lr}$-equivariant} if it admits an isomorphism
$$
\act_{lr}^*K\,\iso\, \pr^*K
$$
satisfying the usual associativity condition and whose restriction to the unit section is the identity (such isomorphism is unique if it exists). One has a similar definition for $\hat\cL(\tilde V)\times\hat\cL(\tilde V)\times H$. 
Let
$$
\act_{\hat G}: \hat G\times \hat\cL(\tilde V)\times\hat\cL(\tilde V)\times H\to \hat\cL(\tilde V)\times\hat\cL(\tilde V)\times H
$$
be the map sending $(g, \hat L, \hat M, h)$ to $(g\hat L, g\hat M, gh)$. This is an action map in the sense of A.2.1, so one has a notion of a $\hat G$-equivariant perverse sheaf on  $\hat\cL(\tilde V)\times\hat\cL(\tilde V)\times H$ (cf. A.4).

 For a scheme $S$ and $K,K'\in \D(S\times H)$ define their convolution $K\ast K'\in \D(S\times H)$ by
$$
K\ast K'=\mult_!(\pr_{12}^*K\otimes \pr_{13}^*K')[n+2-2\dim\cL(\tilde V)], 
$$
where $\pr_{12},\pr_{13}: S\times H\times H\to S\times H$ are the projections, and $\mult: H\times H\to H$ is the product map sending $(h_1,h_2)$ to $h_1h_2$. The above shift is chosen to that the formula (\ref{formula_without_shifts}) below holds without any shift.

 Let 
$$
i_{\triangle}: (\cL(\tilde V)\times H)_{\triangle}\hook{} \cL(\tilde V)\times H
$$ 
be closed subscheme of $(\tilde L, h)\in \cL(\tilde V)\times H$ such that there exist $x\in L, z\in Z(H)$ with $h=\tau_{\tilde L}(x)z$. Let $\alpha_{\triangle}: (\cL(\tilde V)\times H)_{\triangle}\to Z(H)$ be the map sending the above point to $z\in Z(H)$. 

\begin{Th} 
\label{Th_CIO}
There exists an irreducible perverse sheaf $F$ on $\hat\cL(\tilde V)\times\hat\cL(\tilde V)\times H$ with the following properties:
\begin{itemize}
\item for the diagonal map $i: \hat\cL(\tilde V)\times H\to \hat\cL(\tilde V)\times\hat\cL(\tilde V)\times H$ the complex $i^*F$ identifies canonically with the inverse image of 
$$
(i_{\triangle})_!\alpha_{\triangle}^*\cL_{\psi}[2\dim \cL(\tilde V)+n+2]
$$
under the projection $\hat\cL(\tilde V)\times H\to \cL(\tilde V)\times H$;
\item $F$ is $\act_{lr}$-equivariant and $(Z(H),\cL_{\psi})$-equivariant;
\item $F$ is $\hat G$-equivariant;
\item convolution property for $F$ holds, namely for $ij$-th projections
$$
q_{ij}: \hat\cL(\tilde V)\times \hat\cL(\tilde V)\times\hat\cL(\tilde V)\times H \to \hat\cL(\tilde V)\times\hat\cL(\tilde V)\times H
$$
inside the triple $\hat\cL(\tilde V)\times \hat\cL(\tilde V)\times\hat\cL(\tilde V)$ one has $(q_{12}^*F)\ast (q_{23}^*F)\,\iso\, q_{13}^*F$ canonically.
\end{itemize}
\end{Th} 
 
 The proof of Theorem~\ref{Th_CIO} is given in Sections~7.1.3-7.1.5.
 
\medskip\noindent
7.1.3  Remind that $Y_0\subset Y$ is the open subscheme classifying $(\tilde L,\tilde M)\in Y$ such that $\tilde L\cap \tilde M=0$. Define a perverse sheaf $F_0$ on $Y_0\times H$ as follows. Let 
$$
\alpha_0: Y_0\times H\to Z(H)
$$ 
be the map sending $(\tilde L, \tilde M, h)$ to $z$, where $z\in Z(H)$ is uniquely defined by the property that there exist $l\in L, m\in M$ such that $h=\tau_{\tilde L}(l)\tau_{\tilde M}(m)z$. Set
$$
F_0=\alpha_0^*\cL_{\psi}[\dim (Y_0\times H)]
$$
 
 Let $U_{13}$ be the scheme classifying $(\tilde L, \tilde N, \tilde M)\in \cL(\tilde V)^3$ such that $\tilde N\cap\tilde L=\tilde N\cap \tilde M=0$. Let $\nu_{13}: U_{13}\to Y$ send $(\tilde L, \tilde N, \tilde M)$ to $(\tilde L, \tilde M)$. Define $\nu_{12}, \nu_{23}: U_{13}\to Y_0$ by 
$$
\nu_{12}(\tilde L, \tilde N, \tilde M)=(\tilde L, \tilde N)\;\;\;\;\;\;
\nu_{23}(\tilde L, \tilde N, \tilde M)=(\tilde N, \tilde M)
$$  
Let $i_U: Y_0\to U_{13}$ send $(\tilde N, \tilde L)$ to $(\tilde L, \tilde N, \tilde L)$. The open subscheme $U_3\subset U_{13}$ classifies triples of pairwise transverse lagrangians in $\cL(\tilde V)$. Let $p_{Y_0}: Y_0\to \cL(\tilde V)$ send $(\tilde N, \tilde L)$ to $\tilde L$. Remind the map $\pi_U$ defined by (\ref{map_pi_3}). 
 
\begin{Lm}  
\label{Lm_first_true_convolution}
1) The complex 
$$
(\nu_{12}^* F_0)\ast (\nu_{23}^* F_0)[\dim\cL(\tilde V)]
$$
is an irreducible perverse sheaf on $U_{13}\times H$. 
For the map $p_{Y_0}\times\id: Y_0\times H\to \cL(\tilde V)\times H$ one has canonically
$$
i_U^*((\nu_{12}^* F_0)\ast (\nu_{23}^* F_0))\,\iso\, (p_{Y_0}\times\id)^* (i_{\triangle})_!\alpha_{\triangle}^*\cL_{\psi}[2\dim \cL(\tilde V)+n+2]
$$
over $Y_0\times H$. \\
2) There is a canonical isomorphism 
$$
(\nu_{12}^* F_0)\ast (\nu_{23}^*F_0)[\dim\cL(\tilde V)]\,\iso\, \pi_U^*\SSS_{\psi}\otimes \nu_{13}^*F_0[2n^2]
$$
of perverse sheaves over $U_3\times H$. Here, by abuse of notation, $\nu_{13}: U_3\to Y_0$ is the restriction of $\nu_{13}$.
\end{Lm} 
\begin{Prf}
1) Let us give an argument at the classical level (of $K$-groups)
and explain the changes needed for geometrization. 
For a point $(\tilde L,\tilde M)\in Y_0$ write 
$\FF_{\tilde L, \tilde M}$ for the restriction of $F_0$ to this point, so $\FF_{\tilde L, \tilde M}$ is a `function' on $H$. 
Given $(\tilde L, \tilde N, \tilde M)\in U_{13}$, let us calculate the convolution
\begin{equation}
\label{fun_one}
(\FF_{\tilde L, \tilde N}\ast \FF_{\tilde N, \tilde M})(h)=\int_{h_1\in H} \FF_{\tilde L, \tilde N}(hh_1^{-1})\FF_{\tilde N, \tilde M}(h_1)dh_1
\end{equation}
as a function of $h\in H$, here $dh_1$ is a `Haar measure'. Write $N=\tilde N\otimes_R k$, and similarly for $L,M$. Because of equivariance properties of (\ref{fun_one}), we may assume $h=\tau_{\tilde N}(x)$ for $x\in N$. Write $h_1=\tau_{\tilde N}(y)\tau_{\tilde M}(u)z$ with $z\in Z(H), u\in M, y\in N$. Then (\ref{fun_one}) equals the volume of $N\times Z$ multiplied by
\begin{equation}
\label{fun_two}
\int_{u\in M} \FF_{\tilde L, \tilde N}(\tau_{\tilde N}(x)\tau_{\tilde M}(u))du=
\int_{u\in M} \FF_{\tilde L, \tilde N}(\tau_{\tilde M}(u))\psi(\omega(x,u))du
\end{equation}
We have used the equality $\tau_{\tilde N}(x)\tau_{\tilde M}(u)=
\tau_{\tilde M}(u)\tau_{\tilde N}(x)(0, \omega(x,u))$. The formula (\ref{fun_two}) shows that the resulting complex on $N$ is the Fourier transform of a rank one local system on $M$ (the symplectic form induces an isomorphism $M^*\,\iso\, N$). So, our first assertion follows from the fact that the Fourier transform preserves perversity and irreducibility.
 
 Assume further that $\tilde M=\tilde L$ then (\ref{fun_two}) equals
$$
\int_{u\in L} \psi(\omega(x,u))du=\left\{\begin{array}{ll}
0, & \mbox{for}\; x\ne 0\\
vol(L), & \mbox{for}\; x=0
\end{array}
\right.
$$ 
The geometrization is straightforward, our second assertion follows.

\medskip\noindent
2) Assume that $(\tilde L, \tilde N, \tilde M)\in U_3$. Take $x\in N$ and continue the calculation of (\ref{fun_two}) from 1) as follows. Let $r: N\to L$ be the $k$-linear map such that $M=\{r(w)-w\in L\oplus N\mid w\in N\}$. Consider $Q:=Q_{\tilde L, \tilde N, \tilde M}\in Q^!(N^*)$ defined by (\ref{eq_def_for_Q_123}), that is, $Q: N\to R$ is given by
\begin{equation}
\label{formula_for_Q_enchanced_geom}
\tau_{\tilde M}(r(w)-w)\tau_{\tilde N}(w)\tau_{\tilde L}(-r(w))=(0, Q(w))
\end{equation}
for any $w\in N$.  Now  (\ref{fun_two}) equals
\begin{equation}
\label{fun_three}
\int_{w\in N} \FF_{\tilde L, \tilde N}(\tau_{\tilde M}(r(w)-w))\psi(\omega(x, r(w)-w))dw=
\int_{w\in N} \psi(Q(w)+\omega(x, r(w)))dw
\end{equation}
Remind that for $w,w_1\in N$ one has 
$$
Q(w+w_1)=Q(w)+Q(w_1)+\omega(r(w), w_1)
$$ 
(cf. Section~2.5). After the change of variables $w=t+x$, $t\in N$ the expression (\ref{fun_three}) rewrites as 
$$
\psi(Q(x)+\omega(x,r(x)))\int_{t\in N} \psi(Q(t))dt
$$
Now (\ref{formula_for_Q_enchanced_geom}) with $w=x$ also yields
$$
\FF_{\tilde L, \tilde M}(\tau_{\tilde N}(x))=\psi(Q(x)+\omega(x, r(x)))
$$
Remind the notation $C(\tilde L, \tilde N, \tilde M)=\int_{t\in N} \psi(Q(t))dt$ given by (\ref{Gauss_sum_first_classical}). Combining with 1) we get 
$$
(\FF_{\tilde L, \tilde N}\ast \FF_{\tilde N, \tilde M})(h)=vol(N\times Z) C(\tilde L, \tilde N, \tilde M)\FF_{\tilde L, \tilde M}(h)
$$
Because of equivariance properties, the latter formula holds for all $h\in H$. By Remark~\ref{Rem_comparison_two_q_forms}, 
$Q$ is the image of $\pi_U(\tilde L, \tilde N, \tilde M)$ in $Q^!(N^*)$. The above proof goes through also in the geometric setting.
\end{Prf} 
 
\medskip\noindent
7.1.4  Let $\Theta_{\psi}$ be the rank one local system on $\hat Y_0$ defined by $\Theta_{\psi}=S_{\tilde V,\psi}\mid_{\hat Y_0}[-\dim Y_0]$. The map $\hat\nu_{23}$ has been defined in Section~6.3.1. Consider the diagram
$$
\begin{array}{ccc}
\hat Y_0 \getsup{\hat\nu_{12}} & U_3 & \toup{\hat\nu_{13}} \hat Y_0\\
& \downarrow\lefteqn{\scriptstyle \hat\nu_{23}}\\
& \hat Y_0,
\end{array}
$$
where $\hat\nu_{12}, \hat\nu_{13}$ are defined by
$\hat\nu_{13}(\tilde L, \tilde N, \tilde M)=\hat\nu_{23}(\tilde N, \tilde L, \tilde M)$ and $\hat\nu_{12}(\tilde L,\tilde N, \tilde M)=\hat\nu_{23}(\tilde M, \tilde L, \tilde N)$ for $(\tilde L, \tilde N, \tilde M)\in U_3$. From Lemma~\ref{Lm_dihedral} one derives the following.

\begin{Lm} 
\label{Lm_Theta_permutations}
One has canonical isomorphisms $\hat\nu_{12}^*\Theta_{\psi}\;\iso\; \hat\nu_{23}^*\Theta_{\psi}\;\iso\; \hat\nu_{13}^*\Theta_{\psi}^{-1}$ over $U_3$. \QED
\end{Lm}
 
\begin{Def} Let $\hat F_0$ be the perverse sheaf on $\hat Y_0\times H$ given by $\hat F_0=\pr_1^*\Theta_{\psi}^{-1}\otimes F_0$. It is understood that we take the inverse image of $F_0$ under the projection $\hat Y_0\times H\to Y_0\times H$. Let $F$ be the intermediate extension of $\hat F_0$ under $\hat Y_0\times H\hook{} \hat Y\times H$. The restriction of $F$ under the natural map $\hat\cL(\tilde V)\times \hat\cL(\tilde V)\times H\to \hat Y\times H$ is also denoted by $F$. Note that $\hat F_0$ is $\hat G$-equivariant (cf. Remark~\ref{Rem_theta_sheaf_is_hat_G-equiv}).
\end{Def}

 Combining Lemmas~\ref{Lm_first_true_convolution}, \ref{Lm_Theta_permutations} with Corollary~\ref{Cor_theta_descent}, one gets the following. 

\begin{Cor}  
\label{Cor_convolution_with_theta_multiple}
There is a canonical isomorphism over $U_3\times H$
\begin{equation}
\label{formula_without_shifts}
(\hat\nu_{12}^*\hat F_0)\ast (\hat\nu_{23}^*\hat F_0)\,\iso\, (\hat \nu_{23}^*\Theta_{\psi}^2) \otimes\hat\nu_{13}^*\hat F_0
\end{equation}
\end{Cor} 

 Denote also by $\hat\nu_{13}: U_{13}\to \hat Y$ the map $\hat\nu_{13}(\tilde L, \tilde N, \tilde M)=\hat\nu_{23}(\tilde N, \tilde L, \tilde M)$ for $(\tilde L, \tilde N, \tilde M)\in U_{13}$. The cartesian square
$$
\begin{array}{ccc}
U_3\times H  & \hook{} & U_{13}\times H\\
\downarrow\lefteqn{\scriptstyle \hat\nu_{13}\times\id} && \downarrow\lefteqn{\scriptstyle \hat\nu_{13}\times\id}\\
\hat Y_0\times H & \hook{} & \hat Y\times H
\end{array}
$$
together with Lemma~\ref{Lm_first_true_convolution} yield a canonical isomorphism over $U_{13}\times H$
\begin{equation}
\label{iso_F_explicit_formula}
(\hat\nu_{13}\times\id)^*F\,\iso\, (\nu_{12}^* F_0)\ast (\nu_{23}^* F_0)
\end{equation}
obtained by the intermediate extension from $U_3\times H$. This gives an explicit formula for $F$. Further, the diagram is canonically 2-commutative
$$
\begin{array}{ccc}
Y_0 & \toup{i_U} & U_{13}\\
\downarrow\lefteqn{\scriptstyle p_{Y_0}} && \downarrow\lefteqn{\scriptstyle \hat\nu_{13}}\\
\cL(\tilde V) & \toup{\gt} & \hat Y
\end{array}
$$
Restricting (\ref{iso_F_explicit_formula}) under $i_U\times\id: Y_0\times H\to U_{13}\times H$, one gets an isomorphism
$$
(p_{Y_0}\times\id)^*(\gt\times\id)^*F\,\iso\, (p_{Y_0}\times\id)^*
(i_{\triangle})_!\alpha_{\triangle}^*\cL_{\psi}[2\dim \cL(\tilde V)+n+2]
$$
Since $p_{Y_0}$ has connected fibres, the latter isomorphism
descends under $p_{Y_0}\times\id: Y_0\times H\to \cL(\tilde V)\times H$ to an isomorphism
$$
(\gt\times\id)^*F\,\iso\, (i_{\triangle})_!\alpha_{\triangle}^*\cL_{\psi}[2\dim \cL(\tilde V)+n+2]
$$

 By construction, $\hat F_0$ is $\act_{lr}$-equivariant, $(Z(H), \cL_{\psi})$-equivariant, and $\hat G$-equivariant (this property holds over $\hat Y_0\times H$ and is preserved by the intermediate extension).
 
\medskip\noindent
7.1.5 Consider the scheme $Y\times_{\cL(\tilde V)} Y$ classifying $((\tilde L, \tilde N)\in Y, (\tilde N, \tilde M)\in Y)$. Define the gerb  $\hat Y\times_{\cL(\tilde V)} \hat Y$ over it as the $\ZZ/4\ZZ\times \ZZ/4\ZZ$-gerb obtained from  $\hat\cL(\tilde V)\times \hat\cL(\tilde V)\times \hat\cL(\tilde V)$ by extension of the structure group
$$
\ZZ/4\ZZ\times \ZZ/4\ZZ\times \ZZ/4\ZZ\to  
\ZZ/4\ZZ\times \ZZ/4\ZZ, \; (a,b,c)\mapsto (b-a, c-b)
$$  
So, a point of $\hat Y\times_{\cL(\tilde V)} \hat Y$ is a collection $((\hat L: \hat N)\in \hat Y, (\hat N: \hat M)\in \hat Y)$. Extending further the structure group with respect to 
$\ZZ/4\ZZ\times \ZZ/4\ZZ\to\ZZ/4\ZZ, \; (u,v)\mapsto u+v$, one gets the gerb $\hat Y\times \cL(\tilde V)$, the corresponding morphism of stacks 
$$
\gamma': \hat Y\times_{\cL(\tilde V)} \hat Y\to \hat Y\times \cL(\tilde V)
$$
sends $((\hat L:\hat N), (\hat N: \hat M))$ to $((\hat L: \hat M), \tilde N)$. Write $\gamma: \hat Y\times_{\cL(\tilde V)} \hat Y\to \hat Y$ for $\pr_1\comp \gamma'$. A straightforward calculation yields the following.
 
\begin{Lm} 
\label{Lm_correcting_torsor}
Consider the diagram
$$
\begin{array}{cccc}
\hat Y\times_{\cL(\tilde V)} \hat Y & \;\toup{\gamma} & \hat Y & \toup{p} Y\\
\uparrow\lefteqn{\scriptstyle \hat\nu_{12}\times \hat\nu_{23}} & \;\nearrow\lefteqn{\scriptstyle \hat\nu_{13}}\\
U_3,
\end{array}
$$
where $p$ is the structure map sending $(\hat L: \hat M)$ to $(\tilde L, \tilde M)$. The two maps thus obtained from $U_3$ to $Y$ coincide, but the triangle in the above diagram is not 2-commutative. More precisely, the two sections so obtained of the gerb $\hat Y$ over $U_3$ differ by the $\ZZ/4\ZZ$-torsor $\cC(\tilde L, \tilde N, \tilde M)^2$, where $(\tilde L, \tilde N, \tilde M)\in U_3$.\QED
\end{Lm} 

 For $i=1,2$ let $\gamma_i: \hat Y\times_{\cL(\tilde V)} \hat Y\to \hat Y$ denote the projection to the $i$-th term. To finish the proof of Theorem~\ref{Th_CIO}, it remains to establish the convolution property of $F$. We actually prove it in the following form.
 
\begin{Pp} There is a canonical isomorphism over $(\hat Y\times_{\cL(\tilde V)} \hat Y)\times H$
\begin{equation}
\label{iso_F_conv_on_Y_bars}
(\gamma_1^*F)\ast(\gamma_2^*F)\,\iso\,\gamma^*F
\end{equation}
\end{Pp}
\begin{Prf}
\Step 1 Let 
$
(\hat Y\times_{\cL(\tilde V)} \hat Y)_0\subset (\hat Y\times_{\cL(\tilde V)} \hat Y)
$ 
be the open substack obtained by the base change $U_3\subset Y\times_{\cL(\tilde V)} Y$. In view of Lemma~\ref{Lm_correcting_torsor}, the isomorphism of Corollary~\ref{Cor_convolution_with_theta_multiple} descends under the covering 
$$
\hat\nu_{12}\times \hat\nu_{23}: U_3\to (\hat Y\times_{\cL(\tilde V)} \hat Y)_0
$$ 
to the desired isomorphism (\ref{iso_F_conv_on_Y_bars}) over $(\hat Y\times_{\cL(\tilde V)} \hat Y)_0\times H$.

\smallskip

\Step 2 It suffices to show that $(\gamma_1^*F)\ast(\gamma_2^*F)$ is perverse, the intermediate extension under the open immersion
$$
(\hat Y\times_{\cL(\tilde V)} \hat Y)_0\times H\subset (\hat Y\times_{\cL(\tilde V)} \hat Y)\times H
$$ 
Let us first explain the idea informally, at the level of `functions'. For $(\hat L: \hat M)\in\hat Y$ write $F_{\hat L:\hat M}$ for the restriction of $F$ to this point, this is a `function' on $H$.
For $(\tilde L, \tilde M)\in Y_0$ write $\FF_{\tilde L, \tilde M}$ for the restriction of $F_0$ to this point, this is also a `function' on $H$ (as in the proof of Lemma~\ref{Lm_first_true_convolution}).

 Let $((\hat L: \hat N), (\hat N: \hat M))\in \hat Y\times_{\cL(\tilde V)} \hat Y$. Pick any $\tilde S, \tilde T\in \cL(\tilde V)$ such that 
\begin{itemize}
\item $(\tilde L, \tilde S, \tilde N), (\tilde N, \tilde T, \tilde M)\in U_{13}$,  
\item $\hat\nu_{13}(\tilde L, \tilde S, \tilde N)=(\hat L: \hat N)$, 
$\hat\nu_{13}(\tilde N, \tilde T, \tilde M)=(\hat N: \hat M)$,
\item $\tilde S\cap \tilde T=\tilde S\cap\tilde M=0$.
\end{itemize}
By (\ref{iso_F_explicit_formula}), we get (up to some explicit volumes that we omit)
\begin{multline*}
F_{\hat L:\hat N}\ast F_{\hat N:\hat M}=(\FF_{\tilde L, \tilde S}\ast \FF_{\tilde S, \tilde N})\ast (\FF_{\tilde N, \tilde T}\ast \FF_{\tilde T, \tilde M})=C(\tilde S, \tilde N, \tilde T) \FF_{\tilde L, \tilde S}\ast \FF_{\tilde S, \tilde T}\ast \FF_{\tilde T, \tilde M}=\\
C(\tilde S, \tilde N, \tilde T) C(\tilde S, \tilde T, \tilde M) \FF_{\tilde L, \tilde S}\ast\FF_{\tilde S, \tilde M}=C(\tilde S, \tilde N, \tilde T) C(\tilde S, \tilde T, \tilde M) F_{\hat L: \hat M},
\end{multline*}
where $(\hat L: \hat M)=\nu_{13}(\tilde L, \tilde S, \tilde M)$. Now we turn back to the geometric setting. 

\smallskip

\Step 3 Consider the scheme $\cX$ classifying $(\tilde L, \tilde S, \tilde N)\in U_{13}$, $(\tilde N, \tilde T, \tilde M)\in U_{13}$ such that $\tilde S\cap \tilde T=\tilde S\cap\tilde M=0$. Let 
$$
\zeta: \cX\to \hat Y\times_{\cL(\tilde V)} \hat Y
$$ 
be the map sending the above collection to $(\hat\nu_{13}(\tilde L, \tilde S, \tilde N), \hat\nu_{13}(\tilde N, \tilde T, \tilde M))$. It is smooth and surjective. It suffices to show that 
$$
\zeta^*((\gamma_1^*F)\ast(\gamma_2^*F))
$$ 
is a shifted perverse sheaf on $\cX\times H$, the intermediate extension from $\zeta^{-1}((\hat Y\times_{\cL(\tilde V)} \hat Y)_0)\times H$.

 Let $\mu: \cX\to U_{13}$ be the map sending a point of $\cX$ to $(\tilde L, \tilde S, \tilde M)$. Applying (\ref{iso_F_explicit_formula}) several times as in Step 2, we learn that there is a rank one local system $\cJ$ on $\cX$ such that 
$$
\zeta^*((\gamma_1^*F)\ast(\gamma_2^*F))\,\iso\, \cJ\otimes \mu^*\hat\nu_{13} ^*F
$$
over $\cX\times H$. Since $F$ is an irreducible perverse sheaf, our assertion follows.
\end{Prf}

\medskip

 Thus, Theorem~\ref{Th_CIO} is proved.
 
\medskip\noindent
7.1.6 Given $k$-points $\hat L, \hat M\in \hat\cL(\tilde V)$, let $F_{\hat L, \hat M}$ be the $*$-restriction of $F$ under $(\hat L, \hat M)\times\id: H\to \hat\cL(\tilde V)\times \hat\cL(\tilde V)\times H$. Define the functor $\cF_{\hat L, \hat M}: \D\cH_{\tilde M}\to \D\cH_{\tilde L}$ by
$$
\cF_{\hat L, \hat M}(K)=F_{\hat L, \hat M}\ast K
$$
By Theorem~\ref{Th_CIO}, for $\hat L, \hat N,\hat M\in\hat\cL(\tilde V)$ the diagram is canonically 2-commutative
$$
\begin{array}{ccc}
\D\!\cH_{\tilde M} & \toup{\cF_{\hat N, \hat M}} & 
\;\D\!\cH_{\tilde N}\\
& \searrow\lefteqn{\scriptstyle \cF_{\hat L, \hat M}} &  \;\downarrow\lefteqn{\scriptstyle \cF_{\hat L, \hat N}}\\
&& \;\D\!\cH_{\tilde L}
\end{array}
$$
To see that $\cF_{\hat L, \hat M}$ preserves perversity, pick $\tilde N\in \cL(\tilde V)$ such that $\tilde N\cap\tilde M=\tilde N\cap \tilde L=0$ and use the commutativity of the latter diagram.
This reduces the question to the case $\tilde L\cap\tilde M=0$, in the latter case $\cF_{\hat L, \hat M}$ is nothing but the Fourier transform between the dual vector spaces $L, M$ for the perfect pairing $\omega: L\times M\to 2R$. Here $L=\tilde L\otimes_R k$, $M=\tilde M\otimes_R k$. This completes the definition of (\ref{functor_cF_CIO}).
 
\medskip\noindent 
7.1.7 For a $k$-point $\hat M\in \hat\cL(\tilde V)$ let $i_{\hat M}: \hat\cL(\tilde V)\to \hat\cL(\tilde V)\times \hat\cL(\tilde V)\times H$ be the map sending $\hat L$ to $(\hat L, \hat M, 0)$. Let
\begin{equation}
\label{functor_cF_hat_M_def}
\cF_{\hat M}: \D\!\cH_{\tilde M}\to \D(\hat\cL(\tilde V))
\end{equation}
be the functor sending a complex $K$ to
$
i_{\hat M}^*(F\ast \pr_3^*K)[\dim\cL(\tilde V)-\dim H]
$.
For any $k$-points $\hat M, \hat N\in\hat\cL(\tilde V)$ the diagram commutes
$$
\begin{array}{ccc}
\D\!\cH_{\tilde M} & \toup{\cF_{\hat M}} & \;\D(\hat\cL(\tilde V))\\
 & \searrow\lefteqn{\scriptstyle \cF_{\hat N, \hat M}} & \;\uparrow\lefteqn{\scriptstyle \cF_{\hat N}}\\
 && \;\D\!\cH_{\tilde N}
\end{array}
$$ 
 
\begin{Lm}  The functor $\cF_{\hat N}$ is exact for the perverse t-structures.
\end{Lm}
\begin{Prf} 
Let $\tilde V\,\iso\, \tilde M\oplus \tilde M^*$ be 
a decomposition of $\tilde V$ into an orthogonal sum of two free lagrangian submodules. Remind the open subscheme $\cU(\tilde M)=\{\tilde L\in \cL(\tilde V)\mid \tilde L\cap\tilde M=0\}\subset \cL(\tilde V)$, it identifies naturally with $\Sym^{!2}(\tilde M)$. It suffices to show that for any above decomposition
the composition 
\begin{equation}
\label{functor_comp_for_exactness_cF_hat}
\cH_{\tilde M}\toup{\cF_{\hat M}} \D(\hat\cL(\tilde V))\to \D(\cU(\tilde M))
\end{equation}
is exact, where the second arrow is the restriction under the canonical section $\cU(\tilde M)\to \hat\cL(\tilde V)$.

 Let $M=\tilde M\otimes_R k$. The functor $\cH_{\tilde M}\to \P(M^*)$ sending $K$ to $\bar K=(\tau_{\tilde M^*})^*K[-n-2]$ is an equivalence. Remind the map $\bar\pi: M^*\to \Sym^{*2}(\tilde M^*)$ from Section~4.1. One checks that (\ref{functor_comp_for_exactness_cF_hat}) sends $\bar K$ to $\Four_{\psi}(\bar\pi_!\bar K)\in \P(\Sym^{!2}(\tilde M))$. We are done.
\end{Prf}

\begin{Def} \select{The non-ramified Weil category} $W(\hat\cL(\tilde V))$ is the essential image of the functor (\ref{functor_cF_hat_M_def}). This is a full subcategory in $\P(\hat\cL(\tilde V))$ independent of the choice of a $k$-point $\hat M\in \hat\cL(\tilde V)$.
\end{Def}

 The group $\hat G$ acts naturally on $\hat\cL(\tilde V)$, hence also on $\P(\hat\cL(\tilde V))$. This action preserves the full subcategory $W(\hat\cL(\tilde V))$.

\medskip\noindent
7.1.8 Let $\pr: \bar C\to\hat\cL(\tilde V)$ be the vector bundle whose fibre at $\hat L$ is $L$, where $L=\tilde L\otimes_R k$. Let $\act_l: \bar C\times H\to \hat\cL(\tilde V)\times H$ be the map sending $(\hat L, h, x\in L)$ to $(\hat L, \tau_{\tilde L}(x)h)$. A perverse sheaf $K\in\P(\hat\cL(\tilde V)\times H)$ is \select{$\act_l$-equivariant} if it is equipped with an isomorphism
$$
\act_l^*K\,\iso\, \pr^*K
$$
satisfying the usual associativity property, and whose restriction to the unit section is the identity. 

\begin{Def} \select{The Weil category $W(\tilde V)$} is the category of pairs $(K,\sigma)$, where $K\in \P(\hat\cL(\tilde V)\times H)$ is $\act_l$-equivariant and $(Z(H), \cL_{\psi})$-equivariant, and 
$$
\sigma: F\ast \pr_{23}^*K\,\iso\, \pr_{13}^*K
$$
is an isomorphism. Here $\pr_{13}, \pr_{23}: \hat\cL(\tilde V)\times \hat\cL(\tilde V)\times H\to \hat\cL(\tilde V)\times H$ are the corresponding projections. It is required that $\sigma$ is compatible with the associativity constraint and the unit section constraint of $F$.
\end{Def}
 
  The group $\hat G$ acts on $\hat\cL(\tilde V)\times H$ sending $(g\in\hat G, \hat L, h)$ to $(g\hat L, gh)$. This action extends to an action of $\hat G$ on the category $W(\tilde V)$. One has a natural functor $W(\tilde V)\to W(\hat\cL(\tilde V))$, we don't know if this is an equivalence.
   
\bigskip\smallskip

\centerline{\scshape Appendix  A. Generalities on $H$-gerbs}

\bigskip
\noindent
A.1 The notion of a group stack on an arbitrary site is essentially given in \cite{Br}, where it is called a gr-champ. Here is its algebro-geometric version. Let $S$ be a scheme, $\cG$ be an algebraic $S$-stack. Call it a \select{group stack over $S$} if we are given an action map $\nu: \cG\times_S\cG\to\cG$, a unit section $i: S\to\cG$ over $S$, and the associativity 2-morphism $\alpha: \nu\comp (\nu\times\id)\to \nu\comp(\id\times\nu)$ making the following diagram 2-commutative
$$
\begin{array}{ccc}
\cG\times_S\cG\times_S\cG & \toup{\nu\times\id} & \cG\times_S \cG\\
\downarrow\lefteqn{\scriptstyle \id\times\nu} && \downarrow\lefteqn{\scriptstyle \nu}\\
\cG\times_S \cG & \toup{\nu} & \cG
\end{array}
$$
The morphism $\alpha$ should satisfy the following pentagone axiom. For a scheme $T$ and $T$-points $C_1,C_2,C_3, C_4\in \cG$ the diagram commutes 
$$
\begin{array}{ccc}
((C_1C_2)C_3)C_4 & \toup{\alpha\times\id} & (C_1(C_2C_3))C_4\\
\downarrow\lefteqn{\scriptstyle \alpha} && \downarrow\lefteqn{\scriptstyle \alpha}\\
(C_1C_2)(C_3C_4) & \toup{\alpha} C_1(C_2(C_3C_4)) \getsup{\id\times\alpha} & C_1((C_2C_3)C_4)
\end{array}
$$
We should also be given 2-morphisms $\tau_l: \nu\comp(i\times\id)\to\id$ and $\tau_r: \nu\comp(\id\times i)\to\id$ making the following diagrams 2-commutative
$$
\begin{array}{ccc}
\cG & \toup{i\times\id} & \cG\times_S \cG\\
& \searrow\lefteqn{\scriptstyle \id} & \downarrow\lefteqn{\scriptstyle \nu}\\
&& \cG
\end{array}
\;\;\;\;\;\;\;\;\;\;\;
\begin{array}{ccc}
\cG & \toup{\id\times i} & \cG\times_S \cG\\
&  \searrow\lefteqn{\scriptstyle \id} & \downarrow\lefteqn{\scriptstyle \nu}\\
&& \cG
\end{array}
$$
Writing $O$ for the unit object of $\cG$, the restrictions $\tau_l,\tau_r: OO\to O$ should coincide. Further, the morphisms $\alpha,\tau_l, \tau_r$ should be compatible, namely, for a scheme $T$ and $T$-points $C_1,C_2\in\cG$ the diagrams commute
$$
\begin{array}{ccc}
C_1(OC_2) & \getsup{\alpha} & (C_1O)C_2\\
\downarrow\lefteqn{\scriptstyle \tau_l} & \swarrow\lefteqn{\scriptstyle \tau_r}\\
C_1C_2
\end{array}
\;\;\;\;\;
\begin{array}{ccc}
O(C_1C_2) & \getsup{\alpha} & (OC_1)C_2\\
\downarrow\lefteqn{\scriptstyle \tau_l} & \swarrow\lefteqn{\scriptstyle \tau_l}\\
C_1C_2,
\end{array}
\;\;\;\;\;\;
\begin{array}{ccc}
C_1(C_2O) & \getsup{\alpha} & (C_1C_2)O\\
\downarrow\lefteqn{\scriptstyle \tau_r} & \swarrow\lefteqn{\scriptstyle \tau_l}\\
C_1C_2,
\end{array}
$$
where $O$ is the unit section of $\cG$.  Finally, we require that for any scheme $T$ and any $T$-point $C\in \cG$, writing $\cG_T$ for the category fibre of $\cG$ over $T$, the functors 
$\cG\times_S T\to \cG\times_S T$, $D\mapsto CD$ and 
$\cG\times_S T\to \cG\times_S T$, $D\mapsto DC$ are isomorphisms of $T$-stacks. (We do not explicitely choose a quasi-inverse under multiplication, as we are interested only in actions of $\cG$). 

\medskip\noindent
A.2.1 The notion of a torsor over a group stack is defined in (\cite{Br}, Definition~6.1). Let us formulate its version for algebraic stacks.

 Let $S$ be a scheme, $\cG$ be a group stack over $S$, write
$\nu: \cG\times_S \cG\to \cG$ for the product morphism, let $i: S\to \cG$ be the unity morphism.

Let $Y$ be an $S$-scheme and $\hat Y\to Y$ be an algebraic stack over $Y$. \select{An action of $\cG$} on $\hat Y$ over $Y$ is a data of an action map $m: \cG\times_S \hat Y\to \hat Y$ over $Y$, and a $2$-morphism $\mu: m\comp (\nu\times\id)\to m\comp (\id\times m)$ making the following diagram 2-commutative
\begin{equation}
\label{diag_for_action_1}
\begin{array}{ccc}
\cG\times_S \cG\times_S \hat Y & \toup{\nu\times\id} & \cG\times_S \hat Y\\
\downarrow\lefteqn{\scriptstyle \id\times m} && \downarrow\lefteqn{\scriptstyle m}\\
\cG\times_S \hat Y & \toup{m} & \hat Y
\end{array}
\end{equation}
The 2-morphism $\mu$ should satisfy the pentagone axiom. Namely, for any test scheme $T$ and $T$-points $C_1,C_2,C_3\in \cG$, $D\in \hat Y$ the diagram commutes
$$
\begin{array}{ccc}
((C_1C_2)C_3)D & \to & (C_1(C_2C_3))D\\
\downarrow\lefteqn{\scriptstyle\mu} && \downarrow\lefteqn{\scriptstyle\mu}\\
(C_1C_2)(C_3D) & \toup{\mu} C_1(C_2(C_3D)) \getsup{\mu} & C_1((C_2C_3)D),
\end{array}
$$
where the top horizontal arrow is the associativity constraint for the group stack $\cG$. We further should be given a 2-morphism $\lambda: m\comp (i\times\id)\,\iso\,\id$ making the following diagram 2-commutative
$$
\begin{array}{ccc}
\hat Y & \toup{i\times\id}&  \cG\times_S\hat Y\\
& \searrow\lefteqn{\scriptstyle \id} & \downarrow\lefteqn{\scriptstyle m}\\
&& \hat Y
\end{array}
$$
The morphisms $\mu$ and $\lambda$ should be compatible, namely,  for a test scheme $T$ and $T$-points $C\in \cG, D\in \hat Y$ the diagrams commute
$$
\begin{array}{ccc}
C(OD) & \getsup{\mu} & (CO)D\\
\downarrow\lefteqn{\scriptstyle \lambda} & \swarrow\lefteqn{\scriptstyle \tau_r}\\
CD
\end{array}
\;\;\;\;\;\;\;\;\;\;\;
\begin{array}{ccc}
O(CD) & \getsup{\mu} & (OC)D\\
\downarrow\lefteqn{\scriptstyle \lambda} & \swarrow\lefteqn{\scriptstyle \tau_l}\\
CD,
\end{array}
$$
where $O$ is the unit object of $\cG$, and $\tau_l, \tau_r$ are parts of data for the group stack $\cG$. 

 \select{A $\cG$-torsor} over an $S$-scheme $Y$ is a data of an algebraic stack $\hat Y\to Y$ over $Y$, an action of $\cG$ on $\hat Y$ over $Y$ such that two additional conditions hold. First, 
$m\times\pr_2: \cG\times_S \hat Y\to \hat Y\times_Y \hat Y$ should be an isomorphism. Second, after localization in fppf topology in $Y$ there should exist an isomorphism $\hat Y\,\iso\, Y\times_S \cG$ such that the action map becomes isomorphic to the left translations on $\cG$.

 Now let $\hat Y$ and $\hat Y'$ be two $\cG$-torsors over $Y$. A morphism of such $\cG$-torsors is a pair $(f,h)$, where $f: \hat Y\to \hat Y'$ is a morphism over $Y$, and $h: m\comp (\id\times f)\to f\comp m$ is a 2-morphism making the following diagram 2-commutative
$$
\begin{array}{ccc}
\cG\times_S \hat Y & \toup{\id\times f} & \cG\times_S \hat Y'\\
\downarrow\lefteqn{\scriptstyle m} && \downarrow\lefteqn{\scriptstyle m}\\
\hat Y & \toup{f} & \hat Y'
\end{array}
$$
Besides, $\mu,\mu'$ and $h$ should be compatible, namely, for a test scheme $T$ and $T$-points $C_1, C_2\in \cG$, $D\in \hat Y$ the diagram commutes
$$
\begin{array}{ccc}
(C_1C_2)f(D) & \toup{\mu'} C_1(C_2f(D))  \toup{h} & C_1(f(C_2D))\\
\downarrow\lefteqn{\scriptstyle h} && \downarrow\lefteqn{\scriptstyle h} \\
f((C_1C_2)D) & \toup{f(\mu)} & f(C_1(C_2D))
\end{array}
$$

 Given two morphisms $(f_1,h_1)$ and $(f_2,h_2)$ from $\hat Y$ to $\hat Y'$, a natural transformation from $(f_1,h_1)$ to $(f_2,h_2)$ is a 2-morphism $\phi: f_1\to f_2$ such that $\phi$ is compatible with the actions of $\cG$ on $\hat Y$ and $\hat Y'$ , namely, for a test scheme $T$ and $T$-points $C\in \cG, D\in\hat Y$ the diagram commutes
$$
\begin{array}{ccc}
Cf_1(D) & \toup{h_1} & f_1(CD)\\
\downarrow\lefteqn{\scriptstyle \phi} && \downarrow\lefteqn{\scriptstyle \phi}\\
Cf_2(D) & \toup{h_2} & f_2(CD)
\end{array}
$$ 
  
\medskip\noindent
A.2.2 Now assume that $H$ is a commutative group scheme over $S$ and $\cG$ is the $S$-stack $B(H)$, the classifying stack. Remind that for a morphism of schemes $S'\to S$, the $S'$-points of $B(H)$ is the category of $H\times_S S'$-torsors on $S'$. Then $B(H)$ is a commutative group stack over $S$.

 For this particular $\cG$ the notion of a $\cG$-torsor becomes nothing but a $H$-gerb over $Y$. The definition simplifies considerably as follows (\cite{Br2}, Definition~2.9, p. 49). An \select{$H$-gerb\footnote{in loc.cit. it is called an abelian $H$-gerb} over $Y$} is a stack $\hat Y\to Y$ together with a 2-morphism $\rho: \pr_2\to \pr_2$, where
$\pr_2: H\times_S \hat Y\to \hat Y$ is the projection. It is subject to the following condition. For a scheme $T$ and a $T$-point $D\in \hat Y$ write $\Aut D\mid_T$ for the sheaf of groups on $T$ of automorphisms of $D$. It is required that for any $S$-scheme $T$ the map $\rho: H\mid_T\to \Aut D\mid_T$ is an isomorphism of sheaves of groups on $T$ (in fppf topology). 

 Given an $H$-gerb $\hat Y\to Y$, a $T$-point $D\in \hat Y$ yields an isomorphism $\hat Y\times_S T\,\iso\, B(H)\times_S T$, then the action map $m: B(H)\times_S (\hat Y\times_S T)\to \hat Y\times_S T$ becomes the morphism
$$
B(H)\times_S B(H)\times_S T\to B(H)\times_S T
$$
sending $(\cF_1, \cF_2, t)$ to $(\cF_1\otimes \cF_2, t)$. Here $\cF_i$ are $H$-torsors over $T$.

 Given two $H$-gerbs $\hat Y$ and $\hat Y'$ over $Y$, a morphism of $H$-gerbs is a 1-morphism $f: \hat Y\to\hat Y'$ such that for any scheme $T$, any $T$-points $D\in \hat Y$, $\sigma\in H$ we have $f(\rho_{\hat Y}(\sigma))=\rho_{\hat Y'}(\sigma)$. We strengthen that there is no need to provide in addition a 2-morphism $h$ as in A.2.1 (it is constructed uniquely out of the other data).
 
\medskip\noindent
A.3 Assume in addition that $G$ is a group scheme over $S$, and $G$ acts on an $S$-scheme $Y$ over $S$. The $H$-gerb $\hat Y\to Y$ then gives rise to a group stack $\hat G$ over $G$ defined as follows. For an $S$-scheme $T$ the $T$-points of $\hat G$ is the category of pairs $(g, f)$, where $g\in G(T)$ and 
$$
f: g^*(\hat Y\times_S T)\,\iso\, \hat Y\times_S T
$$ 
is an isomorphism of $H$-gerbs over $Y\times_S T$. A morphism from $(g_1,f_1)$ to $(g_2, f_2)$ exists only under the condition $g_1=g_2$ and it is a natural transformation from $f_1$ to $f_2$. 

   The product morphism $\hat G\times_S\hat G\to\hat G$ sends a $T$-point $(g_1,f_1), (g_2,f_2)$ to the $T$-point $(g_1g_2,f)$, where $f$ is the composition
$$
g_2^*g_1^*(\hat Y\times_S T)\toup{g_2^*f_1} g_2^*(\hat Y\times_S T)\toup{f_2} \hat Y\times_S T,
$$   
Assume that the element of $\H^2(Y, H)$ corresponding to $\hat Y$ is stable under $G$, so that $\hat G\to G$ is surjective. Assume also the following condition:
\begin{itemize}
\item[($\ast$)] For any geometric point $s\in S$ any $H$-torsor on $Y\times_S s$ is trivial
\end{itemize}
Then $\hat G$ is algebraic and fits into an exact sequence of group stacks over $k$
$$
1\to B(H)\to \hat G\to G\to 1
$$   

 Besides, $\hat G$ acts naturally on $\hat Y$, and the projection $\hat Y\to Y$ is equivariant with respect to the homomorphism $\hat G\to G$. The action map $\hat G\times_S \hat Y\to \hat Y$ sends a $T$-point $(g,f)\in \hat G(T), D\in \hat Y(T)$ to $p(f^{-1}D)$, where $p$ is given by the diagram, whose square is cartesian
$$
\begin{array}{ccccc}
\hat Y\times_S T & \toup{f^{-1}} & g^*(\hat Y\times_S T) & \toup{p} & \hat Y\times_S T\\
 & \searrow & \downarrow && \downarrow\\
  && Y\times_S T & \toup{g} & Y\times_S T
\end{array}
$$  

\medskip\noindent
A.4 Let $\cG$ be an algebraic group stack over $S$. Assume that $\hat Y$ is an algebraic stack over $S$, and $\cG$ acts on $\hat Y$ over $S$. Assume that for any scheme $T$ and $T$-points $g\in \cG, y\in\hat Y$ the natural map $\Aut_{G_T}(g)\to \Aut_{\hat Y_T}(gy)$ is injective. Here $G_T$ is the category fibre of $G$ over $T$, similarly for $\hat Y_T$. According to (\cite{D}, Section~2.4.4), in this case one might define the stack quotient $\hat Y/\cG$, which is a priori a 2-stack, it turns out to be  representable by a 1-stack (we don't claim anything about algebraicity of the latter).

 Assume that $\cG$ is smooth of finite type over $S$. Let $K$ be an $\ell$-adic perverse sheaf on $\hat Y$ (here $\ell$ is invertible on $S$). Say that $K$ is \select{$\cG$-equivariant} if we are given an isomorphism
$$
\xi: m^*K\,\iso\, \pr_2^*K
$$
of shifted perverse sheaves for the maps $m, \pr_2: \cG\times_S \hat Y\to\hat Y$, where $m$ is the action map. It is subject to the following conditions (in the notation of A.2.1). The diagram of morphisms (of shifted perverse sheaves on $\cG\times_S\cG\times_S \hat Y$) induced by (\ref{diag_for_action_1}) should be commutative
$$
\begin{array}{ccccc}
(\nu\times\id)^*m^*K & \toup{(\nu\times\id)^*\xi}  & \pr_3^*K & \iso &
\pr_{23}^*\pr_2^*K\\
 \downarrow\lefteqn{\scriptstyle \mu^*} &&&& \uparrow\lefteqn{\scriptstyle \pr_{23}^*\xi}\\
 (\id\times m)^* m^* K & \toup{(\id\times m)^*\xi} & (\id\times m)^*\pr_2^* K &  \iso & \pr_{23}^*m^* K
\end{array}
$$
This means that for a scheme $T$ and $T$-points $C_1, C_2\in \cG, D\in\hat Y$ the diagram commutes
$$
\begin{array}{ccc}
K_{(C_1C_2)D} &\toup{\xi} & K_D\\
\downarrow\lefteqn{\scriptstyle \mu^*} && \uparrow\lefteqn{\scriptstyle \xi}\\
K_{C_1(C_2D)} & \toup{\xi} & K_{C_1D}
\end{array}
$$
where $K$ with a subscript denotes the corresponding $*$-restriction. 

 We used the following. Given two stacks $\cY, \cZ$ and an $\ell$-adic complex $K$ on $\cY$, the inverse image of $K$ is a functor 
$$
\{\mbox{category of 1-morphisms} \; \cZ\to \cY\}\; \to\; 
\{\mbox{the derived category of} \; \ell-\mbox{adic sheaves on} \;\cZ\}
$$
 
 Further, for any $T$-point $D\in\hat Y$, writing $O$ for the unit section of $\cG$, the diagram of $\ell$-adic complexes on $T$ should commute
$$
\begin{array}{ccc}
K_{OD} & \toup{\lambda^*} & K_D \\
& \searrow\lefteqn{\scriptstyle \xi} & \parallel\\
&& K_D
\end{array}
$$

\bigskip

\centerline{\scshape Appendix B. Comparison with the local field case}

\bigskip\noindent
B.1 In this section we show that the Weil representation considered in Section~2 is obtained by some reduction from the Weil representation over the non archimedian local field of characteristic zero (and residual characteristic two). 

Let $F$ be a finite unramified extension of $\QQ_2$, $\cO\subset F$ be the ring of integers and $k$ be the residue field of $\cO$. Set $R=\cO/4\cO$, so $R$ is the ring of Witt vectors of length two over $k$. Let $q$ be the number of elements of $k$.

 Let $W$ be a $F$-vector space of dimension $2n$ with symplectic form $\<.,.\>: W\times W\to F$. Write $H(W)$ for the Heisenberg group $W\ttimes F$ with operation 
$$
(w_1,z_1)(w_2,z_2)=(w_1+w_2, z_1+z_2+\frac{1}{2}\<w_1,w_2\>)
$$
The symbol $\ttimes$ will refer to the above product.

As above, $l\ne 2$. Fix an additive character $\chi_0: \QQ_2\to \Qlb^*$ whose conductor is $\ZZ_2$. Assume that the restriction $\chi_0: \frac{1}{4}\ZZ_2/\ZZ_2\to \Qlb^*$ is the character $\psi: \ZZ/4\ZZ\to\Qlb^*$ we fixed in Section~1. Let $\chi: F\to\Qlb^*$ be given by $\chi(z)=\chi_0(\tr z)$, here $\tr: F\to \QQ_2$ is the trace.

 A subgroup in $A\subset W$ is closed iff $\ZZ_2 A\subset A$ (cf. \cite{MVW}, p.32). For a closed subgroup $A\subset W$ let $A^{\perp}=\{w\in W\mid \<w,a\>\in \cO\;\mbox{for all}\; a\in A\}$. Say that a $\cO$-lattice $M\subset W$ is a \select{symplectic lattice} iff $M^{\perp}=M$. 
 
 For a closed subgroup $A\subset W$ with $A^{\perp}=A$ let $\bar A=A\ttimes F\subset H(W)$, this is a subgroup. The group $A\ttimes (\frac{1}{2}\cO/\cO)$ is abelian. The whole difficulty comes from the fact that there are no natural way to extend the character $\chi$ from $F$ to $\bar A$.
 
 Let $M\subset W$ be a symplectic lattice. Consider the group $(M/2M)\ttimes (\frac{1}{2}\cO/\cO)$, it is abelian (and naturally has a structure of a commutative unipotent algebraic group over $k$). Let $\phi: M/2M\to \frac{1}{2}\cO/\cO$ be a quadratic form such that 
$$
\phi(m_1+m_2)-\phi(m_1)-\phi(m_2)=\frac{1}{2}\<m_1,m_2\>
$$
for all $m_i\in M/2M$. It yields a splitting of the exact sequence of abelian groups
$$
0\to \frac{1}{2}\cO/\cO\to (M/2M)\ttimes (\frac{1}{2}\cO/\cO)\to M/2M\to 0
$$
given by $m\mapsto (m, \phi(m))$.
Once such $\phi$ is chosen,  we get a unique character $\chi_{\phi}: \bar M\to\Qlb^*$ extending $\chi$. Namely, for $(m,z)\in M\ttimes F$ we set $\chi_{\phi}(m,z)=\chi(\phi(m)+z)$.

 Then we get a model of the Weil representation of $H(W)$
\begin{multline*}
\cH_{M,\phi}=\{f: H(W)\to\Qlb\mid f(w h)=\chi_{\phi}(w)f(h), w\in \bar M;\\ 
\mbox{there is an open subgroup}\; M_1\subset W \;
\mbox{such that}\; f(h(w,0))=f(h), w\in M_1, h\in H(W)\}
\end{multline*}

 Let $G(F)=\Sp(W)(F)$ and $G(\cO)$ be the stabilizor of $M$ in $G(F)$. The group $G(F)$ acts on $H(W)$ by automorphisms, namely $g\in G(F)$ sends $(w,z)$ to $(gw,z)$. Write $\rho: H(W)\to \Aut(\cH_{M,\phi})$ for the action by right translations. The model $\cH_{M,\phi}$ yields the metaplectic extension 
\begin{equation}
\label{metapl_ext_for_Q2} 
1\to \Qlb^*\to \tilde G(F)\to G(F)\to 1,
\end{equation}
where
$$
\tilde G(F)=\{(g, \xi)\mid g\in G(F), \xi\in \Aut(\cH_{M,\phi});\;
\rho(gh)\comp \xi=\xi\comp \rho(h), \;\mbox{for}\; h\in H(W)\}
$$

 Let $G_{\phi}(\cO)$ be the group of those $g\in G(\cO)$ whose image in $\Sp(M/2M)$ preserves $\phi$. Then for all $g\in G_{\phi}(\cO), h\in\bar M$ we have $\chi_{\phi}(gh)=\chi_{\phi}(h)$. It follows that $G_{\phi}(\cO)$ acts on $\cH_{M,\phi}$, namely $g\in G_{\phi}(\cO)$ sends $f\in \cH_{M,\phi}$ to $gf$ given by $(gf)(h)=f(g^{-1}h)$, $h\in H(W)$. This is a splitting of (\ref{metapl_ext_for_Q2}) over $G_{\phi}(\cO)$.

 The space $\cH_{M,\phi}$ has a distinguished vector $v_M$, which is extension by zero under $\bar M\hook{} H(W)$ of the function $\chi_{\phi}: \bar M\to\Qlb$. It also has a distinguished linear functional $e_M: \cH_{M,\phi}\to\Qlb$ sending $f$ to $f(0)$.
 
\begin{Rem} For any $\phi$ as above there is a lattice $2M\subset N\subset M$ such that $N/2M\subset M/2M$ is a lagrangian, and $\phi$ vanishes identically on $N/2M$.  
\end{Rem}
 
\medskip\noindent 
B.2 If $M_1\subset W$ is a $\cO$-lattice such that $M_1\subset M_1^{\perp}$ then we get the induced symplectic form $\<.,.\>$ on $M_1^{\perp}/M_1$ with values in $F/\cO$. If moreover, $M_1\subset M\subset M_1^{\perp}$ then $M/M_1$ is a lagrangian in $M_1^{\perp}/M_1$.
 
\begin{Lm} Let $M_1\subset M$ be a $\cO$-lattice. The space $\cH_{M,\phi}^{M_1}:=\{f\in \cH_{M,\phi}\mid f(h(w,0))=f(h), \,\mbox{for all}\; w\in M_1\}$ is as follows. \\
1) If $\phi$ vanishes on $M_1/(M_1+2M)\subset M/2M$ then  $\cH_{M,\phi}^{M_1}$ identifies (via extension by zero) with the space
$$
\{f: M_1^{\perp}\ttimes F\to\Qlb\mid f(wh)=\chi_{\phi}(w)f(h), w\in \bar M\}
$$
The latter space is 
$$
\{f: M_1^{\perp}\to\Qlb\mid f(w+y)=\chi(\phi(w)+\frac{1}{2}\<y,w\>)f(y),\;\mbox{for all}\; w\in M, y\in M_1^{\perp}\}
$$ 
2) We have $\cH_{M,\phi}^M=0$.
\end{Lm}
\begin{Prf}
1) Assume that $f\in \cH_{M,\phi}^{M_1}$ does not vanish on $(w,0)\in H(W)$. Then for any $m_1\in M_1$ we have 
$f(w,0)=f((w,0)(m_1,0))=f((m_1, \<w, m_1\>)(w,0))=\chi(\phi(m_1)+\<w,m_1\>)f((w,0))$. So, for any $m_1\in M_1$ we have $\chi(\phi(m_1)+\<w,m_1\>)=1$.\\
2) is easy. 
\end{Prf} 

 \medskip
 
 Let $2M\subset N\subset M$ be a $\cO$-lattice such that $N/2M\subset M/2M$ is lagrangian, and $\phi$ vanishes on $N/2M$. Note that $N^{\perp}=\frac{1}{2}N$ and 
$\dim_{\Qlb} \cH^{N}_{M,\phi}=q^n$. The action of $N^{\perp}\ttimes F$ on $\cH_{M,\phi}$ preserves the subspace $\cH^N_{M,\phi}$. The symplectic form $\<.,.\>: N^{\perp}\times N^{\perp}\to \frac{1}{2}\cO$ is non degenerate, and 
$N^{\perp}\ttimes \frac{1}{4}\cO\subset N^{\perp}\ttimes F$ is a subgroup. Besides, 
$$
N\ttimes\cO\subset  N^{\perp}\ttimes \frac{1}{4}\cO
$$ 
is a normal subgroup, and the quotient will be denoted $H(N^{\perp}/N)$, it acts naturally on $\cH^{N}_{M,\phi}$. We have an exact sequence of groups
$$
1\to \frac{1}{4}\cO/\cO\to H(N^{\perp}/N)\to N^{\perp}/N\to 1
$$
Its push-forward under $\frac{1}{4}\cO/\cO\to \frac{1}{4}\cO/\frac{1}{2}\cO$ admits a splitting.
 Let $\bar M_N$ be the image of $M\ttimes \frac{1}{4}\cO$ in $H(N^{\perp}/N)$. Then $\chi_{\phi}$ yields a character still denoted $\chi_{\phi}: \bar M_N\to\Qlb^*$, and $\cH^N_{M,\phi}$ identifies with the representation of $H(N^{\perp}/N)$ in
\begin{equation}
\label{Weil_repres_H(N^{perp}/N)}
\{f: H(N^{\perp}/N)\to\Qlb\mid f(wh)=\chi_{\phi}(w)f(h),\;\mbox{for all}\, w\in \bar M_N, h\in H(N^{\perp}/N)\}
\end{equation}
acting in the latter space by right translations. In this way $\cH^N_{M,\phi}$ becomes the Schrodinger model of the oscillator representation of $H(N^{\perp}/N)$ for the enhanced lagrangian $M/N\subset N^{\perp}/N$. The enhanced structure $\tau: M/N\to H(N^{\perp}/N)$ is as follows. First, we have a natural surjective homomorphism
$$
\delta: N^{\perp}\ttimes(\frac{1}{4}\cO/\cO)\to H(N^{\perp}/N),
$$
then for $m\in M$ we have $\tau(m\mod N)=\delta(m,-\phi(m))$.

 Now let $G_N\subset G(F)$ be the stabilizor of $N$. Let $\tilde G_N$ be its preimage in $\tilde G(F)$. Since $G_N$ preserves $N\times\{0\}\subset H(W)$, it follows that $\tilde G_N$ acting on $\cH_{M,\phi}$ preserves $\cH^N_{M,\phi}$. The group $G_N$ acts naturally on $N^{\perp}\ttimes \frac{1}{4}\cO$, and the action of $G_N$ on $H(N^{\perp}/N)$ factors through an action of $\Sp(N^{\perp}/2N)$. 
Set $V=N^{\perp}/N$ and $\tilde V=N^{\perp}/2N$.
 
 One gets the `finite' metaplectic extension $1\to\Qlb^*\to \tilde G_{\tilde V} \to \Sp(\tilde V)\to 1$, where
$$
\tilde G_{\tilde V}=\{(g, \xi_0)\mid g\in \Sp(\tilde V), \xi_0\in\Aut(\cH^N_{M,\phi});\;
\rho_0(gh)\comp \xi_0=\xi_0\comp \rho_0(h),\;\mbox{for}\;
h\in H(N^{\perp}/N)\}
$$
Here $\rho_0$ is the action of $H(N^{\perp}/N)$ by right translations in (\ref{Weil_repres_H(N^{perp}/N)}), and $\xi_0$ is an automorphism of the corresponding $\Qlb$-vector space.
 
 Let $G_{N,1}$ be the kernel of $G_N\to \Sp(\tilde V)$. Then $G_{N,1}\subset G_{\phi}(\cO)$, and the composition $G_{N,1}\hook{} G_{\phi}(\cO)\hook{} \tilde G(F)$ realizes $G_{N,1}$ as a normal subgroup in $\tilde G_N$. This normal subgroup acts trivially on $H(N^{\perp}/N)$, hence a morphism of exact sequences
$$
\begin{array}{ccccccc}
1\to & \Qlb^* & \to  & \tilde G_{\tilde V}  & \to & \Sp(\tilde V)& \to 1\\
& \uparrow\lefteqn{\scriptstyle\id} && \uparrow && \uparrow\\
1\to & \Qlb^* & \to & \tilde G_N & \to & G_N & \to 1
\end{array}
$$ 
yielding an isomorphism $(\tilde G_N)/G_{N,1}\,\iso\, \tilde G_{\tilde V}$.

\end{document}